\documentclass[11pt,a4paper]{article}

\usepackage{graphicx, graphics, epsfig}
\usepackage{booktabs}
\usepackage{subcaption}
\usepackage{natbib}

\usepackage{tikz, pgfplots}
\usetikzlibrary{plotmarks,spy}
\pgfplotsset{compat=newest}
\usepackage[margin=1in]{geometry}
\usepackage{hyperref}
\hypersetup{colorlinks=false}
\usepackage{enumitem}

\usepackage{nicefrac}       % compact symbols for 1/2, etc.
\usepackage{microtype}      % microtypography
\usepackage{xcolor}         % colors

\usepackage{algorithm}
\usepackage{algorithmic}

\usepackage{amsmath}
\allowdisplaybreaks
\usepackage{amssymb}
\usepackage{mathtools}
\usepackage{amsthm}

\usepackage[utf8]{inputenc} % allow utf-8 input
\usepackage[T1]{fontenc}    % use 8-bit T1 fonts
\usepackage{url}            % simple URL typesetting
\usepackage{amsfonts}       % blackboard math symbols

\usetikzlibrary{decorations.pathreplacing}
\usetikzlibrary{quotes,angles}
\usepgfplotslibrary{groupplots}

\usepackage{algorithm}
\usepackage{algorithmic}

\usepackage{multirow}
\usepackage{threeparttable}

\usepackage{makecell}
\usepackage[capitalize,noabbrev]{cleveref}

\theoremstyle{plain}
\newtheorem{theorem}{Theorem}[section]

\newtheorem{lemma}[theorem]{Lemma}
\newtheorem{corollary}[theorem]{Corollary}

\theoremstyle{definition}
\newtheorem{definition}[theorem]{Definition}

\theoremstyle{remark}
\newtheorem{remark}[theorem]{Remark}

\title{Debiasing Random Oblique Projections for Subsampled OLS and Fast CUR in High Dimensions}

\author{
 Chengmei Niu \\
  Huazhong University of Science and Technology\\
  \texttt{chengmeiniu@hust.edu.cn} 
  \and
  Sachin Garg\\
  University of Michigan - Ann Arbor\\
  \texttt{sachg@umich.edu}
  \and
  Micha{\l} Derezi\'nski\\
  University of Michigan - Ann Arbor\\
   \texttt{derezin@umich.edu}
   \and
  Zhenyu Liao\footnote{Author to whom any correspondence should be addressed: Zhenyu Liao (\href{mailto:zhenyu_liao@hust.edu.cn}{zhenyu\_liao@hust.edu.cn}).}\\
     Huazhong University of Science and Technology\\
      \texttt{zhenyu\_liao@hust.edu.cn}
}

\DeclareMathOperator{\tr}{tr}

\DeclareMathOperator{\rank}{rank}
\DeclareMathOperator{\diag}{diag}

\DeclareMathOperator{\Bias}{Bias}
\DeclareMathOperator{\Var}{Var}
\DeclareMathOperator{\loglog}{loglog}
\DeclareMathOperator*{\argmin}{arg\,min}

\newcommand{\RR}{{\mathbb{R}}}

\newcommand{\EE}{{\mathbb{E}}}

\newcommand{\A}{\mathbf{A}}
\newcommand{\B}{\mathbf{B}}
\newcommand{\C}{\mathbf{C}}
\newcommand{\D}{\mathbf{D}}
\newcommand{\F}{\mathbf{F}}

\newcommand{\V}{\mathbf{V}}
\newcommand{\U}{\mathbf{U}}
\newcommand{\W}{\mathbf{W}}
\newcommand{\X}{\mathbf{X}}

\newcommand{\Q}{\mathbf{Q}}
\newcommand{\R}{\mathbf{R}}
\newcommand{\E}{\mathbf{E}}

\renewcommand{\S}{\mathbf{S}}
\renewcommand{\H}{\mathbf{H}}
\renewcommand{\P}{\mathbf{P}}
\newcommand{\K}{\mathbf{K}}
\newcommand{\TT}{\mathbf{T}}
\newcommand{\Y}{\mathbf{Y}}

\newcommand{\uu}{\mathbf{u}}
\newcommand{\vv}{\mathbf{v}}

\newcommand{\ee}{\mathbf{e}}

\newcommand{\x}{\mathbf{x}}
\newcommand{\y}{\mathbf{y}}
\newcommand{\z}{\mathbf{z}}
\newcommand{\rr}{\mathbf{r}}

\newcommand{\w}{\mathbf{w}}

\newcommand{\I}{\mathbf{I}}

\newcommand{\cc}{\mathbf{c}}

\newcommand{\bbeta}{\boldsymbol{\beta}}

\newcommand{\SRHT}{ \mathrm{SRHT} }

\newcommand{\HD}{ \mathrm{HD} }
\newcommand{\OLS}{\mathrm{OLS}}
\newcommand{\CUR}{\mathrm{CUR}}
\definecolor{RED}{rgb}{0.7,0,0}
\definecolor{BLUE}{rgb}{0,0,0.69}
\definecolor{GREEN}{rgb}{0,0.6,0}
\definecolor{PURPLE}{rgb}{0.69,0,0.8}
\definecolor{BLACK}{rgb}{0,0,0}
\definecolor{TEAL}{rgb}{0, 0.5, 0.5}
\begin{document}
\maketitle

\begin{abstract}
Random sampling is a fundamental tool in modern machine learning and numerical linear algebra for reducing the computational cost of large-scale matrix problems.
Existing analyses, however, rely primarily on subspace embedding guarantees, which do \emph{not} precisely characterize the statistical bias of \emph{nonlinear} random oblique projections induced by sampling, which arises ubiquitously in subsampled least squares and fast low-rank approximation methods.
Because (pseudo)inversion is \emph{nonlinear}, these random oblique projections can be systematically \emph{biased} even when the underlying sketch is unbiased, thereby introducing hidden bias into downstream least squares and low-rank approximation solutions.

In this work, we develop a unified non-asymptotic theory for random oblique projections in high dimensions.
We show that standard random sampling schemes generally induce a systematic statistical \emph{bias} overlooked by classical subspace embedding-style analyses, and we propose a principled debiasing framework to correct it.
We illustrate the power of the theory through two canonical applications. 
For subsampled least squares, we obtain sharp bias--variance characterizations, reveal previously unrecognized statistical \emph{suboptimality} in widely used sampling schemes, and identify when debiasing yields provable improvements. 
For fast CUR decomposition, we develop a debiased approach with improved approximation accuracy. 
Numerical experiments further validate our theoretical findings.
\end{abstract}

\section{Introduction}

Many problems in modern machine learning (ML) and scientific computing involve data matrices that are too large to be stored or processed efficiently. 
Randomized numerical linear algebra (RandNLA) addresses this challenge through random sketching techniques~\cite{drineas2006sampling,drineas2011faster,Drineas2012fast,avron2017faster,roosta2019sub,bollapragada2019exact,Lacotte2022adaptive,derezinski2024recent,halko2011finding,niu2025fundamental}. 
Given a tall matrix $\X \in \RR^{n\times p} $ with $n \gg p $, one constructs a sketch $\tilde{\X} = \S \X \in \RR^{m\times p}$, with $ m \ll n $, to serve as a computationally efficient proxy for $\X$ in downstream tasks such as linear least squares, low-rank approximation, and iterative second-order optimization. 
The sketching matrix $\S \in \RR^{m\times n} $ is typically generated either through \emph{random sampling}, which randomly selects rows of $\X$, or through \emph{random projection}, which forms random linear combinations of its rows. 

In this work, we focus on the random \emph{oblique projection} $\tilde \P$~\cite{chi2021multip,chi2021aprojector} induced by \emph{random sampling}, together with its associated residual projection $\tilde \P_{\perp}$, defined as: $\tilde \P =\X(\S\X)^{\dagger}\S,$ $\tilde \P_{\perp}= \I_n - \tilde \P$.

These operators are designed to approximate, and should thus be contrasted with, the orthogonal projection $\P = \X \X^\dagger$ onto the column space of $\X$, and its residual projection $\P_\perp = \I_n - \P$.
Random oblique projections arise naturally in a broad range of RandNLA methods, including subsampled least squares~\cite{drineas2006sampling,drineas2011faster}, fast low-rank approximation methods such as CUR decomposition~\cite{wang2016towards,ye2019fast}, and randomized optimization methods~\cite{bartan2022distributed}.

A common principle in the design of sampling schemes is to ensure that the subsampled Gram matrix $\tilde{\X}^\top \tilde{\X}$ is an unbiased or nearly unbiased estimator of the full Gram matrix $\X^\top \X$, i.e., $\EE[\tilde{\X}^\top \tilde{\X}] = \X^\top \X$. 
However, such first-order moment matching at the Gram level does \emph{not} generally imply unbiasedness of the induced oblique projection $\tilde \P = \X(\S\X)^{\dagger}\S$.
Because the Moore--Penrose pseudoinverse is nonlinear, one typically has $\EE[ \tilde\P] \not \approx \P$, even when the underling sketch itself if unbiased.
Existing analyses of random oblique projections (and consequently of subsampled least squares and fast low-rank approximation methods) rely primarily on Johnson--Lindenstrauss (JL)-style analyses~\cite{johnson1984extensions}.
Such results aim to establish that $\tilde{\X}^{\top}\tilde{\X} \approx \X^{\top}\X$ with high probability and therefore that $\tilde{\P} \approx \P$ in some (matrix norm) sense. 
While powerful and easy to use, such guarantees are often too coarse to capture finer statistical behaviors central to the practical performance of subsampled least squares and fast low-rank approximation methods, as simple and fundamental as bias and variance.

In this paper, we address this gap by developing a refined non-asymptotic theory for oblique projections induced by random sampling, and illustrate its implication through the two canonical applications of subsampled least squares and fast CUR decomposition.

\subsection{Our contributions} 

Our main contributions are summarized as follows.

\begin{enumerate}[leftmargin=*, itemsep=1.2pt, topsep=-0.5pt, parsep=0pt, partopsep=0pt]%
  \item We provide in \Cref{theo:de_sam_obliqueprojec} a precise non-asymptotic characterization of the bias induced by random oblique projections under general random sampling schemes, along with a principled debiasing approach.
  
  \item Building on this result, we first establish in \Cref{sec:fin_bias_var_ols} sharp bias--variance characterizations for subsampled ordinary least squares (OLS), including a bias lower bound for classical subsampled OLS (\Cref{theo:lower_bound}) and an bias improved upper bound for the proposed \emph{debiased} subsampled OLS, whose variance matches that of the classical solution (\Cref{theo:de_sam_ols} versus \Cref{theo:variance_lev_least-sq}).
  We further extend this analysis to randomized CUR decomposition (\Cref{theo:cur}), deriving a debiased fast CUR method with provably improved approximation accuracy. 
  
\end{enumerate}

\subsection{Related work}

Our work connects to three closely related lines of research: subsampled OLS, sketching-based CUR decomposition, and recent analyses based on random matrix theory beyond classical asymptotics. We briefly review these directions below and position our contributions relative to prior work.

\paragraph{Classical subsampled OLS.}
Subsampled OLS estimators are a classical tool for accelerating least-squares problems by reducing the problem size while approximately preserving the information in the original design matrix~\cite{drineas2006sampling,drineas2011faster,wangsketchridge2018,bartan2022distributed}. In addition to their computational benefits, a growing body of work has studied subsampled OLS from a statistical perspective, focusing on properties such as bias and variance~\cite{ma2015statistical,garvesh2016astatistic,derezinski2017unbiased,wangsketchridge2018,dobriban2019asymptotics,derezinski2019correcting,derezinski2019minimax,bartan2022distributed,derezinski2022unbiased}.
Despite this progress, existing analyses do not provide a precise characterization of the bias between the subsampled OLS estimator and the full OLS solution, except for when using specialized determinantal sampling schemes, which are computationally expensive \cite{derezinski2017unbiased,derezinski2019correcting,derezinski2022unbiased}. A central difficulty is the lack of a fine-grained understanding of the oblique projection operator $\tilde{\P}$ induced by random sampling. As a result, prior work typically relies on coarse bounds or asymptotic arguments that do not fully capture the mechanism by which sampling introduces bias.
Our work addresses this gap by leveraging precise characterizations of the oblique projector $\tilde{\P}$ to obtain a more refined analysis of the bias in subsampled OLS estimators.

\paragraph{CUR decomposition.}
Traditionally, low-rank matrix approximations are obtained via SVD-based methods, such as truncated SVD or randomized SVD~\citep{halko2011finding}, which operate in the space of matrix singular vectors. 
In contrast, CUR decomposition constructs low-rank approximations directly from selected rows and columns of the data matrix~\cite{michael2009curmatrix}, thereby preserving structural properties such as sparsity, non-negativity, and interpretability. These properties are often important in applications such as signal and image processing~\citep{candes2008Introduction,elad2006Image}, recommendation systems~\citep{koren2009Matrix}, and big data 
analysis~\citep{michael2009curmatrix}. Given a matrix $\X\in\RR^{n\times p}$, a key design choice in CUR decomposition is the middle factor $\U$. Two common options are
 $\U=\C^{\dagger}\X\R^{\dagger}$, where $\C$ and $\R$ consist of sampled columns and rows of $\X$, and $\U=\X_{R,C}^{\dagger}$, where $\X_{R,C}$ denotes their intersection. The former is more robust but computationally expensive, while the latter is more efficient but can be unstable when $\X_{R,C}$ is nearly singular.   Classical approaches mitigate this instability via pivoting-based methods, such as column-pivoted QR and LU with complete pivoting~\cite{golub2013matrix,trefethen2022numerical}, which enjoy strong theoretical guarantees~\cite{gu1996efficient} but are often impractical for large-scale problems. 

To improve scalability, recent work has proposed randomized sketching methods. For $\U=\X_{R,C}^{\dagger}$, pivoting is typically performed on a randomized sketch $\tilde{\X}$ instead of $\X$~\cite{dong2023simpler,park2025accuracy}, with oversampling shown to improve stability and accuracy~\cite{park2025accuracy}.  For another  choice $\U=\C^{\dagger}\X\R^{\dagger}$, existing approaches rely on the sketches of $\C$, $\R$, and $\X$ to form estimators such as $\tilde{\C}^{\dagger}\tilde{\X}\tilde{\R}^{\dagger}$~\cite{wang2016towards,ye2019fast}. However, existing analyses rely largely on subspace embedding arguments and provide only coarse guarantees, without precisely characterizing how random sketching affects the induced projection operators $\C\tilde{\C}^{\dagger}$ and $\tilde{\R}^{\dagger}\R$. In this work, we develop a fine-grained analysis of these random-sampling-based estimators by precisely characterizing the associated projection operators, leading to sharper guarantees and improved approximations of $\C^{\dagger}\X\R^{\dagger}$.

\paragraph{Random matrix theory.}
Our analysis is closely related to recent advances inspired by asymptotic random matrix theory (RMT)~\cite{couillet2022RMT4ML}, which studies the spectral behavior of large-dimensional random matrices~\cite{anderson2010introduction}. RMT-based techniques have recently been adopted in large-scale machine learning to analyze optimization and generalization properties~\cite{pennington2017nonlinear,fan2020spectra,mei2021generalization,couillet2022RMT4ML}. In particular, building on tools such as Stieltjes transforms, several works have developed fine-grained characterizations of inverse Gram matrices. Using Sherman-Morrison rank-one updates, \cite{derezinski2021newtonless,derezinski2021sparse,garg2024distributed,niu2025fundamental} derive precise expressions for $\EE[(\tilde\X^\top\tilde\X)^{-1}]$ beyond classical asymptotic analyses. Our work is closely connected to this line of research and extends these techniques to obtain a fine-grained characterization of the oblique projection operator $\tilde{\P}$ under random sampling.

\paragraph{Notation.}
Scalars, vectors, and matrices are denoted by lowercase letters, bold lowercase letters, and bold uppercase letters, respectively. 
For a matrix $\X\in \RR^{n\times p}$, we denote $\X^\top$, $\X^\dagger$, $\x^\top_{i}\in \RR^{p}$, $\|\X\|$, and $\|\X\|_F$ the transpose, the Moore--Penrose pseudoinverse, $i^{\text{th}}$ row, the spectral norm, and the Frobenius norm of $\X$, respectively.
We denote $\X \preceq \Y$ if $\Y - \X$ is positive semi-definite, and use $\I_p$ for the identity matrix of size $p$.
For a random variable $x$, $\EE[x]$ denotes its expectation, and $\EE_{\zeta}[x]$ denotes its expectation conditioned on the event $\zeta$. 
We use the standard Big-$O$ and Big-$\Omega$ notations, and $\tilde O(\cdot)$ to ignore logarithmic dependence on $n$.

\section{Preliminaries}\label{sec:pre}

In this section, we introduce several definitions that will be used throughout the paper.

\begin{definition}[\textbf{Random sampling}]\label{def:RS}
For a matrix $\X \in \RR^{n \times p}$ with $n\geq p$, a sketch $\tilde \X \in \RR^{m \times p}$ of $\X$ can be constructed by sampling \emph{with replacement} $m$ of the $n$ rows of $\X$ with a sampling distribution, $\{\pi_i\}_{i=1}^n$, $\sum_{i=1}^n \pi_i = 1$, and then rescaling by $1/\sqrt{m \pi_i}$.
This procedure writes $\tilde\X=\S\X$, with \underline{\emph{sampling matrix}} $\S\in \RR^{m\times n}$ having only one nonzero entry per row.
One has $\EE[\tilde\X^\top \tilde\X] = \X^\top \X$.
\end{definition}

\Cref{def:RS} includes commonly used random sampling schemes such as uniform ($\pi_i = 1/n$), row-norm-based ($\pi_i = \|\x_i\|^2 / \sum_{j=1}^n \|\x_j\|^2$), exact and approximate leverage score sampling~\cite{mahoney2011randomized}, as well as hybrid schemes that interpolate between them, e.g., shrinkage leverage score sampling~\cite{ma2015statistical}.

\begin{definition}[\textbf{Leverage score sampling}]\label{def:lev} 
For $\X \in \RR^{n \times p}$ of rank $p$ with $n \geq p$, the $i^{\text{th}}$ \underline{\emph{leverage score}} $\ell_i(\X)$ of $\X$ is defined as $\ell_i(\X) = \x_i^\top(\X^\top\X)^{-1}\x_i, i \in \{ 1, \ldots, n \}$.
The exact leverage score sampling refers to the random sampling approach in \Cref{def:RS} with $\pi_i=\ell_i(\X)/p$.  
\end{definition}

\begin{definition}[\textbf{Importance sampling approximation factor},~\cite{niu2025fundamental}]\label{def:approx_factor}
For $\X \in \RR^{n \times p}$ with $n \geq p$ and a random sampling matrix $\S \in \RR^{m\times n}$ in \Cref{def:RS} with sampling distribution $\{\pi_i \}_{i=1}^n$, the associated minimum and maximum \underline{\emph{importance sampling approximation factors}} are defined as $\theta_{\min}(\X) \equiv \min_{1\leq i \leq n}\ell_i(\X)/(\pi_i p)$ and $\theta_{\max}(\X) \equiv \max_{1\leq i \leq n}\ell_i(\X)/(\pi_i p)$.
\end{definition}

The importance sampling approximation factors in \Cref{def:approx_factor}, introduced in~\cite{ma2015statistical,niu2025fundamental}, quantify how a random sampling scheme deviates from the \emph{exact} leverage score sampling in \Cref{def:lev}, 
In particular, one has $\theta_{\min}(\X) \leq 1 \leq  \theta_{\max}(\X)$ with $\theta_{\min}(\X) = \theta_{\max}(\X) = 1$ for \emph{exact} leverage score sampling.
As we shall see below, these parameters are crucial in our analyses of (pseudo)inverses.

\begin{definition}[\textbf{Subspace embedding},~\cite{drineas2006sampling,mahoney2011randomized}] \label{def:rela_error_approxi}
For $\X \in \RR^{n \times p}$ with $n\geq p$ and a random sampling matrix $\S\in \RR^{m\times n}$ with $m\leq n$, we say that a sketch $\tilde \X=\S\X \in \RR^{m \times p}$ is an $(\epsilon,\delta)$-\emph{subspace embedding}, or an $(\epsilon,\delta)$-approximation for $\X$ if
\begin{equation*}
  (1+\epsilon)^{-1}\X^\top\X \preceq \tilde \X^\top \tilde \X  \preceq(1+\epsilon)\X^\top\X,
\end{equation*}
holds with probability at least $1-\delta$.
\end{definition}
The subspace embedding property in \Cref{def:rela_error_approxi} ensures that, with high probability, $\tilde \X^\top \tilde \X$ provides a reliable approximation of $\X^\top \X$.
This property has been extensively exploited in the literature to establish statistical guarantees for a wide range of sampling schemes $\S$; see~\cite{mahoney2011randomized,david2014sketching}.

Nonetheless, for many applications in ML and scientific computation, the subspace embedding-type guarantee in \Cref{def:rela_error_approxi} is \emph{not} sufficient.
As an instance, while one has $\EE[\tilde\X^\top \tilde\X] = \X^\top \X$ (as in \Cref{def:RS}), the sketched matrix inverse $(\tilde \X^\top \tilde \X)^{-1}$ is \emph{no longer unbiased}, i.e., $\| \EE[(\tilde \X^\top \tilde \X)^{-1}] - (\X^\top\X)^{-1} \| \gg 0$, due to the nonlinear nature of the inverse.

\begin{remark}[\textbf{Inversion bias for Gaussian random projection}]\label{rem:inversion_bias_Gaussian}
In the case of Gaussian random projection, with $\tilde \X = \S \X$ for $\S$ having i.i.d.\@ Gaussian entries with zero mean and variance $1/m$, the inverse $(\tilde \X^\top \tilde \X)^{-1}$ is known to follow the \emph{inverse Wishart distribution}~\cite{haff1979identity} with $\EE[(\tilde \X^\top \tilde \X)^{-1}] = \frac{m}{m-p-1} (\X^\top \X)^{-1}$.
As such, the scalar debiasing $\frac{m}{m-p-1}$ is \emph{exact} for Gaussian projections. Moreover, a similar \emph{scalar} debiasing factor, $\frac{m}{m-p}$, is known to remain effective for both i.i.d.\@ sub-Gaussian and LESS random projections~\cite{derezinski2021sparse}. 

\end{remark}

The \emph{inversion bias} for random sampling has been characterized in the following result.

\begin{theorem}[\textbf{Inversion bias for random sampling},~{\cite[Theorem~3.1]{niu2025fundamental}}]\label{theo:inverse-bias}
For $\X\in \RR^{n\times p}$ of rank $p$ with $n \geq p$, let $\S \in \RR^{m \times n}$ be a random sampling  matrix with sampling distribution $\{ \pi_i \}_{i=1}^n$ as in \Cref{def:RS} and $\theta_{\min}(\X), \theta_{\max}(\X)$ as in \Cref{def:approx_factor}.
Then, for diagonal matrix $\D = \diag\{ D_{ii} \}_{i=1}^n$ the solution to $D_{ii} = \left(1 + \x_i^\top (\X^\top \D \X)^{-1} \x_i/(m \pi_i) \right)^{-1}$,
there exists $C > 0$ independent of $n,p$, so that for $m\geq C\theta_{\max}(\X) p(\log(p/\delta )+ (\log n/ \loglog n)^{2/3}/\epsilon^{2/3}),~\delta\leq m^{-3}$, when conditioned on an event $\zeta$ that holds with probability at least $1-\delta$,
\begin{equation*}
  (1 + \epsilon)^{-1}(\X^\top \D \X)^{-1}  \preceq  \EE_{\zeta} [ (\tilde \X^\top \tilde \X)^{-1} ] \preceq (1 + \epsilon) (\X^\top \D \X)^{-1}.
\end{equation*}
\end{theorem}

\Cref{theo:inverse-bias} shows that \emph{any} random sampling scheme in \Cref{def:RS} admits an \emph{inversion bias}, and that the subsampled inverse $(\tilde \X^\top \tilde \X)^{-1}$ is close, \emph{in expectation}, to $(\X^\top \D \X)^{-1}$, as opposed to $(\X^\top \X)^{-1}$ that one may expect from the subspace embedding guarantee in \Cref{def:rela_error_approxi}.
It has then been shown in~\cite[Proposition~3.2]{niu2025fundamental} that this bias can be effectively corrected, using a \emph{matrix-level} debiasing approach, leading to a debiased sampling matrix $\check\S$ of the same size, for which one has $\check \X = \check\S \X$ and $\EE[(\check \X^\top \check\X)^{-1}]\simeq (\X^\top\X)^{-1}$ in a spectral norm sense.

Since one has, for $\X \in \RR^{n \times p}$ of rank $p$ with $n \geq p$, that the Moore--Penrose pseudoinverse $\X^\dagger$ satisfies $\X^\dagger = (\X^\top \X)^{-1} \X^\top$. 
One should expect, as a consequence of the inversion bias in \Cref{theo:inverse-bias}, that the random oblique projection $\tilde \P=\X(\S\X)^\dagger\S$ is also a \emph{biased} estimator of the projection matrix $\P = \X \X^\dagger$.
It is thus natural to ask whether random oblique projections can also be effectively debiased.
In the following section, we present such results.

\section{Debiased oblique projection for random sampling}\label{sec:main_results}

In the following result, we show that the debiasing sampling matrix $\check \S$ introduced in~\cite{niu2025fundamental} can also be used to effectively debias the random oblique projection.

\begin{theorem}[\textbf{Precise characterizations of debiased oblique projection}]\label{theo:de_sam_obliqueprojec}
For $\X \in \RR^{n\times p}$ of rank $p$  with $n\geq p$, assume that $ \min_{i,\|\x_i\|>0}\|\x_i\|^2 /\|\X\|_F^2 \geq n^{-\alpha}$ for some constant $\alpha>0$. Let $\S$ be a random sampling matrix with sampling distribution $\{ \pi_i \}_{i=1}^n$ as in \Cref{def:RS}, define the debiased sampling matrix $ \check\S \in \RR^{m\times n}$~as\footnote{
For rank-deficient $\X$ (i.e., $\rank(\X) < p$), it suffices to restrict to the column space of $\X$ (that is of dimension smaller than $p$) and \Cref{theo:de_sam_obliqueprojec} naturally extends to this setting.
}
\begin{equation}\label{eq:debias_check_S}
  \check\S = \diag \left\{1/\sqrt{1-\ell_{i_s}(\X)/(m \pi_{i_s}) } \right\}^m_{s=1}\cdot\S.\hspace{-1mm} %\quad i_s\in \{1,\ldots,n\}.
\end{equation}
Then, there exists $C > 0$ independent of $n, p$ so that for $m \geq C\theta_{\max}(\X) p\log (p/\delta)$, $\delta\leq n^{-(3+\alpha)}$, 
and $\theta_{\max}(\X)$ in \Cref{def:approx_factor}, when conditioned on an event $\zeta$ that holds with probability at least $1-\delta$, the debiased oblique projection $\check \P \equiv \X( \check\S \X)^{\dagger} \check\S $ satisfies\footnote{The event $\zeta$ captures the high-probability guarantee that for  $m \geq C\theta_{\max}(\X) p\log (p/\delta)$,  $\S\X$ is an $(\epsilon,\delta)$-subspace embedding of $\X$ (see \Cref{lem:sub_embed} of \Cref{sec:use_lemms}), ensuring a well-behaved sketch throughout the analysis.
\label{footnote:event_zeta}
} 
\begin{align}
 \left\| \EE_\zeta[\check \P] - \P \right\|_{F}^2 &= \epsilon^2 \cdot \|\P_{\perp}\|_{F}^2,\label{eq:expec_obliquepro_bias}\\
 \EE_{\zeta} \left[\| \check \P - \P \|^{2}_{F} \right] &= \tr \left( \P_{\perp} \diag \left\{ \ell_{i}(\X)/(m \pi_i) \right\}_{i=1}^n \right) + \epsilon \cdot \|\P_{\perp}\|^{2}_{F},\label{eq:expec_obliqueprovariance}
\end{align} 
with $\epsilon=O (\sqrt{(\log n/\loglog n)^2 \cdot \theta_{\max}^3(\X) p^3/m^3} )$, for orthogonal projection  $\P = \X\X^{\dagger}$  and residual projection $\P_{\perp}=\I_n - \X\X^{\dagger}$.
\end{theorem}

Similar to \cite{niu2025fundamental}, in \Cref{theo:de_sam_obliqueprojec} we leverage leave-one-out arguments to ``expand'' Moore--Penrose pseudoinverses under random sampling. 
However, establishing the bound in \eqref{eq:expec_obliqueprovariance} requires controlling $\EE[\check\S^\top(\X^\top \check\S^\top )^{\dagger} \X^\top\X( \check\S \X)^{\dagger} \check\S]$ rather than inverse moments such as $\EE[(\X^\top \check\S^\top \check\S \X)^{-2} ]$ considered in \cite{niu2025fundamental}. 
The former involves cross interactions between distinct sampled rows that cannot be decoupled using a single leave-one-out.
To address this issue, we adopt a refined non-asymptotic leave-two-out approach in the proof of \eqref{eq:expec_obliqueprovariance}. 
The detailed proof is provided in \Cref{sec:proof_theo_de_sam_obliqueprojec}.

\Cref{theo:de_sam_obliqueprojec} provides precise characterizations of both first-~and~second-order moments of the debiased oblique projection.
Precisely, we have that the debiased oblique projection $\check \P$ is close, in expectation, to the true projection $\P$, up to an error of order $\tilde O(\theta_{\max}^3(\X) p^3/m^3 )$ for \emph{all} sampling schemes. 
Moreover, the seconder-order moment $\EE_{\zeta} \left[\| \check \P - \P \|^{2}_{F} \right]$ is also precisely characterized, as the sum of the leading-order term $\tr ( \P_{\perp} \diag \{ \ell_{i}(\X)/(m \pi_i) \}_{i=1}^n )$ and some error, also of order $\tilde O( \sqrt{\theta_{\max}^3(\X) p^3/m^3})$.

Notably, it can be checked that for $\theta_{\min}(\X), \theta_{\max}(\X)$ in \Cref{def:approx_factor} one has $\frac{\theta_{\min}(\X) p}{m}\cdot \I_n \preceq \diag \{ \ell_{i}(\X)/(m \pi_i) \}_{i=1}^n \preceq \frac{\theta_{\max}(\X) p}{m} \cdot \I_n$, so that \eqref{eq:expec_obliqueprovariance} further writes $( \frac{\theta_{\min}(\X) p}{m} + \epsilon) \|\P_{\perp}\|^{2}_{F} \leq\EE_{\zeta} \left[\| \check \P - \P \|^{2}_{F} \right] \leq ( \frac{\theta_{\max}(\X) p}{m} + \epsilon ) \|\P_{\perp}\|^{2}_{F}$, with equality for \emph{exact} leverage score sampling.

\begin{remark}[\textbf{Scalar debiasing for exact leverage score sampling}]\label{rem:scalar_debiasing_exact}
Note that under exact leverage score sampling, where $\pi_i = \ell_i(\X)/p$ for all $i$ (see again \Cref{def:lev}), the debiased sampling matrix in \eqref{eq:debias_check_S} reduces to $\check \S = \sqrt{\frac{m}{m-p}} \S$, which coincides with the \emph{scalar} biasing factor appearing for Gaussian, sub-Gaussian, and LESS projections in \Cref{rem:inversion_bias_Gaussian}.
In this setting, for the debiased oblique projection $\check \P \equiv \X( \check\S \X)^{\dagger} \check\S $, the scalar factor $\frac{m}{m-p}$ inside and outside the pseudoinverse cancels out, yielding $\check \P = \tilde \P$.
Consequently, \emph{no} debiasing is needed for exact leverage score sampling.
One may thus expect that a similar phenomenon persists, at least to some extent, when the sampling scheme is close to exact leverage score sampling.
This statement will be made precise in the next section.
\end{remark}

\section{Application to subsampled OLS and fast CUR} \label{sec:fin_bias_var_ols}

Most existing analyses of subsampled OLS rely on subspace-embedding-type arguments and fail to provide precise characterizations of, say its first-~and~second-order moments due to random sampling~\cite{wangsketchridge2018,ma2015statistical,garvesh2016astatistic,dobriban2019asymptotics,derezinski2019minimax}. 
In this section, we show how the debiased oblique projection in \Cref{theo:de_sam_obliqueprojec} applies to establish sharper characterizations of the bias and variance of subsampled OLS (see \Cref{def:bias_and_variance} for their formal definitions).

Given a data matrix $ \X\in\RR^{n\times p}$ and a response vector $ \y \in\RR^n$, the OLS solution is given by
\begin{equation}\label{eq:def_OLS_L}
    \bbeta_{\OLS}=\argmin_{\bbeta\in \RR^p} L(\bbeta) =\X^\dagger \y,~~\text{for}~~L(\bbeta) = \|\y-\X\bbeta\|^2,
\end{equation}
which can be expensive when $n$ or $p$ is large.
We consider the subsampled OLS~\cite{drineas2006sampling,wangsketchridge2018,bartan2022distributed} solution:
\begin{align}\label{eq:hat_bbeta}
    \tilde \bbeta=(\S\X)^{\dagger}\S \y,
\end{align}
where $\S \in \RR^{m \times n}$ is a random sampling matrix as in \Cref{def:RS}.
We also consider its associated \emph{debiased} counterpart:
\begin{equation}\label{eq:debia_sam_ols}
     \check \bbeta=( \check \S\X)^{\dagger} \check \S \y,
\end{equation}
where $\check\S$ is the debiased sampling matrix in \eqref{eq:debias_check_S} of \Cref{theo:de_sam_obliqueprojec}.

Observe that $L(\tilde\bbeta)=\|\y-\X(\S\X)^\dagger \S \y\|^2$ and $ L(\check\bbeta)=\|\y-\X(\check\S\X)^\dagger \check\S \y\|^2$. 
Thus, subsampling replaces the orthogonal projection $\X\X^\dagger$ with the oblique projections $\X(\S\X)^\dagger \S$ and $\X(\check\S\X)^\dagger \check\S$.
% onto the column space of $\X$.  
The deviation of these oblique projections from orthogonal projections governs the statistical behavior of the corresponding estimators relative to OLS.

In this section, we evaluate the statistical behavior of the classical subsampled OLS solution $\tilde\bbeta$ in \eqref{eq:hat_bbeta}, the debiased solution $\check \bbeta$ in \eqref{eq:debia_sam_ols}, versus that of the OLS solution $\bbeta_{\OLS}$ in \eqref{eq:def_OLS_L}, per the following two metrics, in line with a series of previous efforts~\cite{bartan2022distributed,derezinski2019minimax,derezinski2022unbiased,garg2024distributed,garvesh2016astatical}.

\begin{definition}[\textbf{Bias and variance of subsampled OLS}]\label{def:bias_and_variance}
For a data matrix $ \X\in\RR^{n\times p}$ of rank $p$, a response vector $ \y \in\RR^n$, and the OLS solution $\bbeta_{\OLS}$ defined in \eqref{eq:def_OLS_L}, let $\bbeta \in \RR^{p}$ be a random vector (such as the subsampled OLS solution $\tilde \bbeta$ in \eqref{eq:hat_bbeta} or its debiased counterpart $\check \bbeta$).
Define 
\begin{enumerate}[leftmargin=*, itemsep=1.2pt, topsep=-0.5pt, parsep=0pt, partopsep=0pt]%
%[leftmargin=*, itemsep=1.2pt, topsep=-1pt, parsep=-1pt, partopsep=-1pt]
  \item \underline{\emph{Bias}} $\Bias_{\zeta}(\bbeta) = L(\EE_{\zeta}[\bbeta]) - L(\bbeta_{\OLS})$ as the deviation of $\bbeta$ from $\bbeta_{\OLS}$ in expectation; and
  \item \underline{\emph{Variance}} $\Var_{\zeta}(\bbeta) = \EE_{\zeta}[L(\bbeta)] - L(\bbeta_{\OLS})$ that quantifies the fluctuation due to random sampling; %second-order
\end{enumerate}
both measured by the square loss $L(\cdot)$ in \eqref{eq:def_OLS_L} and conditioned on a high probability event $\zeta$ (which is necessary for many random sampling schemes, see \Cref{footnote:event_zeta}).
\end{definition}

\subsection{Bias and variance characterizations of subsampled OLS with and without debiasing}
\label{sub:bias_and_variance_of_subsampled_OLS}

We first establish in \Cref{theo:lower_bound} that the classical subsampled OLS $\tilde \bbeta$ in \eqref{eq:hat_bbeta}, \emph{without debiasing}, can exhibit a \emph{substantial bias} (in the sense of \Cref{def:bias_and_variance}), at least under \emph{certain} sampling schemes.
We then derive in \Cref{theo:variance_lev_least-sq} an upper bound on the variance of $\tilde \bbeta$.
% (which, interestingly, admits the same dominant and higher-order decomposition as that in \Cref{theo:de_sam_ols}). 
Finally, in \Cref{theo:de_sam_ols}, we present precise bias--variance characterizations for the proposed \emph{debiased} subsampled OLS $\check \bbeta$, demonstrating that debiasing not only mitigates the bias but also preserves, if not \emph{improves}, the variance behavior.

\begin{theorem}[\textbf{Bias lower bound for subsampled OLS}]\label{theo:lower_bound}%\textbf{under approximate leverage score sampling}
Let $p=4k$ for any $k\geq 1$ and let $n\geq 2p$. 
Then, there exists a data matrix $\X \in \RR^{n\times p}$, a response vector $\y \in \RR^n$, an approximate leverage score sampling matrix $\S \in \RR^{m\times n}$ with $\{\pi_i\}_{i=1}^n$ satisfying $ \ell_i(\X)/2 \leq \pi_i \leq 2\ell_i(\X)$ for all $i$, such that for all $m\geq p$ and \emph{any} real constant $\gamma$, when conditioned on an event $\zeta$ that holds with probability at least $1-\delta$, 
\begin{equation*}
    %\Bias_{\zeta}(\gamma \cdot \tilde \bbeta) = \Omega \left(\frac{p^2}{m^2} \right) \cdot\|\rr\|^2,
    \Bias_{\zeta}(\gamma \cdot \tilde \bbeta) = \Omega (p^2/m^2 ) \cdot\|\rr\|^2,
\end{equation*}
with $\Bias_{\zeta}(\cdot)$ in \Cref{def:bias_and_variance}, $\tilde\bbeta = (\S\X)^\dagger \S \y$ in \eqref{eq:hat_bbeta}, and residual vector $\rr = \y-\X\beta_{\OLS}$.
\end{theorem}

\Cref{theo:lower_bound} says that, under certain random sampling scheme (that is $1/2$-approximate leverage score sampling), there exists scenarios for which the bias of classical subsampled OLS $\tilde \bbeta$ is \emph{at least} of order $\Omega(p^2/m^2)$.
Moreover, this lower bound is \emph{not} limited to $\tilde \bbeta$ itself, and holds for \emph{all} scalar-debiased solutions of the form $\gamma \cdot \tilde \bbeta$, for any $\gamma \in \RR$.
This result stands in sharp contrast to the Gaussian random projection setting discussed in \Cref{rem:inversion_bias_Gaussian}, where an appropriate scalar suffices to eliminate the inversion bias, and thus the oblique projection in subsampled OLS.

The proof of \Cref{theo:lower_bound} (see \Cref{sec:proof_of_theo:lower_bound}) is based on the construction of a counterexample which imposes a non-uniform bias across the coordinates of the subsampled OLS, so that it \emph{cannot} be corrected via a single scalar rescaling.

In the following result, we derive a variance upper bound for subsampled OLS \emph{with debiasing}.
\begin{theorem}[\textbf{Variance upper bound for subsampled OLS}]\label{theo:variance_lev_least-sq}
For a data matrix $ \X\in\RR^{n\times p}$ of rank $p$, a response vector $ \y \in\RR^n$, and random sampling matrix $\S \in \RR^{m\times p}$ with distribution $\{\pi_i\}_{i=1}^{n}$ as in \Cref{def:RS}, there exists $C > 0$ independent of $n, p$ so that for $m \geq C\theta_{\max}(\X) p\log (p/\delta)$, any $\delta \in (0,1)$, when conditioned on an event $\zeta$ that holds with probability at least $1-\delta$, the subsampled OLS $\tilde \bbeta = (\S\X)^\dagger \S \y$ in \eqref{eq:hat_bbeta} satisfies
\begin{align*}
    \Var_{\zeta}(\tilde \bbeta) \leq (1-\delta)^{-1}(\Delta(\X) + \epsilon \cdot \|\rr\|^2),
\end{align*}
with $\epsilon = O\left( \sqrt{\log(p/\delta)\theta_{\max}^3(\X)p^3/m^3} \right)$, for $\Var_{\zeta}(\cdot)$ in \Cref{def:bias_and_variance}, and
\begin{equation}\label{eq:def_Delta}
  \Delta(\X) = \rr^\top\diag\left\{\ell_i(\X)/(m\pi_i)\right\}_{i=1}^{n}\rr, \quad \rr\equiv \y - \X\bbeta_{\OLS}.
\end{equation}
\end{theorem}

\Cref{theo:variance_lev_least-sq} provides an upper bound for the variance of subsampled OLS $\tilde \bbeta$ \emph{without debiasing}.
Furthermore, recall $\theta_{\min}(\X) p/m \cdot \I_n \preceq \diag \{ \ell_{i}(\X)/(m \pi_i) \}_{i=1}^n \preceq \theta_{\max}(\X) p/m \cdot \I_n$ from the discussion after \Cref{theo:de_sam_obliqueprojec}, the variance (upper bound) in \Cref{theo:variance_lev_least-sq} can be similarly decomposed as the sum of the dominant term $\Delta(\X) = O (\| \rr \|^2 \cdot \theta_{\max}(\X)  p/m )$ and higher-order term of order $ \tilde O (\| \rr\|^2 \cdot \sqrt{\theta_{\max}^3(\X)  p^3/m^3} )$ akin to \Cref{theo:de_sam_obliqueprojec}.
The proof of \Cref{theo:variance_lev_least-sq} is given in \Cref{sec:proof_of_variance_lev_least_sq}.

With Theorems~\ref{theo:lower_bound}~and~\ref{theo:variance_lev_least-sq} (on the bias lower bound and variance upper bound of subsampled OLS) at hand, a natural question is how the \emph{debiased} oblique projection framework developed in \Cref{theo:de_sam_obliqueprojec} applies to subsampled OLS, and whether such debiasing may incur an \emph{decreased bias but increased variance}.
This question is addressed by the following result; see \Cref{sec:proof_theo_sam_ols} for the proof.

\begin{theorem}[\textbf{Precise bias--variance characterizations for debiased subsampled OLS}]\label{theo:de_sam_ols}
For a data matrix $\X \in \RR^{n\times p}$ of rank $p$ with $n\geq p$ and response vector $\y \in \RR^{n}$,  assume that $ \min_{i,\|\x_i\|>0}\|\x_i\|^2 /\|\X\|_F^2 \geq n^{-\alpha}$ for some constant $\alpha>0$. Let $ \check\S \in \RR^{m\times n}$ the debiased sampling matrix in \eqref{eq:debias_check_S}.
Then, 
there exists $C > 0$ independent of $n, p$ so that for $m \geq C\theta_{\max}(\X) p\log (p/\delta)$,  $ \delta\leq n^{-(3+\alpha)}$,  when conditioned on an event $\zeta$ that holds with probability at least $1-\delta$, the debiased subsampled OLS $\check\bbeta$ in \eqref{eq:debia_sam_ols} satisfies 
\begin{align*}
  \Bias_{\zeta}(\check\bbeta) = \epsilon^2 \cdot \|\rr\|^2,~~\text{and}~~\Var_{\zeta}(\check\bbeta) = \Delta(\X) + \epsilon \cdot \|\rr\|^2,
\end{align*}
with $\epsilon= O\Big(\sqrt{(\log n/\loglog n)^2 \cdot \theta_{\max}^3(\X) p^3/m^3}\Big)$, $\Bias_{\zeta}(\cdot)$ and $ \Var_{\zeta}(\cdot)$ in \Cref{def:bias_and_variance}, and $\Delta(\X)$ defined in \eqref{eq:def_Delta}.
\end{theorem}

Comparing the bias-variance characterizations of the \emph{debiased} solution $\check\bbeta$ in \Cref{theo:de_sam_ols} with those of the classical subsampled OLS $\tilde \bbeta$ (without debiasing) in Theorems~\ref{theo:lower_bound}~and~\ref{theo:variance_lev_least-sq}, we find that the \emph{debiased} solution $\check\bbeta$ yields a clear improvement. 
Specially, it achieves a \emph{smaller bias} (of order $\tilde O( \theta_{\max}^3(\X) p^3/m^3)$ compared with the lower bound $\Omega(p^2/m^2)$ in Theorem~\ref{theo:lower_bound} without debiasing), while incurring \emph{no larger variance} (up to constant and some log factor), provided that
$m\geq  C\theta_{\max}(\X) p\log (p/\delta)$.

For sketched OLS, it has been shown in \cite{derezinski2023Algorithmic} that LESS embeddings nearly match sub-Gaussian sketches in variance, achieving $\Var_{\zeta}(\tilde \bbeta)=\frac{p}{m-p-1}\|\rr\|^2+\tilde O(\sqrt{p}/m)\|\rr\|^2$. 
However, such a characterization remains absent for subsampled OLS under popular random sampling schemes such as leverage score sampling and SRHT, for which existing analyses are limited to variance bounds of the form $\Var_{\zeta}(\tilde \bbeta)=\tilde O(p/m)\|\rr\|^2$; see~\cite{chen2019active,david2014sketching}.  
The results in Theorems~\ref{theo:variance_lev_least-sq} and \ref{theo:de_sam_ols} close this gap by providing a finer-grained analysis for subsampled OLS, for both classical and debiased sampling schemes. 
In particular, we show that the previously established coarse bound $\tilde O(p/m)\|\rr\|^2$ can be further refined as $\tilde O(p/m)\|\rr\|^2 =\Delta(\X)+ \tilde O(\sqrt{\theta_{\max}^3(\X)p^3/m^3})\|\rr\|^2 $.

\subsection{Precise bias--variance characterizations for leverage score sampling and SRHT}\label{sub:precise_bias_variance_under_leverage_score_sampling_and_srht}

In this subsection, we focus on the two widely-used special cases of leverage score sampling and subsampled randomized Walsh--Hadamard transform (SRHT, see \Cref{def:srht} below for a formal definition).
We show that the matrix-level debiasing in \eqref{eq:debias_check_S} is indeed \emph{not} necessary (as opposed to general sampling schemes considered in \Cref{sub:bias_and_variance_of_subsampled_OLS}), at least \emph{to some extend}, for these two special cases of random sampling schemes.

Recall from our discussion in \Cref{rem:scalar_debiasing_exact} that in the case of \emph{exact} and close-to-exact leverage score sampling, our proposed debiasing approach in \eqref{eq:debias_check_S} also becomes a simple \emph{scalar} rescaling.
As a consequence, one may expect that classical subsampled OLS solution $\tilde\bbeta$ yields \emph{small bias} in a similar spirit. 
This is given in the following result,  the proof of which is given in \Cref{subsec:auxillary_resul_theodesamols_lev} 
 % of \Cref{sec:proof_theo_sam_ols}.

\begin{corollary}[\textbf{Bias--variance characterizations for debiased subsampled OLS under leverage score sampling}]\label{coro:sub_ols_lev}
Under the settings and notations of \Cref{theo:de_sam_ols}, 
for $\S\in \RR^{m\times n}$ the  sampling matrix with sampling probabilities
$\pi_i\in [\ell_{i}(\X)/(p \theta_{\max}(\X)), \ell_{i}(\X)/(p \theta_{\min}(\X))]$ with $\theta_{\min}(\X)\in [1/2,1]$ as defined in \Cref{def:approx_factor}, there exists $C> 0$, $\alpha>0$, $ \delta \leq  n^{-(3+\alpha)}
$ such that for $m \geq C\theta_{\max}(\X) p\log (p/\delta)$, when conditioned on the event $\zeta$ that holds with the probability at least $1-\delta$, the standard subsampled OLS $\tilde \bbeta$ in \eqref{eq:hat_bbeta} satisfies
\begin{align*}
\Bias_{\zeta}(\tilde \bbeta) = \epsilon^2 \cdot \|\rr\|^2,~~\text{and}~~\Var_{\zeta}(\tilde \bbeta) = \Delta(\X) + \epsilon \cdot \|\rr\|^2,
\end{align*} 
% } 
with $\epsilon= O(\sqrt{(\log n/\loglog n)^2 \cdot \theta_{\max}^3(\X) p^3/m^3}+\epsilon_{\theta} \cdot\theta_{\max}(\X) p/m )$,  $\epsilon_{\theta}=\max\{\theta^{-1}_{\min}(\X)-1,1 - \theta^{-1}_{\max}(\X)\} $,  $\Bias_{\zeta}(\cdot), \Var_{\zeta}(\cdot)$ in \Cref{def:bias_and_variance}, and $\Delta(\X)$ in \eqref{eq:def_Delta}.
\end{corollary}

Besides exact and approximate leverage score sampling discussed above, the subsampled randomized Walsh--Hadamard transform (SRHT)~\cite{ailon2006approximate} is another efficient data-oblivious sketching scheme.  

\begin{definition}[\textbf{Subsampled randomized Walsh--Hadamard transform}, \textbf{SRHT},~\cite{ailon2006approximate}]\label{def:srht}
For a data matrix $\X \in \RR^{n \times p}$ of rank $p$ with $n \geq p$ and a response vector $\y \in \RR^n$, assume without loss of generality that $n = 2^p$ for some integer $p$.
Then, the SRHT of $\X$ and $\y$ are  given by
\begin{equation*}
  \tilde{\X}_{\SRHT} = \S \H_n \D_n \X/\sqrt n \in \RR^{m \times p},\quad \tilde{\y}_{\SRHT} = \S \H_n \D_n \y/\sqrt n \in \RR^{m },
\end{equation*}
respectively, for uniform random sampling matrix $\S \in \R^{m \times n}$ (with $\pi_i = 1/n$ in \Cref{def:RS}), $\H_{n} \in \RR^{n \times n}$ the Walsh--Hadamard matrix of size $n$, and diagonal $\D_n \in \RR^{n \times n}$ having i.i.d.\@ Rademacher random variables on its diagonal.
And we define, similarly to the subsampled OLS solution $\tilde \bbeta$ in \eqref{eq:hat_bbeta}, the subsampled OLS solution using SRHT as 
\begin{equation}\label{eq:def_beta_SRHT}
  \tilde \bbeta_{\SRHT} = \tilde{\X}_{\SRHT}^\dagger \tilde{\y}_{\SRHT}.
\end{equation}
% }
\end{definition}

Since $\H_n^\top \H_n / n = \I_n$ and $\D_n^2 = \I_n$,  one has $ \X^\top \D_n \H_n^\top \H_n \D_n \X /n= \X^\top \X $ and a similar identity holds for $\X^\top \y$.
Moreover, the linear transform $\H_n \D_n$ is known to effectively spread the leverage scores of the data matrix $\X$, in the sense that, with high probability, the leverage scores of $\H_n \D_n \X / \sqrt{n}$ satisfy  $ \ell_i(\H_n \D_n \X/\sqrt{n})= \frac{p}{n}(1\pm \tilde O(\frac{\sqrt{p}}{n}))$, for all $i$, see, for example,~\cite{drineas2011faster}, \cite[Theorems~3.1 and~3.2]{tropp2011improved}, and \cite[Lemma~E.4]{niu2025fundamental}. 
As a consequence, SRHT can be interpreted as an instance of approximate leverage score sampling  applied to the transformed data matrix $\H_n \D_n \X / \sqrt{n}$, and is therefore expected to  be approximately unbiased per \Cref{coro:sub_ols_lev}.

This intuition is formalized in the following result; see \Cref{subsec:auxillary_resul_theodesamols_srht} for the proof.

\begin{corollary}[\textbf{Bias-variance characterizations for subsampled OLS under SRHT}]\label{coro:debiasing_srht} 
Under the setting and notations of \Cref{theo:de_sam_ols}, let $\X_{\HD}=\H_n \D_n\X/\sqrt n$, and   assume that $ \min_{i,\|\x_{\HD,i}\|>0}\|\x_{\HD,i}\|^2 /\|\X_{\HD}\|_F^2 \geq n^{-\alpha}$ for   some constant $\alpha>0$. 
For $\tilde{\X}_{\SRHT} \in \RR^{m \times p}$ the SRHT of $\X$ and $\tilde{\y}_{\SRHT} \in \RR^{m }$ the SRHT of $\y$ as in \Cref{def:srht}, there exists $C> 0$,  $n \exp(-p) < \delta \leq  n^{-(3+\alpha)}
$ such that
for   $m\geq C p \log (p/\delta)$, when conditioned on an event $\zeta$ that holds with probability at least $1-\delta$, the subsampled OLS with SRHT $\tilde\bbeta_{\SRHT}$ in \eqref{eq:def_beta_SRHT} satisfies 
\begin{align*}
\Bias_{\zeta}(\tilde \bbeta_{\SRHT}) = \epsilon^2 \cdot \|\rr\|^2,~~\text{and}~~\Var_{\zeta}(\tilde \bbeta_{\SRHT}) = \Delta(\X_{\HD}) +  \epsilon \cdot \|\rr\|^2,
% &\text{\emph{Bias}}:~\textbf{Bias}(\tilde\bbeta_{\SRHT})= \epsilon^2\|\rr\|^2,\\
% &\text{\emph{Variance}}:~\textbf{Var}(\tilde\bbeta_{\SRHT})=\Delta_{\check\bbeta}+\epsilon\|\rr\|^2.
\end{align*}
% } 
with $ \epsilon= O(\sqrt{(\log n/\loglog n)^2 \cdot p^3/m^3}+ \sqrt{p\log(n/\delta)}/m)$,  $\Bias_{\zeta}(\cdot)$ and $ \Var_{\zeta}(\cdot)$ in \Cref{def:bias_and_variance}, and $\Delta(\cdot)$ in \eqref{eq:def_Delta}.
% where the leverage scores used in $\Delta_{\check\bbeta}$ is $\ell_i(\H_n \D_n \X/\sqrt{n})$. 
\end{corollary}

\Cref{table:compare_ols} summarizes the existing upper bounds on the bias of subsampled OLS estimator established in prior work under uniform sampling, leverage score sampling, and SRHT~\cite{wangsketchridge2018,bartan2022distributed}, and contrasts them with our results on debiased subsampled OLS. 
We see that the proposed debiased approach consistently achieves a substantially smaller bias across all three sampling schemes.

\begin{table*}[t] 
\caption{Bias characterizations for classical subsampled OLS established in previous efforts versus that for the proposed debiased subsampled OLS in \Cref{theo:de_sam_ols},  under  uniform sampling (\textbf{UNI}), (exact/approximate) leverage score sampling (\textbf{Lev}) with $\theta_{\min}(\X)/(2\theta_{\min}(\X)-1)\leq \theta_{\max}(\X) \leq  1/(1 - \sqrt{\theta_{\max}(\X)p/m})$, and \textbf{SRHT} with $p\geq \log(n/\delta)$.
} %with mean accuracy and standard deviation
\label{table:compare_ols}
  \centering
   \setlength{\tabcolsep}{1pt}
  \begin{tabular}{lcccccccc}
    \toprule
    Reference & \textbf{UNI} &   \textbf{Lev} & \textbf{SRHT} \\
    \midrule  
    \cite{bartan2022distributed}     &  $\tilde O\left(\sqrt{\theta^3_{\max}(\X)p^3/m^3}\right)$  &  $\tilde O\left(\sqrt{p^3/m^3}\right)$  &   $\tilde O\left(\sqrt{p^3/m^3}\right)$ \\
\cite{wangsketchridge2018}    &  $\tilde O\left(\theta_{\max}^2(\X)p^2/m^2\right)$  &  $\tilde{O}\left(p^2/m^2 \right)$  &   $\tilde O\left(p^2/m^2\right)$ \\
\underline{\textbf{This work}} &  $ \tilde O\left(\theta_{\max}^3(\X)p^3/m^3\right)$  &  $\tilde O\left(p^3/m^3\right)$     &   $\tilde O\left(p^3/m^3+ p/m^2\right)$ \\
    \bottomrule
  \end{tabular}
\end{table*}

\subsection{Numerical results for subsampled OLS}\label{subsec:num_exper_ols}

\begin{figure}[b]
      \centering 
      \hspace*{-0.1\textwidth}
       \begin{subfigure}[c]{0.3\textwidth}
    \begin{tikzpicture}
    \renewcommand{\axisdefaulttryminticks}{5} 
    \pgfplotsset{every major grid/.style={densely dashed}}       
    \tikzstyle{every axis y label}+=[yshift=-10pt] 
    \tikzstyle{every axis x label}+=[yshift=5pt]
    \pgfplotsset{every axis legend/.append style={cells={anchor=east},fill=none,draw=none, at={(-0.8,1)}, anchor=north west, font=\tiny,legend columns=1,
        transpose legend }}
    \begin{axis}[
    width=1.3\columnwidth,
    height=1.1\columnwidth,
    xlabel style={font=\tiny},
        ylabel style={font=\tiny},
        tick label style={font=\tiny},
     xmin=5000,
      xmax=8000,
      ymin=1.76e-05,
      ymax=0.00067,
                    ymajorgrids=true,
                    scaled ticks=true,
                    % scaled y ticks=false,
%          ytick={0.0840,0.0848,0.0856,0.0864},
% yticklabels={0.0840,0.0848,0.0856,0.0864},
% yticklabel style={
%     /pgf/number format/fixed,
%     /pgf/number format/precision=4
% },
                    xlabel = { Sketch size},
                    ylabel = { Bias},
                    ymode=log,
                    ytick={2e-5, 1e-4, 5e-4},
yticklabels={$2\cdot 10^{-5}$,$ 10^{-4}$,$5\cdot 10^{-4}$},
                    ]

 \addplot[mark=*,color=RED,line width=0.8pt,mark options={solid,fill=RED}] coordinates{
      (5000,4.894449e-05) (6000,3.722920e-05) (7000,2.337374e-05) (8000,1.929648e-05)
        };
        \addlegendentry{{DUNI}}[ font=\tiny];
    
              \addplot[densely dashed, mark=*,color=RED,line width=0.8pt,mark options={solid,fill=RED}] coordinates{
       (5000,0.000662) (6000,0.000477) (7000,0.000353) (8000,0.000281)
        };
        \addlegendentry{{UNI}}[ font=\tiny];

      \addplot[mark=pentagon*,color=BLUE!60!white,line width=0.8pt] coordinates{
     (5000,3.850295e-05) (6000,3.247650e-05) (7000,2.355923e-05) (8000,1.750519e-05)
          };
         \addlegendentry{{DLev}}[ font=\tiny];
                  
                  \addplot[densely dashed, mark=pentagon*,color=BLUE!60!white,line width=0.8pt] coordinates{
                 (5000,4.205343e-05) (6000,3.020662e-05) (7000,2.629135e-05) (8000,1.795434e-05)
                    };
         \addlegendentry{{Lev}}[ font=\tiny];

    \addplot[mark=diamond*,color=GREEN!80!white,line width=0.8pt,mark options={solid,fill=GREEN!80!white}] coordinates{
  (5000,3.507163e-05)  (6000,2.889888e-05)  (7000,2.380111e-05)  (8000,2.337511e-05)
    };
     \addlegendentry{{DSRHT}}[ font=\tiny];
     
    \addplot[densely dashed,mark=diamond*,color=GREEN!80!white,line width=0.8pt,mark options={solid,fill=GREEN!80!white}] coordinates{
 (5000,3.902935e-05)  (6000,3.264855e-05)  (7000,2.572428e-05)  (8000,2.629124e-05)
    };
     \addlegendentry{{SRHT}}[ font=\tiny];
   
      \addplot[mark=square*,color=RED!25!BLUE,line width=0.8pt] coordinates{
                (5000,3.704359e-05) (6000,2.941441e-05) (7000,2.385388e-05) (8000,2.316585e-05)
                    };
           \addlegendentry{{LESS}}[ font=\tiny];
        
    \end{axis}
    \end{tikzpicture}
    % \captionsetup{font=scriptsize}
    % \caption{MSD}
    % \label{subfig:SLev_and_SRHT}
    \end{subfigure}
    \hspace{0.25\textwidth}
      \begin{subfigure}[c]{0.3\textwidth}
      \begin{tikzpicture}
    \renewcommand{\axisdefaulttryminticks}{4} 
    \pgfplotsset{every major grid/.style={densely dashed}}       
    \tikzstyle{every axis y label}+=[yshift=-10pt] 
    \tikzstyle{every axis x label}+=[yshift=5pt]
    \pgfplotsset{every axis legend/.append style={cells={anchor=east},fill=none,draw=none, at={(1,0.95)}, anchor=north east, font=\tiny},legend columns=1,transpose legend}
    \begin{axis}[
    width=1.3\columnwidth,
    height=1.1\columnwidth,
    xlabel style={font=\tiny},
        ylabel style={font=\tiny},
        tick label style={font=\tiny},
     xmin=5000,
      xmax=8000,
      ymin=0.01124,
      ymax=0.02147,
      ymajorgrids=true,
      % scaled ticks=true,
       scaled y ticks=false,
                ytick={0.012,0.015,0.02},
yticklabels={$0.012$,$0.015$,$0.02$},
%          ytick={0.0840,0.0848,0.0856,0.0864},
% yticklabels={0.0840,0.0848,0.0856,0.0864},
% yticklabel style={
%     /pgf/number format/fixed,
%     /pgf/number format/precision=2
% },
      xlabel = { Sketch size},
      ylabel = { Variance},
      % ymode=log
      ]

              \addplot[mark=pentagon*,color=BLUE!60!white,line width=0.8pt] coordinates{
   (5000,0.018717) (6000,0.015840) (7000,0.013455) (8000,0.011617)
          };
         % \addlegendentry{{DLev}}[ font=\tiny];
         
                     \addplot[densely dashed, mark=pentagon*,color=BLUE!60!white,line width=0.8pt] coordinates{
             (5000,0.018727) (6000,0.015760) (7000,0.013341) (8000,0.011668)
                    };
         % \addlegendentry{{Lev}}[ font=\tiny];

         \addplot[mark=*,color=RED,line width=0.8pt,mark options={solid,fill=RED}] coordinates{
     (5000,0.020157) (6000,0.016451) (7000,0.013908) (8000,0.012066)
        };
        % \addlegendentry{{DUNI}}[ font=\tiny];
    
               \addplot[densely dashed, mark=*,color=RED,line width=0.8pt,mark options={solid,fill=RED}] coordinates{
     (5000,0.021467) (6000,0.017345) (7000,0.014488) (8000,0.012668)
        };
        % \addlegendentry{{UNI}}[ font=\tiny];

    \addplot[mark=diamond*,color=GREEN!80!white,line width=0.8pt,mark options={solid,fill=GREEN!80!white}] coordinates{
 (5000,0.018360) (6000,0.015264) (7000,0.013115) (8000,0.011508)
    };
     % \addlegendentry{{DSRHT}}[ font=\tiny];

    \addplot[densely dashed,mark=diamond*,color=GREEN!80!white,line width=0.8pt,mark options={solid,fill=GREEN!80!white}] coordinates{
  (5000,0.018174) (6000,0.015333) (7000,0.013159) (8000,0.011476)
    };
      % \addlegendentry{{SRHT}}[ font=\tiny];

                       \addplot[mark=square*,color=RED!25!BLUE,line width=0.8pt] coordinates{
             (5000,0.018393) (6000,0.015203) (7000,0.013196) (8000,0.011248)

                    };
           % \addlegendentry{{LESS}}[ font=\tiny];

    \end{axis}
    \end{tikzpicture}
    % \captionsetup{font=scriptsize}
    % \caption{MSD}
    % \label{subfig:Uni}
    \end{subfigure}
    \captionsetup{skip=3pt}
    \caption{{Bias and variance as functions of the sketch size $m$, comparing debiased sampling (in \textbf{solid} lines) and standard sampling (in \textbf{dashed} lines) 
    % under  uniform sampling (DUni/Uni), shrinkage leverage score (DSLev/SLev) sampling,  and SRHT (DSRHT/SRHT); and   LESS, 
    on MSD dataset. 
    \textbf{DUNI}, \textbf{DLev}, \textbf{DSRHT} are the corresponding debiased versions.
    % Results are obtained by averaging over $500$ independent runs. 
    Expectation are estimated from $500$ independent runs. 
    }}
    % approximate leverage score (DLev/Lev),
    \label{fig:error_ols}
    \end{figure}
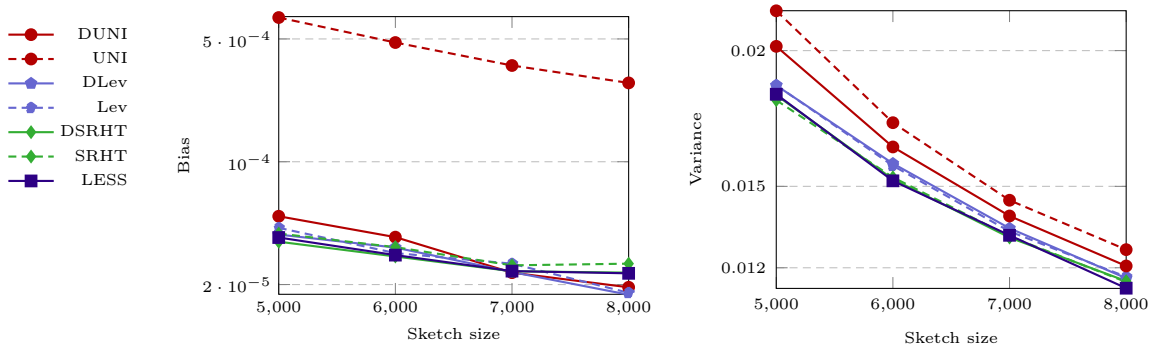

\Cref{fig:error_ols} provides empirical evidence supporting our theoretical findings on classical and debiased subsampled OLS, by evaluating the effect of sketch size $m$ on the (relative) subsampled OLS bias, $(L(\hat{\EE}[\bbeta])-L(\bbeta_{\OLS}))/L(\bbeta_{\OLS})$ and variance, $(\hat{\EE}[L(\bbeta])]-L(\bbeta_{\OLS}))/L(\bbeta_{\OLS})$, with $L(\bbeta)= \|\y - \X\bbeta\|^2$ as in \eqref{eq:def_OLS_L}.
% We consider the following OLS problem: 
% $ \min_{\bbeta\in \RR^p } \|\X\bbeta-\y\|^2,$
The data matrix $\X \in \RR^{n \times p}$ is sampled from the Million Song Year Prediction Dataset (MSD)~\cite{bertin2011million}, and $\y\in \RR^n$ is the corresponding response vector. 
See \Cref{sec:imple_detail_nmerical_exper} for implementation details and further experiments.

We see from the left panel of \Cref{fig:error_ols} that approximate leverage score sampling (\textbf{Lev}) and \textbf{SRHT} (see \Cref{def:srht}) have \emph{nearly identical} bias with and without debiasing, consistent with our \Cref{coro:sub_ols_lev} and \Cref{coro:debiasing_srht}, respectively. Uniform sampling (\textbf{UNI}) exhibits the largest bias, while its debiased counterpart using exact leverage scores (\textbf{DUNI}) substantially reduces the bias, bringing its performance close to that of Lev, SRHT, their debiased variants (\textbf{DLev} and \textbf{DSRHT}), and the random projection approach LESS~\cite{derezinski2021newtonless,garg2024distributed}.
The right panel of \Cref{fig:error_ols} shows that  Lev, SRHT, DLev, DSRHT, and LESS have \emph{comparable} variances, whereas those of UNI and DUNI are slightly larger, with DUNI marginally lower than UNI.

Overall, we see from \Cref{fig:error_ols} that debiasing empirically reduces bias without increasing variance. 
Similar trends are observed on the Flight Delay dataset; see \Cref{subsec:addition_num_ols}.

\subsection{ Extension to fast CUR decomposition}\label{sec:cur}

 CUR decompositions construct approximations from selected rows and columns of the data matrix~\cite{michael2009curmatrix}, which preserves structural properties such as sparsity, non-negativity, and interpretability.  

 Given a data matrix $\X \in \RR^{n \times p}$, a CUR decomposition is obtained
by selecting $c$ columns of $\X$ to form $\C \in \RR^{n \times c}$ and $r$ rows of $\X$ to form $\R \in \RR^{r \times p}$, and then computing $\U \in \RR^{c \times r}$ that minimizes the reconstruction error $\|  \C \U \R- \X  \|_F^2$, namely
\begin{equation}\label{eq:def_exact_U_CUR}
\U_{\CUR}
= \argmin_{\U \in \RR^{c \times r}} \|  \C \U \R- \X  \|_F^2
= \C^\dagger \X \R^\dagger,
\end{equation}
as in line with~\cite{michael2009curmatrix,sorensen2016adeim}.
Although $\U_{\CUR}$ in \eqref{eq:def_exact_U_CUR} typically achieves high accuracy,  computing it requires pseudo-inverses and large matrix multiplications, which can be costly at large scale.

Note that the full OLS estimator in \eqref{eq:def_OLS_L} is a  one-sided counterpart of  $\U_{\CUR}$  in \eqref{eq:def_exact_U_CUR} under the specialization $\X=\y\in\RR^{n}$ and  $\R=1$. Correspondingly, subsampled OLS in \eqref{eq:hat_bbeta} improves efficiency by sketching only on the left. Motivated by the same idea,~\cite{wang2016towards,ye2019fast} proposed the \emph{fast} CUR decomposition, which sketches both sides and computes
\begin{equation}\label{eq:def_U_fast}
\tilde{\U}= \argmin_{\U \in \RR^{c \times r}} \| \S_{\C} \C \U \R\S_{\R}^{\top}- \S_{\C}\X \S_{\R}^{\top} \|_F^2=(\S_{\C}\C)^\dagger \S_{\C}\X\S_{\R}^{\top}(\R\S_{\R}^{\top})^\dagger,
\end{equation}
where $\S_{\C}\in \RR^{m_c \times n}$ and $\S_{\R}\in \RR^{m_r \times p}$ are sampling matrices; see \Cref{alg:fast_cur} in \Cref{sec:proof_theo_cur} for details. 
As discussed in \Cref{sec:main_results}, the sketches   in the random oblique projections defining $\tilde{\U}$ may introduce systematic bias and could degrade  the quality of fast CUR. Motivated by this connection,  we extend our  single sketch-induced oblique-projection analysis and debiasing strategy from subsampled OLS in \eqref{eq:debia_sam_ols} to CUR,
 and propose  the following \emph{debiased} CUR solution:
\begin{equation}\label{eq:def_check_U_fast}
    \check{\U}=(\check\S_{\C}\C)^\dagger 
   (\check\S_{\C} \X \check\S^\top_{\R})
   (\R\check\S^\top_{\R})^\dagger.
\end{equation}
Here, $\check\S_{\C} \in \RR^{m_c \times n}$ and $\check\S_{\R} \in \RR^{m_r \times p}$ are the debiased sampling matrices as in \Cref{theo:de_sam_obliqueprojec}.

The following result characterizes the bias of the proposed debiased solution $\check{\U}$. 
Full bias--variance characterizations and their proofs are deferred to \Cref{sec:proof_theo_cur}.

\begin{theorem}[\textbf{Precise bias characterization for debiased fast CUR decomposition}]\label{theo:cur}
For a data matrix $\X \in \RR^{n \times p}$, let $\C\in  \RR^{n \times c}$ and $\R\in \RR^{r \times p}$ be the pre-selected columns and rows of $\X$, and let $\U_{\CUR}$ be defined as in \eqref{eq:def_exact_U_CUR}. 
For standard random sampling matrices $\S_{\C}$ and $\S_{\R} $ used in \Cref{alg:fast_cur}  in \Cref{sec:proof_theo_cur}, define their corresponding \emph{debiased} sampling matrices $ \check\S_{\C} \in \RR^{ m_{c}  \times n}$ and $ \check\S_{\R} \in \RR^{ m_{r}\times p}$ as in \eqref{eq:debias_check_S}. 
Then, there exists $C > 0$ independent of $n, p, c, r$ so that for $m_{c} \geq C\theta_{\max}(\C)c\log (c/\delta)$ and $m_{r} \geq C \theta_{\max}(\R^\top)r\log (r/\delta)$, when conditioned on an event $\zeta$ that holds with probability at least $1-\delta$, the debiased matrix $\check{\U}$ defined in \eqref{eq:def_check_U_fast} satisfies, 
\begin{align*}
    \Bias_{\zeta}(\check{\U}) = \epsilon^2 \cdot L(\U_{\CUR})
\end{align*}
with $\epsilon= \tilde O(\sqrt{ \theta_{\max}^3(\C)c^3/ m_c^3}) + \tilde  O(\sqrt{ \theta^3_{\max}(\R^\top)r^3/m_r^3}) $, for $L(\U) \equiv \|\C \U \R-\X\|^2_F$, the associated bias $\Bias_{\zeta}(\check{\U}) \equiv L(\EE_{\zeta}[\check{\U}]) - L(\U_{\CUR})$, in a spirit similar to that of subsampled OLS in \Cref{def:bias_and_variance}. Here, $\theta_{\max}(\C)$ and $\theta_{\max}(\R^\top)$ as in \Cref{def:approx_factor} denote the maximum importance sampling approximation
factors for $\C$ and $\R^\top$, respectively.
\end{theorem}

In a spirit similar to Corollaries~\ref{coro:sub_ols_lev}~and~\ref{coro:debiasing_srht}, Corollaries~\ref{coro:sub_cur_lev} and~\ref{coro:debiasing_cur_srht} in \Cref{subsec:cur_lev_srht} also show that debiasing is unnecessary for fast CUR under leverage score sampling and SRHT. For completeness, we report the numerical results for fast CUR in \Cref{subsec:numerical_cur}.

\section{Conclusion and perspectives} \label{sec:conclusion_and_perspectives}

We develop a non-asymptotic framework for random oblique projections and reveal a systematic bias arising from the nonlinear interaction between random sampling and the Moore--Penrose pseudoinverse, a phenomenon beyond the reach of classical subspace embedding analyses. 
To address this issue, we introduce a unified debiasing strategy that corrects this structural bias while preserving computational efficiency. 
When applied to subsampled least squares, our theory yields sharp bias and variance characterizations and shows that widely used schemes, including leverage score sampling and SRHT, are statistically suboptimal in general. 
For fast CUR decomposition, we construct a debiased estimator with provable improvements in approximation accuracy.

More broadly, our results highlight the need for bias-aware analysis in sketching-based algorithms. 
Extending this perspective to other methods, including sketch-and-solve Newton algorithms, randomized preconditioning, and streaming low-rank approximation, may further advance the statistical understanding of randomized numerical linear algebra.

\subsection*{Acknowledgments}
Z.~Liao and C.~Niu were supported by the National Key Research and Development Program of China (No.~2025YFA1018600), the National Natural Science Foundation of China (via fund NSFC-12571561), and the Fundamental Research Support Program of HUST (2025BRSXB0004).
MD and SG were supported in part by NSF CAREER Grant CCF-233865 and a Google ML and Systems Junior Faculty Award.

\bibliography{liao}
\bibliographystyle{plain}

%%%%%%%%%%%%%%%%%%%%%%%%%%%%%%%%%%%%%%%%%%%%%%%%%%%%%%%%%%%%
\newpage
\appendix

\begin{center}
  \textbf{\LARGE Supplementary Material of} \\ 
  \textbf{\Large Debiasing Random Oblique Projections for \\Subsampled OLS and Fast CUR in High Dimensions}
\end{center}

% \section{Supplementary material of }

 % \textbf{\large Supplementary Material of} \\ 
 %  \textbf{Debiasing Random Oblique Projections in High Dimensions \\
 %  with Applications to Subsampled OLS and Fast CUR Decomposition}

The technical appendices are organized as follows.
\begin{itemize}

\item \Cref{sec:use_lemms} collects the technical lemmas used throughout the paper.
\item \Cref{sec:proof_of_theo:lower_bound} presents the lower bound on the bias of the classical subsampled OLS estimator in \Cref{theo:lower_bound}.
\item \Cref{sec:proof_of_variance_lev_least_sq} proves the variance upper bound for the classical subsampled OLS estimator    in \Cref{theo:variance_lev_least-sq}.
\item \Cref{sec:proof_theo_sam_ols} derives  the bias and variance results for the debiased subsampled OLS estimator in \Cref{theo:de_sam_ols}.
\item \Cref{sec:proof_theo_de_sam_obliqueprojec} shows the proof of the  statistical characterizations of the debiased oblique projection in \Cref{theo:de_sam_obliqueprojec}.
\item \Cref{sec:proof_theo_cur} contains the full bias-variance characterizations and the proof of the debiased CUR decomposition in \Cref{theo:cur} and \Cref{theo:cur_bias_variance}.
\item 
\Cref{sec:imple_detail_nmerical_exper} provides implementation details for the numerical experiments in \Cref{subsec:num_exper_ols}, additional results for subsampled OLS, and numerical results for fast CUR  decomposition.
\end{itemize}

\section{Useful lemmas} \label{sec:use_lemms}

In this section, we introduce a few technical lemmas to be used in subsequent sections.

\begin{lemma}[Sherman--Morrison formula]\label{lemm:sherman-morrison}
    For an invertible matrix $\A\in \RR^{n\times n}$ and two vectors $\uu, \vv \in \RR^n$, $\A+\uu\vv^\top$ is invertible if and only if $1+ \vv^\top\A^{-1}\uu\neq 0$ and 
    \begin{align}
        (\A+ \uu\vv^\top)^{-1}=\A^{-1} -\frac{\A^{-1}\uu\vv^\top\A^{-1}}{1+ \vv^\top\A^{-1}\uu}. \nonumber
    \end{align}
    Besides, it also follows that
    \begin{align}
        (\A+ \uu\vv^\top)^{-1}\uu=\frac{\A^{-1}\uu}{1+ \vv^\top\A^{-1}\uu}.\nonumber
    \end{align}
\end{lemma}

\begin{lemma}[Properties of the Moore--Penrose Pseudoinverse]\label{lemm:matrix_produc_pesedo}
Let  $\A \in \RR^{n \times p}$, $\B \in \RR^{p \times d}$. Then,
\begin{align*}
   (\A\B)^{\dagger} = \B^{\dagger}\A^{\dagger},
\end{align*}
if any one of the following conditions holds: (1) $\A^\top\A=\I_p$; (2) $\B^\top\B=\I_d$; or (3) $\rm{rank}(\A)={\mathrm{rank}}(\B)=\mathit{p}$.
\end{lemma}

\begin{lemma}[Subspace embedding for random sampling,{~\cite[Lemma~2.7]{niu2025fundamental}}] \label{lem:sub_embed}
Given $\X\in \RR^{n\times p}$ of rank $p$ with $n\geq p$,
let $\S\in \RR^{m\times n}$ be a random sampling matrix as in \Cref{def:RS}.
Then, there exists $C > 0$ independent of $n,p$ such that for $m\geq C\theta_{\max}(\X) p \log(p/\delta )/\epsilon^2$, failure probability $ \delta \in (0,1/2)$, $\epsilon > 0$, and $\theta_{\max}(\X)$ in \Cref{def:approx_factor}, $ \X^\top\S^\top\S\X$ is an $(\epsilon,\delta)$-approximation of $\X^\top\X$. 
\end{lemma}

\Cref{lem:sub_embed} corresponds to the special case $\C=0$ of  the regularized formulation in~\cite[Lemma~2.7]{niu2025fundamental}. It establishes explicit conditions on the sample size $m$  under which  $\S\X$ forms an  $(\epsilon,\delta)$-subspace embedding of $\X$.  This result serves as a technical foundation, ensuring that the sketch $\S\X$ is well behaved throughout the subsequent analysis.

\section{Proof of Theorem~\ref{theo:lower_bound}}\label{sec:proof_of_theo:lower_bound}

% \zhenyu{The following notations to align:}
% \Sachin{Done}
We begin with the following result from \cite{derezinski2021sparse}, which is used in the proof of~\Cref{theo:lower_bound}.
\begin{lemma}[{\cite[Lemma 36]{derezinski2021sparse}}]\label{l:anti_con}
    There is a universal constant $C>0$ such that for any positive integer $b$, if $x\sim$Binomial$(b, 0.5)$ then,
\begin{align*}
    \EE\Big[\frac{1}{x+b/2}\Big] \geq \Big(1+\frac{1}{Cb}\Big)\cdot\frac{1}{b}.
\end{align*}
\end{lemma}
For notational convenience, in this proof we will denote $\bbeta_{\OLS}$ by $\bbeta^*$. Consider a matrix $\X$ whose first $2p$ rows are given as follows:
\begin{align*}
    \X_{[1:2p]} = \begin{bmatrix}
        \frac{1}{2} &0 &0 &\cdots &0\\
        \frac{\sqrt{3}}{2} &0 &0 &\cdots &0\\
        0 &\frac{1}{2} &0 &\cdots &0\\
        0 &\frac{\sqrt{3}}{2} &0 &\cdots &0\\
        0  &0 &\frac{1}{2} &\cdots &0\\
        0  &0 &\frac{\sqrt{3}}{2} &\cdots &0\\
          0  &0 &0 &\cdots &\frac{1}{2}\\
        0  &0 &0 &\cdots &\frac{\sqrt{3}}{2}\\
    \end{bmatrix}.
\end{align*}
Also, if $n > 2p$, let rows subsequent to $(2p)^{th}$ row be all zeros. In addition, consider the response vector $\y \in \RR^n$ as follows:
\begin{align*}
    [\y]_i = \begin{cases}
        -1,  \ \ &\text{$i\leq p$ and $i$ is odd,}\\
        1, \ \  &\text{$i\leq p$ and $i$ is even,}\\
        1, \ \  &\text{$p<i\leq 2p$,}\\
        0, \ \  &\text{$i>2p$}.
    \end{cases}
\end{align*}
Note that $\X^\top\X=\I_p$, the exact least squares estimator $\bbeta^*$ can be found as $\bbeta^* = \X^\top\y$. 
The entries of $\bbeta^*$ are given as
\begin{align*}
    \bbeta^*_i = \begin{cases}
        \frac{\sqrt{3}}{2}-\frac{1}{2}, \ \ \text{$i \leq \frac{p}{2}$},\\
        \frac{\sqrt{3}}{2} + \frac{1}{2}, \ \ \text{$\frac{p}{2}<i \leq p$}.
    \end{cases}
\end{align*}
Also, note that $\|\y\|^2 =2p$, $\|\X\bbeta^*-\y\|^2=\Theta(p)$, $\|\X\bbeta^*\|^2=\Theta(p)$. The exact leverage scores of $\X$ lie in the set $\{\frac{1}{4},\frac{3}{4},0\}$. Consider the probability distribution $\{\pi_i\}_{i=1}^{n}$, where $\pi_i = \frac{1}{2p}$ for $i\leq 2p$ and $\pi_i=0$ for $i > 2p$. Note that $\{\pi_i\}_{i=1}^{n}$, provides a $\frac{1}{2}$-approximation to exact leverage score sampling. Let $m\geq p$, $\S \in \RR^{m\times n}$ be sampling matrix sampled from the distribution $\{\pi_i\}_{i=1}^{n}$. Let $\tilde\bbeta = (\X^\top\S^\top\S\X)^\dagger\X^\top\S^\top\S\y$ be the sketched least squares estimator and $\bar\bbeta = \EE[\tilde\bbeta]$. We can specify the distribution of every entry in $\tilde\bbeta$. To see that, we start by noting that $\X^\top\S^\top\S\X$ is a diagonal matrix where
\begin{align*}
    [\X^\top\S^\top\S\X]_{ii} = \frac{3}{4}\cdot s_{2i} + \frac{1}{4}\cdot s_{2i-1}, \ \ \text{for} \ \ 1\leq i \leq p.
\end{align*}
where $s_{2i}$ and $s_{2i-1}$ denote the number of times $(2i)^{th}$ and $(2i-1)^{th}$ rows are selected while sampling $\S$ from $\{\pi_i\}_{i=1}^{n}$. Therefore,
\begin{align*}
    [(\X^\top\S^\top\S\X)^\dagger]_{ii} = \begin{cases}
        \frac{4}{3s_{2i} + s_{2i-1}}, \ \ &\text{if $s_{2i} + s_{2i-1} >0$,}\\
        0,  \ \ &\text{otherwise.}
    \end{cases}
\end{align*}
The $i^{th}$ entry of $\tilde\bbeta$ is given as
\begin{align*}
    [\tilde\bbeta]_i = \begin{cases}
        \frac{2(\sqrt{3}s_{2i}-s_{2i-1})}{3s_{2i}+s_{2i-1}},  \ \ \ &\text{$i \leq \frac{p}{2}$ and $s_{2i} + s_{2i-1} >0$,}\\
        0, \ \ \  &\text{$i \leq \frac{p}{2}$ and $s_{2i} + s_{2i-1} =0$,}\\
        \frac{2(\sqrt{3}s_{2i}+s_{2i-1})}{3s_{2i}+s_{2i-1}}, \ \ \ &\text{$\frac{p}{2}<i \leq p$ and $s_{2i} + s_{2i-1} >0$,}\\
        0, \ \ \ &\text{$\frac{p}{2} <i \leq p$ and $s_{2i} + s_{2i-1} =0$}.
    \end{cases}
\end{align*}
Note that, for $1\leq i \leq \frac{p}{2}$, the distribution of $[\tilde\bbeta]_i's$ are identical, and for $i>\frac{p}{2}$, the distribution of $[\tilde\bbeta]_i's$ are identical. Due to this observation, $\EE[\tilde\bbeta]_i = \EE[\tilde\bbeta]_1$ for all $i \leq \frac{p}{2}$ and $\EE[\tilde\bbeta]_i =\EE[\tilde\bbeta]_p$ for all $i > \frac{p}{2}$. We first provide an upper bound on $\EE[\tilde\bbeta]_1$. Let $\zeta$ denote the event $ 
\{s_{2i}+s_{2i-1}>0, \forall i\}
$. We have,
\begin{align*}
    \EE_{\zeta}[\tilde\bbeta]_1 =  \EE_{\zeta}\Big[\frac{2(\sqrt{3}s_{2}-s_{1})}{3s_{2}+s_{1}} \Big]= \EE_{\zeta}\Big[\frac{2\big((1+\sqrt{3})s_2 - (s_1+s_2)\big)}{2s_2+(s_1+s_2)} \Big].
\end{align*}
Let $b:=s_1+s_2$ and $x:=s_2$, we get
\begin{align*}
    \EE_{\zeta}[\tilde\bbeta]_1 &=\EE_{\zeta}\Big[\frac{2\big((1+\sqrt{3})x - b\big)}{2x+b} \Big]\\
    &=\EE_{\zeta}\Big[\frac{(1+\sqrt{3})(2x+b-b)}{2x+b}-\frac{2b}{2x+b} \Big]\\
    &=\EE_{\zeta}\Big[1+\sqrt{3} - \frac{(3+\sqrt{3})b}{2x+b} \Big]\\
    &=(1+\sqrt{3})\cdot\EE_{\zeta}\Big[1-\frac{b\sqrt{3}}{2}\cdot\frac{1}{x+b/2} \Big].
\end{align*}
Note that for any given fixed $b>0$, $q \sim \text{Binomial}(x,0.5)$. Using   Lemma \ref{l:anti_con}, we get $\EE\Big[\frac{1}{x+b/2}\Big] \geq \Big(1+\frac{1}{Cb}\Big)\cdot\frac{1}{b}$ for some universal constant $C$. Therefore,
\begin{align*}
    \EE_{\zeta}[\tilde\bbeta]_1 &\leq (1+\sqrt{3})\cdot\EE_{\zeta}\Big[1 - \frac{\sqrt{3}}{2}\big(1+\frac{1}{Cb}\big)\Big]\\
    &= (1+\sqrt{3})\cdot\Big(1-\frac{\sqrt{3}}{2}- \frac{\sqrt{3}}{2C}\cdot\EE_{\zeta}\Big[\frac{1}{b} \Big]\Big)\\
    &=(1+\sqrt{3})\cdot\Big(1-\frac{\sqrt{3}}{2}- \frac{\sqrt{3}}{2C}\cdot\sum_{t=1}^{m}{\frac{1}{t}\Pr(b=t)} \Big)\\
    &\leq (1+\sqrt{3})\cdot\Big(1-\frac{\sqrt{3}}{2}- \frac{\sqrt{3}}{2C}\cdot\sum_{t=1}^{2m/p}{\frac{1}{t}\Pr(b=t)} \Big)\\
    &\leq (1+\sqrt{3})\cdot\Big(1-\frac{\sqrt{3}}{2}- \frac{p\sqrt{3}}{4mC}\cdot\sum_{t=1}^{2m/p}{\Pr(b=t)} \Big)\\
    &\leq (1+\sqrt{3})\cdot\Big(1-\frac{\sqrt{3}}{2}- \frac{p\sqrt{3}}{8mC} \Big),
\end{align*}
where the last inequality holds because $\sum_{t=1}^{2m/p}{\Pr(b=t)} \geq \frac{1}{2}$, as $b \sim \text{Binomial}(m,\frac{1}{p})$. Therefore, we get,
\begin{align*}
    \EE_{\zeta}[\tilde\bbeta]_1 &\leq 1 - \frac{\sqrt{3}}{2} + \sqrt{3}-\frac{3}{2} - \frac{3+\sqrt{3}}{8C}\cdot\frac{p}{m}=\frac{\sqrt{3}}{2}-\frac{1}{2}-\frac{3+\sqrt{3}}{8C}\cdot\frac{p}{m}\\
    &=[\bbeta^*]_1 - \frac{3+\sqrt{3}}{8C}\cdot\frac{p}{m}.
\end{align*}
Similarly, we can show that,
\begin{align*}
     \EE_{\zeta}[\tilde\bbeta]_p \geq \Big([\bbeta^*]_p + \frac{3-\sqrt{3}}{8C}\cdot\frac{p}{m}\Big).
\end{align*}
Therefore, we get
% \Chengmei{Below, $\EE[\tilde\bbeta]_p$ should be $\EE[\tilde\bbeta]_i$, no?}
\begin{align*}
    \EE_{\zeta}[\tilde\bbeta]_i &= \Big([\bbeta^*]_{i} -c_1\Big) &\text{for $i \leq \frac{p}{2}$,} \\
    &= \Big([\bbeta^*]_{i} + c_2\Big) &\text{for $i > \frac{p}{2}$},% \label{e:lb_1}
\end{align*}
for some $c_1, c_2 = \Omega(\frac{p}{m})$. Let $\bar\bbeta =\EE_\zeta[\tilde\bbeta]$ and $\gamma$ be any scalar factor. We have,
\begin{align}
    \|\X\gamma\bar\bbeta-\y\|^2 &= \|\X\bbeta^*-\y\|^2 + \|\X(\gamma\bar\bbeta-\bbeta^*)\|^2 \nonumber\\
    &=\|\X\bbeta^*-\y\|^2 + \|\X_{[1:p,1:\frac{p}{2}]}(\gamma\bar\bbeta_{[1:\frac{p}{2}]} - \bbeta^*_{[1:\frac{p}{2}]})\|^2 \nonumber \\ 
    &+ \|\X_{[p+1:2p,\frac{p}{2}+1:p]}(\gamma\bar\bbeta_{[\frac{p}{2}+1:p]} - \bbeta^*_{[\frac{p}{2}+1:p]})\|^2, \label{e:lb_2}
\end{align}
where by $\X_{[p+1:2p,\frac{p}{2}+1:p]}$ we mean the sub matrix of $\X$ formed by $(p+1)^{th}$ to $(2p)^{th}$ rows and $(\frac{p}{2}+1)^{th}$ to $p^{th}$ columns. The other subscripts can be understood similarly. First, consider the term $\|\X_{[1:p,1:\frac{p}{2}]}(\gamma\bar\bbeta_{[1:\frac{p}{2}]} - \bbeta^*_{[1:\frac{p}{2}]})\|^2$. We get,
% \Chengmei{In the first line, $\bbeta^\star$ should be $\bbeta^*_{[1:\frac{p}{2}]}$, no? }\Sachin{Yes, fixed}
\begin{align*}
    \Big\|\X_{[1:p,1:\frac{p}{2}]}(\gamma\bar\bbeta_{[1:\frac{p}{2}]} - \bbeta^*_{[1:\frac{p}{2}]})\Big\|^2 &= \Big\|\X_{[1:p,1:\frac{p}{2}]}\big(\gamma\cdot(\bbeta^*_{[1:\frac{p}{2}]} - c_1\mathbf{1}) - \bbeta^*_{[1:\frac{p}{2}]}\big)\Big\|^2\\
    &=\Big\|\X_{[1:p,1:\frac{p}{2}]}\big(\gamma\cdot c_{1}\mathbf{1} + (1-\gamma)\cdot \bbeta^*_{[1:\frac{p}{2}]}\big)\Big\|^2\\
    &=\Big(\gamma\cdot c_1 + (1-\gamma)\cdot [\bbeta^*]_{1}\Big)^2\cdot\|\X_{[1:p,1:\frac{p}{2}]}\mathbf{1}\|^2,
\end{align*}
where $\mathbf{1}$ denotes a vector of $1$'s. As $\|\X_{[1:p,1:\frac{p}{2}]}\mathbf{1}\|^2 =\Theta(p)$, and also recall $\|\X\bbeta^*-\y\|^2 = \Theta(p)$, we get
\begin{align*}
    \Big\|\X_{[1:p,1:\frac{p}{2}]}(\gamma\bar\bbeta_{[1:\frac{p}{2}]} - \bbeta^*_{[1:\frac{p}{2}]})\Big\|^2 &= \Big(\gamma\cdot c_1 + (1-\gamma)\cdot [\bbeta^*]_{1}\Big)^2\cdot\Theta(\|\X\bbeta^*-\y\|^2).
\end{align*}
% For notational convenience, we denote $\gamma$ as $z$.
If $0 < \gamma\leq 1$, then we get
\begin{align*}
    \Big(\gamma\cdot c_1 + (1-\gamma)\cdot [\bbeta^*]_{1}\Big)^2 > \max\Big\{\gamma^2\cdot\Omega\big(\frac{p^2}{m^2}\big) , (1-\gamma)^2\cdot[\bbeta^*]_1^2\Big\}.
\end{align*}
At least one of $\gamma$ or $1-\gamma$ is at least $\frac{1}{2}$. In either case we get,
\begin{align*}
    \Big(\gamma\cdot c_1 + (1-\gamma)\cdot [\bbeta^*]_{1}\Big)^2 = \Omega\big(\frac{p^2}{m^2}\big).
\end{align*}
Therefore, the theorem holds for any $0<\gamma \leq 1$. In case  $\gamma<0$ or $\gamma \geq 1$,  consider the last term in relation (\ref{e:lb_2}). We have,
\begin{align*}
\|\X_{[p+1:2p,\frac{p}{2}+1:p]}(\gamma\bar\bbeta_{[\frac{p}{2}+1:p]} - \bbeta^*_{[\frac{p}{2}+1:p]})\|^2 &= \|\X_{[p+1:2p,\frac{p}{2}+1:p]}(\gamma\cdot(\bbeta^*+c_2\mathbf{1})_{[\frac{p}{2}+1:p]} - \bbeta^*_{[\frac{p}{2}+1:p]})\|^2\\
&=\|\X_{[p+1:2p,\frac{p}{2}+1:p]}(\gamma\cdot c_{2}\mathbf{1} - (1-\gamma)\cdot\bbeta^*_{[\frac{p}{2}+1:p]})\|^2\\
&=\Big(\gamma\cdot c_2 - (1-\gamma)\cdot[\bbeta^*]_p\Big)^2\cdot\Theta(\|\X\bbeta^*-\y\|^2)\\
&=\Omega\big(\frac{p^2}{m^2}\big)\cdot\Theta(\|\X\bbeta^*-\y\|^2).
\end{align*}
The last inequality holds as either $\gamma<0$ or $\gamma \geq 1$. 
This completes the proof of Theorem \ref{theo:lower_bound}

\section{Proof of Theorem~\ref{theo:variance_lev_least-sq}}\label{sec:proof_of_variance_lev_least_sq}
Let $\U \in \RR^{n\times p}$ has orthonormal columns spanning the column space of $\X$. Then, $\bbeta_{\OLS} = \U^\top\y$ and $\tilde\bbeta = (\U^\top\S^\top\S\U)^\dagger\U^\top\S^\top\S\y$. Furthermore,
\begin{align*}
    L(\tilde\bbeta) - L(\bbeta_{\OLS}) &=\big\|\U(\tilde\bbeta-\bbeta_{\OLS})\big\|^2\\
    &=\big\|\tilde\bbeta-\bbeta_{\OLS}\big\|^2.
\end{align*}
Let $\zeta$ denote the event that
\begin{align*}
   (1+\epsilon)^{-1}\cdot\I \preceq \U^\top\S^\top\S\U \preceq (1+\epsilon)\cdot\I,
\end{align*}
for $\epsilon = \sqrt{\frac{3p\theta_{\max}(\X)\log(2p/\delta)}{m}}$. First, note that for $m > 12p\theta_{\max}\log(2p/\delta) $, $\Pr(\zeta) \geq 1-\delta$. This claim is essentially the subspace embedding guarantee provided by the approximate leverage sampling matrix $\S$ (Lemma  \ref{lem:sub_embed}). In what follows, we will upper bound $\EE_\zeta[L(\tilde\bbeta)- L(\bbeta_{\OLS})]$. We have,
\begin{align*}
    \EE_\zeta[L(\tilde\bbeta)- L(\bbeta_{\OLS})] &= \EE_\zeta[\|\tilde\bbeta-\bbeta_{\OLS}\|^2]\\
    &=\EE_\zeta\Big[\Big\|(\U^\top\S^\top\S\U)^{-1}\big(\U^\top\S^\top\S\y-(\U^\top\S^\top\S\U)\bbeta_{\OLS}\big)\Big\|^2\Big]\\
    &=\EE_\zeta\Big[\Big\|(\U^\top\S^\top\S\U)^{-1}\U^\top\S^\top\S\big(\y-\U\bbeta_{\OLS}\big)\Big\|^2\Big]\\
    &\leq \EE_\zeta\Big[\big\|(\U^\top\S^\top\S\U)^{-1}\big\|^2\cdot\big\|\U^\top\S^\top\S\rr\big\|^2\Big]\\
    &\leq (1+3\epsilon)^2\cdot\EE_\zeta\Big[\|\U^\top\S^\top\S\rr\||^2\Big],
\end{align*}
where $\rr$ denotes $\y-\U\bbeta_{\OLS}$. The remaining term to upper bound is $\EE_\zeta\Big[\|\U^\top\S^\top\S\rr\|^2\Big]$. It is straightforward to show that $\EE\Big[\|\U^\top\S^\top\S\rr\|^2] = \rr^\top\diag\Big(\frac{\ell_i(\X)}{m\pi_i}\Big)\rr$. Note that,
\begin{align*}
    \EE\Big[\|\U^\top\S^\top\S\rr\|^2\Big] = \EE_\zeta\Big[\|\U^\top\S^\top\S\rr\|^2\Big]\cdot\Pr(\zeta) + \EE_{\neg \zeta}\Big[\|\U^\top\S^\top\S\rr\|^2\Big]\cdot\Pr(\neg \zeta).
\end{align*}
Therefore,
\begin{align*}
\EE_\zeta\Big[\|\U^\top\S^\top\S\rr\|^2\Big] \leq  \frac{1}{\Pr(\zeta)}\cdot\rr^\top\diag\Big(\frac{\ell_i(\X)}{m\pi_i}\Big)\rr.
\end{align*}
Combining everything we get,
\begin{align}
     \EE_\zeta[L(\tilde\bbeta)- L(\bbeta_{\OLS})]&\leq (1+3\epsilon)^2\cdot\frac{1}{\Pr(\zeta)}\cdot\rr^\top\diag\Big\{\frac{\ell_i(\X)}{m\pi_i}\Big\}_{i=1}^{n}\rr \nonumber\\
     &\leq (1+15\epsilon)\cdot\frac{1}{\Pr(\zeta)}\cdot\rr^\top\diag\Big\{\frac{\ell_i(\X)}{m\pi_i}\Big\}_{i=1}^{n}\rr\nonumber\\
    &\leq  \Big(\Delta + 15\epsilon'\|\rr\|^2\Big)\cdot\frac{1}{\Pr(\zeta)},\nonumber%\label{e:variance_hp}
\end{align}
        where $\Delta =\rr^\top\diag\Big\{\frac{\ell_i(\X)}{m\pi_i}\Big\}_{i=1}^{n}\rr$ and $\epsilon' = \sqrt{\frac{3p^3\theta^3_{\max}(\X)\log(2p/\delta)}{m^{3}}} $. This proves Theorem \ref{theo:variance_lev_least-sq}.

\section{Proof of Theorem~\ref{theo:de_sam_ols}} \label{sec:proof_theo_sam_ols}

In this section, we first provide  RMT intuition for the bias of the debiased estimator in \Cref{subsec:RMT_intuition_ols_bias}.
The corresponding intuition for the variance follows analogously and is omitted for brevity.
We then present the detailed proof of \Cref{theo:de_sam_ols} in \Cref{subsec:detail_proof_ols_bias} and \Cref{subsec:detail_proof_ols_variance}.
We next provide the proofs of auxiliary results related to \Cref{theo:de_sam_ols} in \Cref{subsec:auxillary_resul_theodesamols_lev} and \Cref{subsec:auxillary_resul_theodesamols_srht}, with additional lemmas from \cite{niu2025fundamental} deferred to \Cref{subsec:auxil_results_proof_de_sam_pls}.

\subsection{RMT intuition for Theorem~\ref{theo:de_sam_ols}}
\label{subsec:RMT_intuition_ols_bias}

In this subsection, we provide a heuristic derivation of the bias of the debiased estimator  in \Cref{theo:de_sam_ols}.
We begin by recalling some notations from \Cref{theo:de_sam_ols}.
 Let $\z^\top_s=\ee^\top_{i_s}/\sqrt{\pi_{i_s}}\X \in \RR^p$ and $\w_s^\top=\ee^\top_{i_s}/\sqrt{\pi_{i_s}}$ so that $\EE[\z_s{\z}^{\top}_s]=\X^\top\X$ and $\EE[\z_s\w_s^\top]=\X^\top$.
Denote 
% \Sachin{$\Q_{-s}$ should also have $F_{i_s}$ terms, no?}\Chengmei{Done}
\begin{equation*}%\label{eq:def_H_Q_proof_theo:sam_ols}
\Q=(\X^\top\check \S^\top\check\S \X )^{-1}= \left(\sum^{m}_{s=1}\frac{1}{m}F_{i_s,i_s}\z_s{\z}^{\top}_s \right)^{-1},~~  \Q_{-s}= \left(\sum_{k\neq s}\frac{1}{m}F_{i_s,i_s}\z_k{\z}^{\top}_k  \right)^{-1},
\end{equation*}
where $\Q_{-s}$ is independent of  $\z_s$, and $F_{i_s,i_s}$ is a deterministic term whose explicit form will be specified later.
Recalling $\rm{rank}(\X)=\mathit{p}$, we  rewrite  the  estimator %\Sachin{Shouldn't this be $\check{\S}$}\Chengmei{I agreeed.}
$( \check\S\X)^\dagger \check\S=\Q\X^\top   \check{\S}^\top \check{\S}$ and the pseudoinverse $\X^\dagger=(\X^\top\X)^{-1}\X^\top$.
Here, our objective is to choose the diagonal entries  $F_{ii}$  such that $\EE[\Q\X^\top   \check{\S}^\top \check{\S}]\simeq \X^\dagger$.
Considering 
\begin{equation*}
 \sum^{m}_{s=1}  \frac{1}{m} \EE[F_{i_si_s} \z_{s}{\z}^{\top}_{s}]= \sum^{n}_{i=1}  F_{ii}{\x}_{{i}}{\x}_{{i}}^\top,
\end{equation*}
together with  Sherman-Morrison formula in \Cref{lemm:sherman-morrison}, 
we first obtain 
 \begin{align*}
 \EE[{\Q}\X^\top   \check{\S}^\top \check{\S}]&=\EE \left[\frac{{\Q}_{-s}F_{i_si_s}\z_s\w_s^\top}{1+F_{i_si_s}\z_s^\top{\Q}_{-s}\z_s/m}\right]= \sum^{n}_{i=1}\EE\left[\frac{{\Q}_{-s}F_{ii}\X^\top\ee_{i}\ee_{i}^\top}{1+F_{ii}\ee_{i}^\top\X{\Q}_{-s}\X^\top\ee_{i}/m\pi_{i}}\right]
 % & \sum^{n}_{i=1}\EE\left[\frac{\check{\Q}_{-i}F_{ii}\X^\top\ee_{i}\ee_{i}^\top\y}{1+F_{ii}\ee_{i}^\top\X(\X^\top\X)^{-1}\X^\top\ee_{i}/m\pi_{i}}\right].\label{eq:equive_qassaq}
 \end{align*}
Let $\ell_{i}(\X)=\ee_{i}^\top\X(\X^\top\X)^{-1}\X^\top\ee_{i}$  the $i^{\rm{th}}$ leverage score of $\X$.   Applying the rank-one perturbation formula for matrix inversion (see, e.g., \cite[Lemma~2.6]{silverstein1995empirical}), we obtain 
\begin{align*}
 \EE[{\Q}\X^\top   \check{\S}^\top \check{\S}] & \simeq  \sum^{n}_{i=1}\EE\left[\frac{{\Q}_{-s}F_{ii}\X^\top\ee_{i}\ee_{i}^\top}{1+F_{ii}\ee_{i}^\top\X(\X^\top\X)^{-1}\X^\top\ee_{i}/m\pi_{i}}\right]\simeq \EE[ \Q]  \X^\top\sum^{n}_{i=1}\frac{F_{ii}\ee_{i}\ee_{i}^\top}{1+F_{ii}\ell_{i}(\X)/m\pi_{i}}.
    % &\simeq (\X^\top\X+\lambda\I)^{-1}\X^\top\sum^{n}_{i=1}\frac{F_{ii}\ee_{i}\ee_{i}^\top}{1+F_{ii}\ell_{i}(\X)/m\X\pi_{i}}\Y.
\end{align*}
We now choose 
\begin{align}\label{eq:F_ii}
   F_{ii}=  \frac{m\pi_{i}}{m\pi_{i} - \ell_{i}(\X) }
\end{align}
so that 
\begin{align*}
    \sum^{n}_{i=1}\frac{F_{ii}\ee_{i}\ee_{i}^\top}{1+F_{ii}\ell_{i}(\X)/m\pi_{i}}=\I_n.
\end{align*}
Combined with Proposition~3.2 in~\cite{niu2025fundamental} and the choice of
$F_{ii}$ in \eqref{eq:F_ii}, which  together imply  that  $ \EE[{\Q}]\simeq (\X^\top\X)^{-1}$, 
we get
\begin{equation*}
  \EE[{\Q}\X^\top   \check{\S}^\top \check{\S}] \simeq(\X^\top\X)^{-1}\X^\top.
\end{equation*}
This choice leads to the following debiased sampling matrix
 \begin{align}
  \check\S =\diag \left\{\sqrt{m\pi_{i_s}/(m\pi_{i_s}-\ell_{i_s}(\X))} \right\}^m_{s=1}\cdot\S, \quad i_s\in \{1,\ldots,n\}.\nonumber
\end{align}

\subsection{Detailed proof of the bias of debiased estimator  in  Theorem~\ref{theo:de_sam_ols}}
\label{subsec:detail_proof_ols_bias}
For ease of the subsequent analysis, we first rewrite the bias of the debiased estimator:
\begin{align*}
   \textbf{Bias}(\check\bbeta)&=\EE_{\zeta}[\check\bbeta]^\top\X^\top\X\EE_{\zeta}[\check\bbeta]-2\y^\top \X\EE_{\zeta}[\check\bbeta]- \bbeta_{OLS}^\top \X^\top\X\bbeta_{OLS} +2\y^\top \X\bbeta_{OLS}\\
     &=\EE_{\zeta}[\check\bbeta]^\top\X^\top\X\EE_{\zeta}[\check\bbeta]-2\bbeta_{OLS}^\top\X^\top \X\EE_{\zeta}[\check\bbeta]- \bbeta_{OLS}^\top \X^\top\X\bbeta_{OLS} +2\bbeta_{OLS}^\top\X^\top \X\bbeta_{OLS}\\
     &=\EE_{\zeta}[\check\bbeta]^\top\X^\top\X\EE_{\zeta}[\check\bbeta]-2\bbeta_{OLS}^\top\X^\top \X\EE_{\zeta}[\check\bbeta] +\bbeta_{OLS}^\top\X^\top \X\bbeta_{OLS}\\
     &=(\EE_{\zeta}[\check\bbeta]-\bbeta_{OLS})^\top \X^\top \X(\EE_{\zeta}[\check\bbeta]-\bbeta_{OLS})=(\EE_{\zeta}[\check\bbeta]-\bbeta_{OLS})^\top \H(\EE_{\zeta}[\check\bbeta]-\bbeta_{OLS})\\
     &=\Big\|\H^{\frac{1}{2}}(\EE_{\zeta}[\check\bbeta]-\bbeta_{OLS})\Big \|^2,
 \end{align*}
with $\H=\X^\top\X$. %\Sachin{$\H$ in place of $\bar\H$ in the last equality relation.}\Chengmei{Done}
This reformulation shows that  the proof of the bias bound for the debiased estimator proceeds in two main steps:
\begin{enumerate}
  \item construct a high probability event $\zeta$, based on subspace-embedding-type results in \Cref{lem:sub_embed}, under which the inverse $(\X^\top \check\S^\top \check\S \X)^{-1}$ is well conditioned; and
  \item conditioned on that event $\zeta$, bound the quantity  $ \|\H^{1/2}(\EE_{\zeta}[\check\bbeta]-\bbeta_{OLS}) \|^2$ using a ``leave-one-out'' type analysis.
\end{enumerate}

We begin by constructing a  high-probability event $\zeta$.
Without loss of generality, assume that $t=m/4$ is an integer. We partition the index set $\{1,\ldots,m\}$ into four blocks of size $t$, and define the events:
\begin{align}\label{eq:events}
    \zeta_j:\sum^{tj}_{s=t(j-1)+1}\frac{1}{t}\z_s{\z}^{\top}_s\succeq \frac{1}{2}   \X^\top\X ,~~~j=1,2,3,4,~~~ \zeta=\bigcap^4_{j=1}\zeta_j.
\end{align}
By $F_{i_si_s}\geq 1$, each event $ \zeta_j$ further implies
\begin{align}
 \sum^{tj}_{s=t(j-1)+1}\frac{1}{t}F_{i_si_s}\z_s{\z}^{\top}_s   \succeq  \sum^{tj}_{s=t(j-1)+1}\frac{1}{t}\z_s{\z}^{\top}_s   \succeq \frac{1}{2}   \X^\top\X, ~~~j=1,2,3,4.\nonumber
\end{align}
Intuitively, the event  $\zeta_j$  ensures that the (weighted) average of the rank-one  matrices $\z_s\z_s^\top$ over the  $j$-th block forms a sketch of size $t$ that provides a   ``lower'' $1/2$-spectral-approximation of $\X^\top\X$, in the sense of \Cref{def:rela_error_approxi}. 

Under random sampling in  \Cref{def:RS}, the events $\zeta_1$, $\zeta_2$, $\zeta_3$,  $\zeta_4$ are mutually independent. Consequently, for any index $s\in \{1,\ldots, m\}$,  there exists a block index $j=j(s)\in\{1,2,3,4\}$ such that
\begin{enumerate}
    \item the event $  \zeta_j$ is independent of $\z_{s}$; and 
    \item conditioning  on $  \zeta_j$, it follows that  $ \Q\preceq \Q_{-s}\preceq 8\H^{-1}$.
\end{enumerate}

Denote $\check\z^\top_s=\ee^\top_{i_s}/\sqrt{\pi_{i_s}}\X_{\H} \in \RR^p$ with $\X_{\H}=\X\H^{-1/2}$, let 
\begin{equation*}
\check\Q=\H^{\frac{1}{2}}\Q\H^{\frac{1}{2}}=(\X_{\H}^\top\check \S^\top\check\S \X_{\H} )^{-1}= \left(\sum^{m}_{s=1}\frac{1}{m}F_{i_s,i_s}\check\z_s\check{\z}^{\top}_s \right)^{-1},
\end{equation*}
and
\begin{equation*}
  \check\Q_{-s}=\H^{\frac{1}{2}}\Q_{-s}\H^{\frac{1}{2}}= \left(\sum_{k\neq s}\frac{1}{m}F_{i_k,i_k}\check\z_k\check{\z}^{\top}_k  \right)^{-1}.
\end{equation*}

Letting $\rr=\y-\X\bbeta_{OLS}$ and $\check a_s=\ee^\top_{i_s}/\sqrt{\pi_{i_s}}\rr$, we now rewrite
\begin{align*}
    \EE_{\zeta}[\check\bbeta]-\bbeta_{OLS}=  \EE_{\zeta}[(\X^\top\check \S^\top\check\S \X)^{-1}\X^\top\check \S^\top\check\S \rr],
\end{align*}
which gives 
\begin{align*}
   \H^{\frac{1}{2}}(\EE_{\zeta}[\check\bbeta]-\bbeta_{OLS}) &= \H^{\frac{1}{2}} \EE_{\zeta}[(\X^\top\check\S^\top\check\S\X)^{-1}\X^\top\check\S^\top\check\S\rr]=\EE_{\zeta}[(\X_{\H}^\top\check\S^\top\check\S\X_{\H})^{-1}\X_{\H}^\top\check\S^\top\check\S\rr]\\
   &=\EE_{\zeta}[\check\Q\X_{\H}^\top\check\S^\top\check\S\rr]=\frac{1}{m}\sum^m_{s=1}\EE_{\zeta}[\check\Q F_{i_si_s}\check\z_s\check a_s]=\EE_{\zeta}\left[\frac{\check\Q_{-s}F_{i_si_s}\check\z_s\check a_s}{1+\frac{1}{m}F_{i_si_s}\check\z_s^\top\check\Q_{-s}\check\z_s}\right]\\
 &=\underbrace{\EE_{\zeta}\left[\check\Q_{-s}\check\z_s \check a_s\left(\frac{F_{i_si_s}}{1+\frac{1}{m}F_{i_si_s}\check\z_s^\top\check\Q_{-s}\check\z_s}-1\right)\right]}_{\TT_1}+\underbrace{\EE_{\zeta}[\check\Q_{-s}\check\z_s \check a_s]}_{\TT_2}.
\end{align*}
 This decomposition yields
\begin{align}
   \|\H^{\frac{1}{2}}(\EE_{\zeta}[\check\bbeta]-\bbeta_{OLS}) \|^2 \leq \|\TT_1\|^2+ \|\TT_2\|^2+2\|\TT_1\|\|\TT_2\|. \label{eq:H_expe_check_beta_beta_ols}
\end{align}
Consequently, to establish the bias bound for the debiased estimator in \Cref{theo:de_sam_ols}, it suffices to derive upper bounds on the norms $\|\TT_1\|$ and  $\|\TT_2\|$ in \eqref{eq:H_expe_check_beta_beta_ols}.
We first consider to bound the term $\|\TT_1\|$.
Without loss of generality, assume that the events $\zeta_1$, $\zeta_2$, $\zeta_3$ are independent of $\z_{s}$. Define $\zeta^{'} =\bigcap^3_{j=1} \zeta_j$  and $\delta_4=\Pr(\neg\zeta_4)$.
We first consider the first term $ \|\TT_1\|$. 
Letting $\check \gamma_{i_s}=1+\frac{1}{m}F_{i_si_s}\check \z_s^\top\check \Q_{-s}\check \z_s$ and  $\bar \gamma_{i}=1+\frac{1}{m\pi_i}F_{ii}\x_{\H_i}^\top\check \Q_{-s}\x_{\H_i}$ with $\x^\top_{\H_i}$ the $i^{\text{th}}$ row of $\X_{\H}$, along with the fact that the event $\zeta^{'}$ indicates $\check \Q_{-s}= \H^{1/2} \Q_{-s} \H^{1/2} \preceq 8\I_p$  and $\X_{\H}^\top\X_{\H}= \I_p$ %\Sachin{Why $\X_{\H}^\top\X_{\H} \preceq \I$}\Chengmei{This is a mistake; it should be $=$, not a $ \preceq$}
,  we obtain  
\begin{align*}
    \|\TT_1\|&=\left\|\EE_{\zeta}\left[\check \Q_{-s}\check \z_s \check a_s\left(\frac{F_{i_si_s}}{\check \gamma_{i_s}}-1\right)\right]\right\|=\max\limits_{\|\uu\|=1}\EE_{\zeta}\left[\uu^\top\check \Q_{-s}\check \z_s\check a_s\left(\frac{F_{i_si_s}}{\check \gamma_{i_s}}-1\right)\right]\\
    &\overset{(a)}{\leq}\max\limits_{\|\uu\|=1}\sqrt{\EE_{\zeta}\left[\uu^\top\check \Q_{-s}\check \z_s\check \z_s^\top\check \Q_{-s}\uu\left|\frac{F_{i_si_s}}{\check \gamma_{i_s}}-1\right|\right]\EE_{\zeta}\left[\check a_s^2\left|\frac{F_{i_si_s}}{\check \gamma_{i_s}}-1\right|\right]}\\
    &\overset{(b)}{\leq} \max\limits_{\|\uu\|=1}2\sqrt{\EE_{\zeta^{'}}\left[\uu^\top\check \Q_{-s}\check \z_s\check \z_s^\top\check \Q_{-s}\uu\left|\frac{F_{i_si_s}}{\check \gamma_{i_s}}-1\right|
    \right]\EE_{\zeta^{'}}\left[ \check a_s^2\left|\frac{F_{i_si_s}}{\check \gamma_{i_s}}-1\right|\right]}\\
     % &\leq \max\limits_{\|\uu\|=1}4\sqrt{\EE_{\zeta^{'}}\left[\uu^\top\check \Q_{-s}\check \z_s\check \z_s^\top\check \Q_{-s}\uu\left(\frac{F_{i_si_s}}{\check \gamma_{i_s}}-1\right)^2\right]\EE_{\zeta^{'}}[ \check a_s^2]}\\
  &  \leq \max\limits_{\|\uu\|=1}2\|\rr\|\sqrt{\EE_{\zeta^{'}}\left[\max\limits_{1\leq i\leq n}\left|\frac{F_{ii}}{\bar \gamma_{i}}-1\right|\uu^\top\check \Q_{-s}\X_{\H}^\top\X_{\H}\check \Q_{-s}\uu\right]\EE_{\zeta^{'}}\left[\max\limits_{1\leq i\leq n}\left|\frac{F_{ii}}{\bar \gamma_{i}}-1\right|\right]}\\
  & \leq 16\|\rr\|\EE_{\zeta^{'}}\left[\max\limits_{1\leq i\leq n}\left|\frac{F_{ii}}{\bar \gamma_{i}}-1\right|\right],
\end{align*}
where in $(a)$ we apply Cauchy-Schwarz inequality, in    $(b)$ we use  \Cref{lem:niu_proba_space}.
% and in $(c)$ we invoke  $m\pi_i\geq 2\ell_i(\X)$ under the condition  $m>2\theta_{\max,\X}p$.
Recalling the result \eqref{eq:check_Q_quadratic} in \Cref{lem:niu_auxiliary}, we then get
\begin{align*}
    \|\TT_1\|=O\left(\frac{\log n}{\loglog n}\sqrt{\frac{\theta^3_{\max}(\X)p^3}{m^3}}\right)\cdot  \|\rr\|.
\end{align*}
We next  explore the another term $\|\TT_2\|$. Using the fact $\X_{\H}^\top\rr=0$, we get
\begin{align*}
  \|\TT_2\|&=\| \EE_{\zeta}[\check\Q_{-s}\check\z_s \check a_s]\|=\frac{1}{1-\delta_4} (\|\EE_{\zeta^{'}}[\check\Q_{-s}\check\z_s \check a_s]-\EE_{\zeta^{'}}[\check\Q_{-s}\check\z_s \check a_s\cdot\mathbf{1}_{\neg\zeta_4} ]\|)\\
 &=\frac{1}{1-\delta_4}( \|\EE_{\zeta^{'}}[\check\Q_{-s}\X_{\H}^\top\rr]-\EE_{\zeta^{'}}[\check\Q_{-s}\check\z_s \check a_s\cdot\mathbf{1}_{\neg\zeta_4} ]\|)=-\frac{1}{1-\delta_4}\|\EE_{\zeta^{'}}[\check\Q_{-s}\check\z_s \check a_s\cdot\mathbf{1}_{\neg\zeta_4} ]\|\\
 &\leq 2 \|\EE_{\zeta^{'}}[\check\Q_{-s}\check\z_s \check a_s\cdot\mathbf{1}_{\neg\zeta_4} ]\|= \max\limits_{\|\uu\|=1}2 \EE_{\zeta^{'}}[\uu^\top\check\Q_{-s}\check\z_s \check a_s\cdot\mathbf{1}_{\neg\zeta_4} ]\\
 &\leq \max\limits_{\|\uu\|=1}2 \sqrt{\EE_{\zeta^{'}}[\uu^\top\check\Q_{-s}\check\z_s\check\z_s^\top  \check\Q_{-s}\uu\cdot\mathbf{1}_{\neg\zeta_4} ]\EE_{\zeta^{'}}[ \check a^2_s  ]}\\
 &\leq \max\limits_{\|\uu\|=1}2\|\rr\| \sqrt{\EE_{\zeta^{'}}[\uu^\top\check\Q_{-s}\check\z_s\check\z_s^\top\check\Q_{-s} \uu \cdot\mathbf{1}_{\neg\zeta_4} ]}.
\end{align*} 
Taking $\check \Q_{-s}= \H^{1/2} \Q_{-s} \H^{1/2} \preceq8\I_p$  again, we  derive 
%\Sachin{Again, I don't see why $\X_H^\top\X_H \preceq \I$}
\begin{align*}
    \EE_{\zeta^{'}}[\uu^\top\check\Q_{-s}\check\z_s\check\z_s^\top\check\Q_{-s} \uu ]= \sum^n_{i=1}  \EE_{\zeta^{'}}[\uu^\top\check\Q_{-s}\x_{\H_i}\x_{\H_i}^\top 
\check\Q_{-s} \uu ]= \EE_{\zeta^{'}}[\uu^\top\check\Q_{-s}\X_{\H}^\top\X_{\H} 
\check\Q_{-s} \uu ]\leq 64,
\end{align*}
and 
\begin{align*}
    \Var_{\zeta^{'}}[\uu^\top\check\Q_{-s}\check\z_s\check\z_s^\top\check\Q_{-s} \uu]&\leq  \EE_{\zeta^{'}}[(\uu^\top\check\Q_{-s}\check\z_s\check\z_s^\top\check\Q_{-s} \uu )^2]=\sum^n_{i=1}  \frac{1}{\pi_i} \EE_{\zeta^{'}}[(\uu^\top\check\Q_{-s}\x_{\H_i}\x_{\H_i}^\top 
 \check\Q_{-s}\uu )^2]\\
 &\leq 8^4\max\limits_{1\leq i\leq n} \frac{\ell_{i}(\X)}{\pi_i} \leq 8^4\theta_{\max}(\X)p.
\end{align*}
Applying  the Chebyshev's inequality, we obtain, for $z> 64$, 
\begin{align*}
    \Pr(\uu^\top\check\Q_{-s}\check\z_s\check\z_s^\top\check \Q_{-s} \uu \geq z~|~\zeta^{'})\leq\frac{8^4\theta_{\max}(\X)p}{z^2}. 
\end{align*}
This together with $\delta_4\leq m^{-3}$ further leads to 
\begin{align*}
&   \EE_{\zeta^{'}}[\uu^\top\check \Q_{-s}\check \z_s\check \z_s^\top\check \Q_{-s} \uu \cdot\mathbf{1}_{\neg\zeta_4} ]=\int^\infty_0 \Pr(\uu^\top\check \Q_{-s}\check \z_s\check \z_s^\top\check \Q_{-s} \uu \cdot \mathbf{1}_{\neg\zeta_4}\geq z~|~\zeta^{'})dz\\
 &=\int^{2m^2}_0 \Pr(\uu^\top\check \Q_{-s}\check \z_s\check \z_s^\top\check \Q_{-s} \uu \cdot \mathbf{1}_{\neg\zeta_4}\geq z~|~\zeta^{'})dz+\int^\infty_{2m^2} \Pr(\uu^\top\check \Q_{-s}\check \z_s\check \z_s^\top\check \Q_{-s} \uu \cdot \mathbf{1}_{\neg\zeta_4}\geq z~|~\zeta^{'})dz\\
 &\leq 2m^2\delta_4+\int^\infty_{2m^2} \Pr(\uu^\top\check \Q_{-s}\check \z_s\check \z_s^\top\check \Q_{-s} \uu \geq z~|~\zeta^{'})dz\\
 &\leq \frac{2}{m}+\frac{8^4\theta_{\max}(\X)p}{2m^2}= O\left(\frac{1}{m}+\frac{\theta_{\max}(\X)p}{m^2}\right).
\end{align*}
We thus get 
\begin{align*}
     \|\TT_2\|= O\left(\frac{1}{m}+\frac{\theta_{\max}(\X)p}{m^2}\right)\cdot \|\rr\|.
\end{align*}

\subsection{Detailed proof of the variance of debiased estimator  in  Theorem~\ref{theo:de_sam_ols}}
\label{subsec:detail_proof_ols_variance}
A heuristic derivation of the variance of the debiased estimator proceeds analogously to that of the bias, and is omitted here for clarity.

Building on the proof strategy for the bias of the debiased estimator, the analysis of its variance proceeds in two main steps:
\begin{enumerate}
  \item construct a high probability event $\zeta$ as in \eqref{eq:events}; and
  \item conditioning  on that event $\zeta$, bound the quantity  $\textbf{Var}(\check\bbeta)-\Delta(\X)$ using a ``leave-one-out'' type analysis.
\end{enumerate}

To complete the  analysis of the variance of the debiased estimator, 
we  begin by rewriting 
\begin{align*}
    & ~\textbf{Var}(\check\bbeta)-\Delta(\X)=\EE_{\zeta}[(\check\bbeta-\bbeta_{OLS})^\top \X^\top \X(\check\bbeta-\bbeta_{OLS})]-\Delta(\X)\\
    &=\EE_{\zeta}[\rr^\top\check\S^\top\check\S \X_{\H}\check\Q^2\X_{\H}^\top\check\S^\top\check\S\rr]-\Delta(\X) =\frac{1}{m^2}\sum^m_{s=1}\sum^{m}_{k=1}\EE_{\zeta}[\check\z_k^\top\check\Q^2\check\z_sF_{i_si_s}F_{i_ki_k}\check a_s\check a_k]-\Delta(\X)\\
     &=\frac{1}{m^2}\sum^m_{s=k=1}\EE_{\zeta}[\check\z_s^\top\check\Q^2\check\z_sF_{i_si_s}^2\check a_s^2]+\frac{1}{m^2}\sum^m_{s\neq k}\EE_{\zeta}[\check\z_k^\top\check\Q^2\check\z_sF_{i_si_s}F_{i_ki_k}\check a_s\check a_k]-\Delta(\X).
\end{align*}
Recalling \begin{align*}
\Delta(\X)= \rr^\top\diag\Big\{\frac{\ell_i(\X)}{m\pi_{i}}\Big\}^{n}_{i=1}\rr, 
\end{align*}
we then consider
\begin{align*}
 &  \frac{1}{m^2}\sum^m_{s=k=1}\EE_{\zeta}[\check\z_s^\top\check\Q^2\check\z_sF_{i_si_s}^2\check a_s^2]-\Delta(\X) =\sum^m_{s=k=1}\EE_{\zeta}\left[\frac{\check\z_s^\top\check\Q^2_{-s}\check\z_sF_{i_si_s}^2 \check a_s^2}{m^2\check\gamma_{i_s}^2}\right]-\Delta(\X) \\ 
 &=\EE_{\zeta}\left[\frac{\check\z_s^\top\check\Q^2_{-s}\check\z_s F_{i_si_s}^2\check a_s^2}{m\check\gamma_{i_s}^2}\right]-\Delta(\X) \\
 &=\underbrace{\frac{1}{m}\EE_{\zeta}\left[\check\z_s^\top\check\Q^2_{-s}\check\z_s \check a_s^2\left(\frac{F_{i_si_s}^2}{\check\gamma_{i_s}^2}-1\right)\right]}_{W_1}+\underbrace{\EE_{\zeta}\left[\frac{1}{m}\check\z_s^\top\check\Q^2_{-s}\check\z_s\check a_s^2-\rr^\top \diag\left\{\frac{\x_{\H_i}^\top\check\Q^2_{-s}\x_{\H_i}}{m\pi_i}\right\}^n_{i=1}\rr\right]}_{W_2}\\
 &+\underbrace{\rr^\top \diag\left\{\frac{\x_{\H_i}^\top\EE_{\zeta}[\check\Q^2_{-s}-\check\Q^2]\x_{\H_i}}{m\pi_i}\right\}^n_{i=1}\rr}_{W_3}+\underbrace{\rr^\top\diag\Big\{\frac{\x_{\H_i}^\top\EE_{\zeta}[\check\Q^2]\x_{\H_i}-\ell_i(\X)-\check\Delta_{\check\Q^2,i}}{m\pi_{i}}\Big\}^{n}_{i=1}\rr}_{W_4}\\
 &+\underbrace{\rr^\top\diag\Big\{\frac{\check\Delta_{\check\Q^2,i}}{m\pi_{i}}\Big\}^{n}_{i=1}\rr}_{W_5},
\end{align*}
where  $\check\gamma_{i_s}=1+F_{i_si_s}\check\z_{i_s}^\top\check\Q_{-s}\check\z_{i_s}$, $\check\Delta_{\check\Q^2,i}=\x_{\H_i}^\top\X_{\H}^\top\bar \F\X_{\H}\x_{\H_i}$, with diagonal matrix $\bar \F=\diag\{\bar F_{ii}\}^n_{i=1}$ and diagonal entries
\begin{equation*}
    \bar F_{ii}= \frac{\x_{\H_i}^\top\EE_{\zeta}[\check\Q^{2}]\x_{\H_i}}{m\pi_i}.
\end{equation*}

In the following, we  bound the term $W_1$.  Considering   $\F_{i_si_s}<2$ under the condition $m\geq 2\theta_{\max}(\X)p$,  together with \Cref{lem:niu_proba_space},
we get
\begin{align*}
    W_1&\leq  \frac{1}{m}\EE_{\zeta}\left[\check\z_s^\top\check\Q^2_{-s}\check\z_s \check a_s^2\left|\frac{F_{i_si_s}^2}{\check\gamma_{i_s}^2}-1\right|\right]\leq  \frac{2}{m}\EE_{\zeta^{'}}\left[\check\z_s^\top\check\Q^2_{-s}\check\z_s \check a_s^2\left|\frac{F_{i_si_s}^2}{\check\gamma_{i_s}^2}-1\right|\right]\\
    &\leq  \frac{2}{m}\EE_{\zeta^{'}}\left[\check\z_s^\top\check\Q^2_{-s}\check\z_s \check a_s^2\frac{F_{i_si_s}}{\check\gamma_{i_s}}\left|\frac{F_{i_si_s}}{\check\gamma_{i_s}}-1\right|\right]+\frac{2}{m}\EE_{\zeta^{'}}\left[\check\z_s^\top\check\Q^2_{-s}\check\z_s \check a_s^2\left|\frac{F_{i_si_s}}{\check\gamma_{i_s}}-1\right|\right]\\
     &\leq  \frac{4}{m}\EE_{\zeta^{'}}\left[\check\z_s^\top\check\Q^2_{-s}\check\z_s \check a_s^2\left|\frac{F_{i_si_s}}{\check\gamma_{i_s}}-1\right|\right]+\frac{2}{m}\EE_{\zeta^{'}}\left[\check\z_s^\top\check\Q^2_{-s}\check\z_s \check a_s^2\left|\frac{F_{i_si_s}}{\check\gamma_{i_s}}-1\right|\right]\\
     &=\frac{6}{m}\EE_{\zeta^{'}}\left[\check\z_s^\top\check\Q^2_{-s}\check\z_s \check a_s^2\left|\frac{F_{i_si_s}}{\check\gamma_{i_s}}-1\right|\right]=\sum^{n}_{i=1}\frac{6}{m\pi_i}\EE_{\zeta^{'}}\left[\x_{\H_i}^\top\check\Q^2_{-s}\x_{\H_i} r_i^2\left|\frac{ F_{ii}}{\bar\gamma_i }-1\right|\right],
\end{align*}
which combined  with $\check\Q_{-s}\preceq 8\I_p$ , $\x_{\H_i}^\top\check\Q_{-s}\x_{\H_i}/(m\pi_i)\leq 8\theta_{\max}(\X)p/m$,   and \Cref{lem:niu_auxiliary}  gives 
\begin{align*}
  W_1&\leq \sum^{n}_{i=1}\frac{384\ell_{i}(\X)}{m\pi_i}\EE_{\zeta^{'}}\left[ r_i^2\left|\frac{F_{ii}}{\bar\gamma_i }-1\right|\right]\leq  \frac{384\theta_{\max}(\X)p}{m} \EE_{\zeta^{'}}\left[\max\limits_{1\leq i\leq n}\left|\frac{ F_{ii}}{\bar\gamma_i }-1\right|\sum^{n}_{i=1}r_i^2\right]\\
  &\leq \frac{384\theta_{\max}(\X)p\|\rr\|^2}{m}\EE_{\zeta^{'}}\left[\max\limits_{1\leq i\leq n}\left|\frac{ F_{ii}}{\bar\gamma_i }-1\right|\right] = O\left(\frac{\log n}{\loglog n}\sqrt{\frac{\theta^3_{\max}(\X)p^3}{m^3}}\right)\cdot  \|\rr\|^2.
\end{align*}

Next, we explore the term $W_2$,
\begin{align*}
    W_2&=\frac{1}{1-\delta_4} \EE_{\zeta^{'}}\left[\frac{1}{m}\check\z_s^\top\check\Q^2_{-s}\check\z_s \check a_s^2-\rr^\top \diag\left\{\frac{\x_{\H_i}^\top\check\Q^2_{-s}\x_{\H_i}}{m\pi_i}\right\}^n_{i=1}\rr\right]\\
    &-\frac{1}{1-\delta_4} \EE_{\zeta^{'}}\left[\left(\frac{1}{m}\check\z_s^\top\check\Q^2_{-s}\check\z_s \check a_s^2-\rr^\top \diag\left\{\frac{\x_{\H_i}^\top\check\Q^2_{-s}\x_{\H_i}}{m\pi_i}\right\}^n_{i=1}\rr\right)\cdot\mathbf{1}_{\neg\zeta_4}\right ]\\
   & \leq 2\EE_{\zeta^{'}}\left[\frac{1}{m}\check\z_s^\top\check\Q^2_{-s}\check\z_s \check a_s^2\cdot\mathbf{1}_{\neg\zeta_4}\right ]+2\EE_{\zeta^{'}}\left[\rr^\top \diag\left\{\frac{\x_{\H_i}^\top\check\Q^2_{-s}\x_{\H_i}}{m\pi_i}\right\}^n_{i=1}\rr\cdot\mathbf{1}_{\neg\zeta_4}\right ]\\
   & \leq 2\EE_{\zeta^{'}}\left[\frac{1}{m}\check\z_s^\top\check\Q^2_{-s}\check\z_s \check a_s^2\cdot\mathbf{1}_{\neg\zeta_4}\right ]+\frac{128\theta_{\max}(\X)p\|\rr\|^2\delta_4}{m}\\
   &=2\EE_{\zeta^{'}}\left[\frac{1}{m}\check\z_s^\top\check\Q^2_{-s}\check\z_s \check a_s^2\cdot\mathbf{1}_{\neg\zeta_4}\right ]+O\left(\frac{\theta_{\max}(\X)p}{m^2}\right)\|\rr\|^2.
\end{align*}
Using $\check \Q_{-s}\preceq 8\I_p $, we obtain
\begin{align*}
   \EE_{\zeta^{'}}\left[\frac{1}{m}\check\z_s^\top\check\Q^2_{-s}\check\z_s \check a_s^2\right ]=\sum^n_{i=1}\frac{1}{m\pi_i}   \EE_{\zeta^{'}}\left[\x_{\H_i}^\top\check\Q^2_{-s}\x_{\H_i}  r_i^2\right ]\leq   \frac{64\theta_{\max}(\X)p}{m}\|\rr\|^2,
\end{align*}
and 
\begin{align*}
 \Var_{\zeta^{'}}\left[\frac{1}{m}\check\z_s^\top\check\Q^2_{-s}\check\z_s \check a_s^2\right ]& \leq \frac{1}{m^2}  \EE_{\zeta^{'}}\left[(\check\z_s^\top\check\Q^2_{-s}\check\z_s \check a_s^2)^2\right ]\leq \sum^n_{i=1}\frac{1}{m^2\pi_i^3}\EE_{\zeta^{'}}\left[(\x_{\H_i}^\top\check\Q^2_{-s}\x_{\H_i}  r_i^2)^2\right ] \\ 
 &\leq\frac{8^4\theta^3_{\max}(\X)p^3}{m^2} \sum^n_{i=1}\frac{ r_i^4}{\ell_i}\leq \frac{8^4\theta^3_{\max}(\X)p^3}{m^2\ell_{\min}(\X)}\|\rr\|^4,
\end{align*}
where we denote $\ell_{\min}(\X)=\min\limits_{1\leq i\leq m,\ell_i(\X)\neq 0}\ell_i(\X)$. %$\geq 1/\operatorname{poly}(n)$
Following the Chebyshev's inequality results in, for $z> 128\theta_{\max}(\X)p\|\rr\|^2/m $,  we get
\begin{align*}
    \Pr\left(\frac{1}{m}\check\z_s^\top\check\Q^2_{-s}\check\z_s \check a_s^2 \geq z~|~\zeta'\right)\leq \frac{8^4\theta^3_{\max}(\X)p^3}{m^2\ell_{\min}(\X)z^2}\|\rr\|^4. 
\end{align*}
Considering 
\begin{align*}
\ell_{\min}(\X)=
 \min\limits_{i,\|\x_i\|>0} \x_i^\top(\X^\top\X)^{-1}\x_i
\geq
\min\limits_{i,\|\x_i\|>0}\frac{\|\x_i\|^2}{\|\X\|^2} \geq
 \min\limits_{i,\|\x_i\|>0} \frac{\|\x_i\|^2}{\|\X\|_F^2} \geq \frac{1}{n^{\alpha}},
\end{align*}
for some constant $\alpha>0$, together with  $\delta_4\leq n^{-2-\alpha}$, 
we further derive 
\begin{align*}
  \EE_{\zeta^{'}}\left[\frac{1}{m}\check\z_s^\top\check\Q^2_{-s}\check\z_s \check a_s^2\cdot\mathbf{1}_{\neg\zeta_4}\right ]  & =\int^\infty_0 \Pr\left( \frac{1}{m}\check\z_s^\top\check\Q^2_{-s}\check\z_s \check a_s^2 \cdot \mathbf{1}_{\neg\zeta_4}\geq z~|~\zeta'\right)dz\\
 &=\int^{2m\|\rr\|^2/\ell_{\min}(\X)}_0 \Pr\left(\frac{1}{m}\check\z_s^\top\check\Q^2_{-s}\check\z_s \check a_s^2 \cdot \mathbf{1}_{\neg\zeta_4}\geq z~|~\zeta'\right)dz\\
 &+\int^\infty_{2m\|\rr\|^2/\ell_{\min}(\X)} \Pr\left(\frac{1}{m}\check\z_s^\top\check\Q^2_{-s}\check\z_s \check a_s^2 \cdot \mathbf{1}_{\neg\zeta_4}\geq z~|~\zeta'\right)dz\\
 &\leq \frac{2m\|\rr\|^2\delta_4}{\ell_{\min}(\X)}+\int^\infty_{2m\|\rr\|^2/\ell_{\min}(\X)} \Pr\left(\frac{1}{m}\check\z_s^\top\check\Q^2_{-s}\check\z_s \check a_s^2\geq z~|~\zeta'\right)dz \\
 &\leq  2mn^{\alpha}\delta_4\|\rr\|^2+\frac{8^4\theta^3_{\max}(\X)p^3\|\rr\|^2}{2m^3}=O\left(\frac{1}{m}+\frac{\theta^3_{\max}(\X)p^3}{m^3}\right)\cdot \|\rr\|^2.
\end{align*}
This further yields 
\begin{align*}
    W_2 =  O\left(\frac{1}{m}+\frac{\theta^3_{\max}(\X)p^3}{m^3}\right)\cdot \|\rr\|^2.
\end{align*}
Recalling the result \eqref{eq:check_Q_check_Q_s} in \Cref{lem:niu_auxiliary}, it follows that 
\begin{align*}
    W_3\leq \frac{\theta_{\max}(\X)p}{m}\|\rr\|^2 \|\EE_{\zeta}[\check\Q^2_{-s}-\check\Q^2]\| = O\left(\frac{\theta_{\max}(\X)p}{m^2}\right)\cdot \|\rr\|^2.
\end{align*}
Then, applying the result \eqref{eq:check_Q_second_moment} in \Cref{lem:niu_auxiliary} leads to 
\begin{align*}
    W_4&\leq \max\limits_{1\leq i\leq n}\frac{\|\rr\|^2\ell_i(\X) }{m\pi_i}\left\| \EE_\zeta[\check\Q^2]-\I_p-\X_{\H}^\top\bar \F\X_{\H}\right\| \\ 
    &=  O\left(\frac{\log n}{\loglog n}\sqrt{\frac{\theta^3_{\max}(\X)p^3}{m^3}}\cdot \frac{\theta_{\max}(\X)p}{m}\right)\cdot\|\rr\|^2\\
    &=  O\left(\frac{\log n}{\loglog n}\sqrt{\frac{\theta^5_{\max}(\X)p^5}{m^5}}\right)\cdot\|\rr\|^2.
\end{align*}
Then, we turn to bound the term  $W_5$. Taking $\check \Q \preceq 8\I_p $ and 
\begin{align*}
  \|   \X_{\H}^\top\bar \F\X_{\H}\|\leq  \max\limits_{1\leq i\leq n}\bar F_{ii} \leq \frac{64\theta_{\max}(\X)p}{m},
\end{align*}
it follows that 
\begin{align*}
   W_5\leq   \max\limits_{1\leq i\leq n}\frac{\ell_i(\X) }{m\pi_i}\|\rr\|^2\left\|\X_{\H}^\top\bar \F\X_{\H}\right\| = O\left(\frac{\theta^2_{\max}(\X)p^2}{m^2}\right)\cdot  \|\rr\|^2.
\end{align*}

Now, it remains to bound the cross term $m^{-2}\sum^m_{s\neq k}\EE_{\zeta}[\check\z_k^\top\check\Q^2\check\z_sF_{i_si_s}F_{i_ki_k}\check a_s\check a_k]$.  Without loss of generality, assume that the events $\zeta_1$, $\zeta_2$ are independent of both $\z_{s}$ and $\z_k$. Define $\zeta^{''} =\zeta_1\bigcap \zeta_2$, $\zeta^{'''} =\zeta_3\bigcap \zeta_4$  and $\delta_{\zeta^{'''}}=\Pr(\neg \zeta^{'''})$.
Let $\check\gamma_{i_{s,k}}=1+\frac{1}{m}F_{i_ki_k}\check\z_k^\top\check\Q_{-sk}\check\z_k$.
We first rewrite
\begin{align*}
  &  \frac{1}{m^2}\sum^m_{s\neq k}\EE_{\zeta}[\check\z_k^\top\check\Q^2\check\z_sF_{i_si_s}F_{i_ki_k}\check a_s\check a_k]=\sum^m_{s\neq k}\EE_{\zeta}\left[\frac{\check\z_k^\top\check\Q_{-s}^2\check\z_sF_{i_si_s}F_{i_ki_k}\check a_s\check a_k}{m^2\check\gamma_{i_s}}\right]\\
  &-\sum^m_{s\neq k}\EE_{\zeta}\left[\frac{\check\z_k^\top\check\Q_{-s}\check\z_s\check\z_s^\top\check\Q_{-s}^2\check\z_sF^2_{i_si_s}F_{i_ki_k}\check a_s\check a_k}{m^3\check\gamma_{i_s}^2}\right]=\sum^m_{s\neq k}\EE_{\zeta}\left[\frac{\check\z_k^\top\check\Q_{-sk}^2\check\z_sF_{i_si_s}F_{i_ki_k}\check a_s\check a_k}{m^2\check\gamma_{i_s}\check\gamma_{i_{s,k}}}\right]\\
    &-\sum^m_{s\neq k}\EE_{\zeta}\left[\frac{\check\z_k^\top\check\Q^2_{-sk}\check\z_k\check\z_k^\top\check\Q_{-sk}\check\z_sF_{i_si_s}F^2_{i_ki_k}\check a_s\check a_k}{m^3\check\gamma_{i_s}\check\gamma^2_{i_{s,k}}}\right]-\sum^m_{s\neq k}\EE_{\zeta}\left[\frac{\check\z_k^\top\check\Q_{-s}\check\z_s\check\z_s^\top\check\Q_{-s}^2\check\z_sF^2_{i_si_s}F_{i_ki_k}\check a_s\check a_k}{m^3\check\gamma_{i_s}^2}\right]\\
    &=\frac{m(m-1)}{m^2}\EE_{\zeta}\left[\frac{\check\z_k^\top\check\Q_{-sk}^2\check\z_sF_{i_si_s}F_{i_ki_k}\check a_s\check a_k}{\check\gamma_{i_s}\check\gamma_{i_{s,k}}}\right]\\
    &-\frac{m(m-1)}{m^2}\EE_{\zeta}\left[\frac{\check\z_k^\top\check\Q^2_{-sk}\check\z_k\check\z_k^\top\check\Q_{-sk}\check\z_sF_{i_si_s}F^2_{i_ki_k}\check a_s\check a_k}{m\check\gamma_{i_s}\check\gamma^2_{i_{s,k}}}\right]-\sum^m_{s\neq k}\EE_{\zeta}\left[\frac{\check\z_k^\top\check\Q_{-sk}\check\z_s\check\z_s^\top\check\Q_{-s}^2\check\z_sF^2_{i_si_s}F_{i_ki_k}\check a_s\check a_k}{m^3\check\gamma_{i_s}^2\check\gamma_{i_{s,k}}}\right]\\
    &=\frac{m(m-1)}{m^2(1-\delta_{\zeta^{'''}})}\left(\EE_{\zeta^{''}}\left[\frac{\check\z_k^\top\check\Q_{-sk}^2\check\z_sF_{i_si_s}F_{i_ki_k}\check a_s\check a_k}{\check\gamma_{i_s}\check\gamma_{i_{s,k}}}\right]-\EE_{\zeta^{''}}\left[\frac{\check\z_k^\top\check\Q_{-sk}^2\check\z_sF_{i_si_s}F_{i_ki_k}\check a_s\check a_k}{\check\gamma_{i_s}\check\gamma_{i_{s,k}}}\cdot\mathbf{1}_{\neg\zeta^{'''}}\right]\right)\\ 
    &-\frac{m(m-1)}{m^2(1-\delta_{\zeta^{'''}})}\left(\EE_{\zeta^{''}}\left[\frac{\check\z_k^\top\check\Q^2_{-sk}\check\z_k\check\z_k^\top\check\Q_{-sk}\check\z_sF_{i_si_s}F^2_{i_ki_k}\check a_s\check a_k}{m\check\gamma_{i_s}\check\gamma^2_{i_{s,k}}}\right]\right.\\
    &\left.-\EE_{\zeta^{''}}\left[\frac{\check\z_k^\top\check\Q^2_{-sk}\check\z_k\check\z_k^\top\check\Q_{-sk}\check\z_sF_{i_si_s}F^2_{i_ki_k}\check a_s\check a_k}{m\check\gamma_{i_s}\check\gamma^2_{i_{s,k}}}\cdot\mathbf{1}_{\neg\zeta^{'''}}\right]\right)\\
    &-\sum^m_{s\neq k}\EE_{\zeta}\left[\frac{\check\z_k^\top\check\Q_{-sk}\check\z_s\check\z_s^\top\check\Q_{-s}^2\check\z_sF^2_{i_si_s}F_{i_ki_k}\check a_s\check a_k}{m^3\check\gamma_{i_s}^2\check\gamma_{i_{s,k}}}\right]=  \underbrace{\frac{m(m-1)}{m^2(1-\delta_{\zeta^{'''}})}\EE_{\zeta^{''}}\left[\frac{\check\z_k^\top\check\Q_{-sk}^2\check\z_sF_{i_si_s}F_{i_ki_k}\check a_s\check a_k}{\check\gamma_{i_s}\check\gamma_{i_{s,k}}}\right]}_{W_6}\\
    &-\underbrace{\frac{m(m-1)}{m^2(1-\delta_{\zeta^{'''}})}\EE_{\zeta^{''}}\left[\frac{\check\z_k^\top\check\Q_{-sk}^2\check\z_sF_{i_si_s}F_{i_ki_k}\check a_s\check a_k}{\check\gamma_{i_s}\check\gamma_{i_{s,k}}}\cdot\mathbf{1}_{\neg\zeta^{'''}}\right]}_{W_7}\\ 
    &-\underbrace{\frac{m(m-1)}{m^2(1-\delta_{\zeta^{'''}})}\EE_{\zeta^{''}}\left[\frac{\check\z_k^\top\check\Q^2_{-sk}\check\z_k\check\z_k^\top\check\Q_{-sk}\check\z_sF_{i_si_s}F^2_{i_ki_k}\check a_s\check a_k}{m\check\gamma_{i_s}\check\gamma^2_{i_{s,k}}}\right]}_{W_8}\\
    &+\underbrace{\frac{m(m-1)}{m^2(1-\delta_{\zeta^{'''}})}\EE_{\zeta^{''}}\left[\frac{\check\z_k^\top\check\Q^2_{-sk}\check\z_k\check\z_k^\top\check\Q_{-sk}\check\z_sF_{i_si_s}F^2_{i_ki_k}\check a_s\check a_k}{m\check\gamma_{i_s}\check\gamma^2_{i_{s,k}}}\cdot\mathbf{1}_{\neg\zeta^{'''}}\right]}_{W_9}\\
    &-\underbrace{\sum^m_{s\neq k}\frac{1}{1-\delta_{\zeta^{'''}}}\EE_{\zeta^{''}}\left[\frac{\check\z_k^\top\check\Q_{-sk}\check\z_s\check\z_s^\top\check\Q_{-s}^2\check\z_sF^2_{i_si_s}F_{i_ki_k}\check a_s\check a_k}{m^3\check\gamma_{i_s}^2\check\gamma_{i_{s,k}}}\right]}_{W_{10}}\\
    &+  \underbrace{\sum^m_{s\neq k}\frac{1}{1-\delta_{\zeta^{'''}}}\EE_{\zeta^{''}}\left[\frac{\check\z_k^\top\check\Q_{-sk}\check\z_s\check\z_s^\top\check\Q_{-s}^2\check\z_sF^2_{i_si_s}F_{i_ki_k}\check a_s\check a_k}{m^3\check\gamma_{i_s}^2\check\gamma_{i_{s,k}}}\cdot\mathbf{1}_{\neg\zeta^{'''}}\right]}_{W_{11}}.
\end{align*}
Next, we first bound  the term $W_6$. Define  $\bar\gamma_{i_{s,k,j}} =1+\frac{1}{m\pi_j}F_{jj}\x_{\H_j}^\top\check\Q_{-sk}\x_{\H_j}$ and $\check\gamma_{i_{k,s}} =1+\frac{1}{m}F_{i_si_s}\check\z_s^\top\check\Q_{-sk}\check\z_s$. We rewrite
\begin{align*}
   W_6 &=\underbrace{\frac{m(m-1)}{m^2(1-\delta_{\zeta^{'''}})}\EE_{\zeta^{''}}\left[\frac{\check\z_k^\top\check\Q_{-sk}^2\check\z_sF_{i_si_s}F_{i_ki_k}\check a_s\check a_k}{\check\gamma_{i_{k,s}}\check\gamma_{i_{s,k}}}\right]}_{W_{61}}\\
   &+\underbrace{\frac{m(m-1)}{m^2(1-\delta_{\zeta^{'''}})}\EE_{\zeta^{''}}\left[\frac{\check\z_k^\top\check\Q_{-sk}^2\check\z_sF_{i_si_s}F_{i_ki_k}\check a_s\check a_k}{\check\gamma_{i_{s,k}}}\left(\frac{1}{\check\gamma_{i_s}}-\frac{1}{\check\gamma_{i_{k,s}}}\right)\right]}_{W_{62}}.
\end{align*}

For  the term $W_{61}$, we get 
\begin{align*}
    W_{61}&=\underbrace{\frac{m(m-1)}{m^2(1-\delta_{\zeta^{'''}})}\EE_{\zeta^{''}}\left[\frac{\check\z_k^\top\check\Q_{-sk}^2\check\z_sF_{i_ki_k}\check a_s\check a_k}{\check\gamma_{i_{s,k}}}\left(\frac{F_{i_si_s}}{\check\gamma_{i_{k,s}}}-1\right)\right]}_{W_{611}} \\
     &+\underbrace{\frac{m(m-1)}{m^2(1-\delta_{\zeta^{'''}})}\EE_{\zeta^{''}}\left[\check\z_k^\top\check\Q_{-sk}^2\check\z_s\check a_s\check a_k\left(\frac{F_{i_ki_k}}{\check\gamma_{i_{s,k}}}-1\right)\right]}_{W_{612}}\\
    &+\underbrace{\frac{m(m-1)}{m^2(1-\delta_{\zeta^{'''}})}\EE_{\zeta^{''}}\left[\check\z_k^\top\check\Q_{-sk}^2\check\z_s\check a_s\check a_k\right]}_{W_{613}}.
\end{align*}
 Denote $\bar\gamma_{i_{k,s,i}}=1+\frac{1}{m\pi_i}F_{ii}\x_{\H_i}^\top\check\Q_{-sk}\x_{\H_i}$. Together with  $\check\Q_{-sk}\preceq 8\I_p$ and $F_{jj}<2$,
we then obtain
\begin{align*}   W_{611}&=\sum^n_{j=1}\sum^n_{i=1}\frac{m(m-1)}{m^2(1-\delta_{\zeta^{'''}})}\EE_{\zeta^{''}}\left[\frac{\x_{\H_j}^\top\check\Q_{-sk}^2\x_{\H_i}F_{jj} r_i  r_j }{\bar\gamma_{i_{s,k,j}}}\left(\frac{F_{ii}}{\bar\gamma_{i_{k,s,i}}}-1\right)\right]\\
&=\frac{m(m-1)}{m^2(1-\delta_{\zeta^{'''}})}\EE_{\zeta^{''}}\left[\rr^\top\diag\left\{\frac{ F_{jj}}{\bar\gamma_{i_{s,k,j}}}\right\}^n_{j=1} \X_{\H}\check\Q_{-sk}^2\X_{\H}^\top\diag\left\{\left(\frac{F_{ii}}{\bar\gamma_{i_{k,s,i}}}-1\right)\right\}^n_{i=1}\rr\right]\\
&\leq256\EE_{\zeta^{''}}\left[\max\limits_{1\leq i\leq n}\left|\frac{F_{ii}}{\bar\gamma_{i_{k,s,i}}}-1\right|\right]\|\rr\|^2 = O\left(\frac{\log n}{\loglog n}\sqrt{\frac{\theta^3_{\max}(\X)p^3}{m^3}}\right)\|\rr\|^2.
\end{align*}
where the last line follows from \eqref{eq:check_Q_quadratic_two} in \Cref{lem:niu_auxiliary}. Using  $\X_{\H}^\top \rr=0$, 
we next consider the term $W_{612}$:
\begin{align*}   W_{612}&=\sum^n_{j=1}\sum^n_{i=1}\frac{m(m-1)}{m^2(1-\delta_{\zeta^{'''}})}\EE_{\zeta^{''}}\left[\x_{\H_j}^\top\check\Q_{-sk}^2\x_{\H_i} r_i  r_j\left(\frac{ F_{jj}}{\bar\gamma_{i_{s,k,j}}}-1\right)\right]\\
&=\frac{m(m-1)}{m^2(1-\delta_{\zeta^{'''}})}\EE_{\zeta^{''}}\left[\rr^\top\diag\left\{\left(\frac{ F_{jj}}{\bar\gamma_{i_{s,k,j}}}-1\right)\right\}^n_{j=1} \X_{\H}\check\Q_{-sk}^2\X_{\H}^\top\rr\right]=0.
\end{align*}
Recalling $\X_{\H}^\top \rr=0$ again, we now consider $W_{613}$:
\begin{align*}
 W_{613}=\frac{m(m-1)}{m^2(1-\delta_{\zeta^{'''}})}\sum^n_{j=1}\sum^n_{i=1}\EE_{\zeta^{''}}\left[\x_{\H_j}^\top\check\Q_{-sk}^2\x_{\H_i} r_i r_j\right]=\frac{m(m-1)}{m^2(1-\delta_{\zeta^{'''}})}\EE_{\zeta^{''}}\left[\rr^\top \X_{\H}\check\Q_{-sk}^2\X_{\H}^\top\rr\right]=0.
\end{align*}

Now, we turn  to explore the  term $W_{62}$. Noting $1<F_{ii}<2$, $\check\Q_{-sk}\preceq8\I_p$, and 
\begin{align*}
\check\gamma_{i_s}=1+  \frac{F_{i_si_s}\check\z_s^\top\check\Q_{-sk}\check\z_s  }{m}-\frac{F_{i_si_s}F_{i_ki_k}\check\z_s^\top\check\Q_{-sk}\check\z_k\check\z_k^\top\check\Q_{-sk}\check\z_s  }{m^2\check\gamma_{i_{s,k}}} \geq  1/2,
\end{align*}
under $m\geq 256\theta_{\max}p$,
we get
\begin{align}
 W_{62}&=   \frac{m(m-1)}{m^2(1-\delta_{\zeta^{'''}})}\EE_{\zeta^{''}}\left[\frac{\check\z_k^\top\check\Q_{-sk}^2\check\z_s\check\z_s^\top\check\Q_{-sk}\check\z_k\check\z_k^\top\check\Q_{-sk}\check\z_sF^2_{i_si_s}F^2_{i_ki_k}\check a_s\check a_k}{m^2\check\gamma^2_{i_{s,k}}\check\gamma_{i_{k,s}}\check\gamma_{i_{s}}}\right]\nonumber\\
 &= \frac{m(m-1)}{m^2(1-\delta_{\zeta^{'''}})}\sum^n_{i=1}\EE_{\zeta^{''}}\left[\frac{\check\z_k^\top\check\Q_{-sk}^2\x_{\H_i}\x_{\H_i}^\top\check\Q_{-sk}\check\z_k\check\z_k^\top\check\Q_{-sk}\x_{\H_i}F^2_{ii}F^2_{i_ki_k}r_i\check a_k}{m^2\bar\gamma_{i_{k,s,i}}\check\gamma^2_{i_{s,k}}\bar\gamma_{i}\pi_i}\right]
 \nonumber\\
 &\overset{(a)}{\leq}\sum^n_{i=1}\EE_{\zeta^{''}}\left[\frac{\x^\top_{\H_i}\check\Q_{-sk}^2\check\z_k\check\z_k^\top\check\Q_{-sk}^2\x_{\H_i}(\x_{\H_i}^\top\check\Q_{-sk}\check\z_k\check\z_k^\top\check\Q_{-sk}\x_{\H_i})^{1/2}F^2_{ii}F^2_{i_ki_k} r^2_i}{m^2\bar\gamma_{i_{k,s,i}}\check\gamma^2_{i_{s,k}}\bar\gamma_{i}\pi^{3/2}_i}\right]\nonumber\\
 &+\sum^n_{i=1}\EE_{\zeta^{''}}\left[\frac{F^2_{ii}F^2_{i_ki_k}\check a^2_k(\x_{\H_i}^\top\check\Q_{-sk}\check\z_k\check\z_k^\top\check\Q_{-sk}\x_{\H_i})^{3/2}}{m^2\bar\gamma_{i_{k,s,i}}\check\gamma^2_{i_{s,k}}\bar\gamma_{i}\pi^{1/2}_i}\right]\nonumber\\
 &\leq \frac{128\theta_{\max}(\X)p}{m}\sum^n_{i=1}\EE_{\zeta^{''}}\left[\frac{\x^\top_{\H_i}\check\Q_{-sk}^2\check\z_k\check\z_k^\top\check\Q_{-sk}^2\x_{\H_i} r^2_i}{m\pi_i}\right]\nonumber \\ 
 &+\frac{128\theta_{\max}(\X)p}{m}\sum^n_{i=1}\EE_{\zeta^{''}}\left[\frac{\check a^2_k\check\z_k^\top\check\Q_{-sk}\x_{\H_i}\x_{\H_i}^\top\check\Q_{-sk}\check\z_k}{m}\right]\nonumber\\
 &\leq \frac{128\theta_{\max}(\X)p}{m}\sum^n_{i=1}\sum^n_{j=1}\EE_{\zeta^{''}}\left[\frac{\x^\top_{\H_i}\check\Q_{-sk}^2\x_{\H_j}\x_{\H_j}^\top\check\Q_{-sk}^2\x_{\H_i} r^2_i}{m\pi_i}\right] \nonumber \\ 
 &+\frac{128\theta_{\max}(\X)p}{m}\EE_{\zeta^{''}}\left[\frac{\check a^2_k\check\z_k^\top\check\Q_{-sk}\X_{\H}^\top\X_{\H}\check\Q_{-sk}\check\z_k}{m}\right]\nonumber\\
 &\leq \frac{128 \theta_{\max}(\X)p}{m}\sum^n_{i=1}\EE_{\zeta^{''}}\left[\frac{\x^\top_{\H_i}\check\Q_{-sk}^2\X_{\H}^\top\X_{\H}\check\Q_{-sk}^2\x_{\H_i}r^2_i}{m\pi_i}\right]+\frac{8192\theta^2_{\max}(\X)p^2}{m^2}\EE_{\zeta^{''}}\left[\check a^2_k\right]\nonumber\\
 &\leq \frac{2\times 8^6
 \theta^2_{\max}(\X)p^2}{m^2}\sum^n_{i=1}r^2_i+\frac{8192\theta^2_{\max}(\X)p^2}{m^2}\sum^n_{j=1}r^2_j = O\left(\frac{\theta^2_{\max}(\X)p^2}{m^2}\right)\cdot\|\rr\|^2.\label{eq_W_62}
\end{align}
where in $(a)$ we use the inequality $2xy\leq x^2+y^2 $. Thus, we conclude 
\begin{align*}
 W_6 =  O\left(\frac{\log n}{\loglog n}\sqrt{\frac{\theta^3_{\max}(\X)p^3}{m^3}}\right)\|\rr\|^2.
\end{align*}

Next, we proceed to analyze the term $W_7$. Applying
\begin{align}
   & \EE_{\zeta^{''}}\left[\left|\frac{\check\z_k^\top\check\Q_{-sk}^2\check\z_sF_{i_si_s}F_{i_ki_k}\check a_s\check a_k}{\check\gamma_{i_s}\check\gamma_{i_{s,k}}}\right|\right]\leq  4\EE_{\zeta^{''}}\left[|\check\z_k^\top\check\Q_{-sk}^2\check\z_s\check a_s\check a_k|\right]\nonumber\\
   &\leq 4\sqrt{\EE_{\zeta^{''}}\left[\check\z_k^\top\check\Q_{-sk}^2\check\z_s\check\z_s^\top\check\Q_{-sk}^2\check\z_k\right]\EE_{\zeta^{''}}\left[\check  a^2_s\check a_k^2\right]}=4\|\rr\|^2 \sqrt{\EE_{\zeta^{''}}\left[\check\z_k^\top\check\Q^2_{-sk}\X_{\H}^\top\X_{\H}\check\Q^2_{-sk}\check\z_k\right]}\nonumber\\
   &\leq4\|\rr\|^2\sqrt{\sum^n_{j=1}\EE_{\zeta^{''}}\left[\x_{\H_j}^\top\check\Q^2_{-sk}\X_{\H}^\top\X_{\H}\check\Q^2_{-sk}\x_{\H_j}\right]}\leq 256p^{1/2}\|\rr\|^2,\label{eq:exp_zq2zaa}
\end{align}
and 
\begin{align}
   & \Var_{\zeta^{''}}\left[\left(\frac{\check\z_k^\top\check\Q_{-sk}^2\check\z_sF_{i_si_s}F_{i_ki_k}\check a_s\check a_k}{\check\gamma_{i_s}\check\gamma_{i_{s,k}}}\right)^2\right]\leq \sum^n_{j=1}\sum^n_{i=1}\EE_{\zeta^{''}}\left[\frac{(\x_{\H_j}^\top\check\Q_{-sk}^2\x_{\H_i} )^2F_{ii}^2F_{jj}^2r^2_i r^2_j }{\pi_j\pi_i\bar \gamma^2_{i}\bar\gamma_{i_{s,k,j}}^2}\right]\nonumber\\
   &\leq \sum^n_{j=1}\sum^n_{i=1}\frac{16}{\pi_j\pi_i}\EE_{\zeta^{''}}\left[(\x_{\H_j}^\top\check\Q_{-sk}^2\x_{\H_i} )^2 r^2_i  r^2_j\right]\leq 2\times 8^5\theta_{\max}^2(\X)p^2\|\rr\|^4,\label{eq:var_zq2zaa}
\end{align}
together with the  Chebyshev's inequality, it follows that    for $z> 512p^{1/2}\|\rr\|^2 $, 
\begin{align*}
    \Pr\left(\left|\frac{\check\z_k^\top\check\Q_{-sk}^2\check\z_s F_{i_si_s}F_{i_ki_k}\check a_s\check a_k}{\check\gamma_{i_s}\check\gamma_{i_{s,k}}}\right|\geq z~|~\zeta^{''}\right)\leq\frac{2\times 8^5\theta^2_{\max}(\X)p^2\|\rr\|^4}{z^2}. 
\end{align*}
With $\delta_{\zeta^{''}}\leq m^{-3}$, we then deduce
\begin{align*}
   &  \EE_{\zeta^{''}}\left[\frac{\check\z_k^\top\check\Q_{-sk}^2\check\z_sF_{i_si_s}F_{i_ki_k}\check a_s\check a_k}{\check\gamma_{i_s}\check\gamma_{i_{s,k}}}\cdot\mathbf{1}_{\neg\zeta^{'''}}\right ]\leq \EE_{\zeta^{''}}\left[\left|\frac{\check\z_k^\top\check\Q_{-sk}^2\check\z_sF_{i_si_s}F_{i_ki_k}\check a_s\check a_k}{\check\gamma_{i_s}\check\gamma_{i_{s,k}}}\right|\cdot\mathbf{1}_{\neg\zeta^{'''}}\right ]\\
   &=\int^\infty_0 \Pr\left(\left| \frac{\check\z_k^\top\check\Q_{-sk}^2\check\z_sF_{i_si_s}F_{i_ki_k}\check a_s\check a_k}{\check\gamma_{i_s}\check\gamma_{i_{s,k}}}\right|\cdot\mathbf{1}_{\neg\zeta^{'''}}\geq z~|~\zeta^{''}\right)dz\\
 &=\int^{2m^2\|\rr\|^2}_0 \Pr\left(\left|\frac{\check\z_k^\top\check\Q_{-sk}^2\check\z_sF_{i_si_s}F_{i_ki_k}\check a_s\check a_k}{\check\gamma_{i_s}\check\gamma_{i_{s,k}}}\right|\cdot\mathbf{1}_{\neg\zeta^{'''}}\geq z~|~\zeta^{''}\right)dz\\
 &+\int^\infty_{2m^2\|\rr\|^2} \Pr\left(\left|\frac{\check\z_k^\top\check\Q_{-sk}^2\check\z_sF_{i_si_s}F_{i_ki_k}\check a_s\check a_k}{\check\gamma_{i_s}\check\gamma_{i_{s,k}}}\right|\cdot\mathbf{1}_{\neg\zeta^{'''}} \geq z~|~\zeta^{''}\right)dz\\
 &\leq \frac{2\|\rr\|^2}{m}+\int^\infty_{2m^2\|\rr\|^2} \Pr\left(\left|\frac{\check\z_k^\top\check\Q_{-sk}^2\check\z_sF_{i_si_s}F_{i_ki_k}\check a_s\check a_k}{\check\gamma_{i_s}\check\gamma_{i_{s,k}}}\right|\cdot\mathbf{1}_{\neg\zeta^{'''}}\geq z~|~\zeta^{''}\right)dz \\
 &\leq  \frac{2\|\rr\|^2}{m}+\frac{2\times 8^5\theta^2_{\max}(\X)p^2\|\rr\|^2}{2m^2}=O\left(\frac{1}{m}+\frac{\theta^2_{\max}(\X)p^2}{m^2}\right)\cdot\|\rr\|^2.
\end{align*}
This thus gives
\begin{align*}
    W_{7} = O\left(\frac{1}{m}+\frac{\theta^2_{\max}(\X)p^2}{m^2}\right)\cdot\|\rr\|^2.
\end{align*}
Subsequently, we continue to bound the term $W_{8}$. We  get
\begin{align*}
  W_{8}&=  \underbrace{\frac{m(m-1)}{m^2(1-\delta_{\zeta^{'''}})}\EE_{\zeta^{''}}\left[\frac{\check\z_k^\top\check\Q^2_{-sk}\check\z_k\check\z_k^\top\check\Q_{-sk}\check\z_sF_{i_si_s}F^2_{i_ki_k}\check a_s\check a_k}{m\check\gamma_{i_{k,s}}\check\gamma^2_{i_{s,k}}}\right]}_{W_{81}}\\
  &+  \underbrace{\frac{m(m-1)}{m^2(1-\delta_{\zeta^{'''}})}\EE_{\zeta^{''}}\left[\frac{\check\z_k^\top\check\Q^2_{-sk}\check\z_k\check\z_k^\top\check\Q_{-sk}\check\z_sF_{i_si_s}F^2_{i_ki_k}\check a_s\check a_k}{m\check\gamma^2_{i_{s,k}}}\left(\frac{1}{\check\gamma_{i_s}}-\frac{1}{\check\gamma_{i_{k,s}}}\right)\right]}_{W_{82}}.
\end{align*}
In  the following, we consider the first term $W_{81}$:
\begin{align*}
    W_{81} &=\underbrace{\frac{m(m-1)}{m^2(1-\delta_{\zeta^{'''}})}\EE_{\zeta^{''}}\left[\frac{\check\z_k^\top\check\Q^2_{-sk}\check\z_k\check\z_k^\top\check\Q_{-sk}\check\z_s F^2_{i_ki_k}\check a_s\check a_k}{m\check\gamma^2_{i_{s,k}}}\left(\frac{F_{i_si_s}}{\check\gamma_{i_{k,s}}}-1\right)\right]}_{W_{811}}\\
    &+\underbrace{\frac{m(m-1)}{m^3(1-\delta_{\zeta^{'''}})}\EE_{\zeta^{''}}\left[\check\z_k^\top\check\Q^2_{-sk}\check\z_k\check\z_k^\top\check\Q_{-sk}\check\z_s\check a_s\check a_k\left(\frac{ F^2_{i_ki_k}}{\check\gamma^2_{i_{s,k}}}-1\right)\right]}_{W_{812}}\\
    &+\underbrace{\frac{m(m-1)}{m^3(1-\delta_{\zeta^{'''}})}\EE_{\zeta^{''}}\left[\check\z_k^\top\check\Q^2_{-sk}\check\z_k\check\z_k^\top\check\Q_{-sk}\check\z_s\check a_s\check a_k\right]}_{W_{813}}.
\end{align*}
Recalling  $F_{ii}\leq 2$ and $\check\Q_{-sk}\preceq 8\I_p$ again, together with the result \eqref{eq:check_Q_quadratic_two} in \Cref{lem:niu_auxiliary}, we get 
\begin{align*}   W_{811}&=\sum^n_{j=1}\sum^n_{i=1}\frac{m(m-1)}{m^2(1-\delta_{\zeta^{'''}})}\EE_{\zeta^{''}}\left[\frac{\x_{\H_j}^\top\check\Q^2_{-sk}\x_{\H_j}\x_{\H_j}^\top\check\Q_{-sk}\x_{\H_i}F^2_{jj}r_i r_j}{m\pi_j\bar\gamma_{i_{s,k,j}}^2}\left(\frac{ F_{ii}}{\bar\gamma_{i_{k,s,i}}}-1\right)\right]\\
&=\sum^n_{j=1}\sum^n_{i=1}\frac{m(m-1)}{m^2(1-\delta_{\zeta^{'''}})}\EE_{\zeta^{''}}\left[\frac{\x_{\H_j}^\top\check\Q^2_{-sk}\x_{\H_j}\x_{\H_j}^\top\check\Q_{-sk}\x_{\H_i} F_{jj}^2r_ir_j }{m\pi_j\bar\gamma_{i_{s,k,j}}^2}\left(\frac{ F_{ii}}{\bar\gamma_{i_{k,s,i}}}-1\right)\right]\\
&=\frac{m(m-1)}{m^2(1-\delta_{\zeta^{'''}})}\EE_{\zeta^{''}}\left[\rr^\top\diag\left\{\frac{  F_{jj}^2\x_{\H_j}^\top\check\Q^2_{-sk}\x_{\H_j}}{m\pi_j\bar\gamma_{i_{s,k,j}}^2}\right\}^n_{j=1} \X_{\H}\check\Q_{-sk}\X_{\H}^\top\diag\left\{\frac{ F_{ii}}{\bar\gamma_{i_{k,s,i}}}-1\right\}^n_{i=1}\rr\right]\\
&\leq  \frac{4096\theta_{\max}(\X)p}{m}\EE_{\zeta^{''}}\left[\max\limits_{1\leq i\leq n}\left|\frac{F_{ii}}{\bar\gamma_{i_{k,s,i}}}-1\right|\right]\cdot\|\rr\|^2 = O\left(\frac{\log n}{\loglog n}\sqrt{\frac{\theta^5_{\max}(\X)p^5}{m^5}}\right)\cdot\|\rr\|^2.
\end{align*}
Now, we turn to analyze the term $W_{812}$. Applying $\X_{\H}^\top \rr=0$  again leads to
\begin{align*} 
W_{812}=\sum^n_{j=1}\sum^n_{i=1}\frac{m(m-1)}{m^2(1-\delta_{\zeta^{'''}})}\EE_{\zeta^{''}}\left[\frac{ \x_{\H_j}^\top\check\Q^2_{-sk}\x_{\H_j}\x_{\H_j}^\top\check\Q_{-sk}\x_{\H_i}r_i r_j}{m\pi_j}\left(\frac{ F_{jj}^2}{\bar\gamma_{i_{s,k,j}}^2}-1\right)\right]=0.
\end{align*}
 Similarly,  we get
\begin{align*}
 W_{813}&= \frac{m(m-1)}{m^2(1-\delta_{\zeta^{'''}})}\sum^n_{j=1}\sum^n_{i=1}\frac{1}{m\pi_j}\EE_{\zeta^{''}}\left[\x_{\H_j}^\top\check\Q^2_{-sk}\x_{\H_j}\x_{\H_j}^\top\check\Q_{-sk}\x_{\H_i} r_i r_j\right]=0.
\end{align*}
We then consider the term $W_{82}$, with $F_{ii}<2$ and $\check\Q_{-sk}\preceq 8\I_p$,
 \begin{align}
  W_{82}&\leq  \frac{m(m-1)}{m^2(1-\delta_{\zeta^{'''}})}\EE_{\zeta^{''}}\left[\frac{\check\z_k^\top\check\Q^2_{-sk}\check\z_k\check\z_k^\top\check\Q_{-sk}\check\z_s\check\z_s^\top\check\Q_{-sk}\check\z_k\check\z_k^\top\check\Q_{-sk}\check\z_sF^2_{i_si_s}F^3_{i_ki_k}\check a_s\check a_k}{m^3\check\gamma_{i_s}\check\gamma_{i_{k,s}}\check\gamma^3_{i_{s,k}}}\right]\nonumber\\
 &= \frac{m(m-1)}{m^2(1-\delta_{\zeta^{'''}})}\sum^n_{i=1}\EE_{\zeta^{''}}\left[\frac{\check\z_k^\top\check\Q^2_{-sk}\check\z_k\check\z_k^\top\check\Q_{-sk}\x_{\H_i}\x_{\H_i}^\top\check\Q_{-sk}\check\z_k\check\z_k^\top\check\Q_{-sk}\x_{\H_i}F^2_{ii}F^3_{i_ki_k}r_i\check a_k}{m^3\bar\gamma_{i_{k,s,i}}\check\gamma^3_{i_{s,k}}\bar\gamma_{i}\pi_i}\right]
\nonumber \\
 &\leq\sum^n_{i=1}\EE_{\zeta^{''}}\left[\frac{\check\z_k^\top\check\Q^2_{-sk}\check\z_k\x_{\H_i}^\top\check\Q_{-sk}\check\z_k\check\z_k^\top\check\Q_{-sk}\x_{\H_i}(\x_{\H_i}^\top\check\Q_{-sk}\check\z_k\check\z_k^\top\check\Q_{-sk}\x_{\H_i})^{1/2}F^2_{ii}F^3_{i_ki_k} r^2_i}{m^3\bar\gamma_{i_{k,s,i}}\check\gamma^3_{i_{s,k}}\bar\gamma_{i}\pi^{3/2}_i}\right]\nonumber\\
 &+ \sum^n_{i=1}\EE_{\zeta^{''}}\left[\frac{F^2_{ii}F^3_{i_ki_k}\check a^2_k\check\z_k^\top\check\Q^2_{-sk}\check\z_k(\x_{\H_i}^\top\check\Q_{-sk}\check\z_k\check\z_k^\top\check\Q_{-sk}\x_{\H_i})^{3/2}}{m^3\bar\gamma_{i_{k,s,i}}\check\gamma^3_{i_{s,k}}\bar\gamma_{i}\pi^{1/2}_i}\right]\nonumber\\
 &\leq \frac{16384\theta^2_{\max}(\X)p^2}{m^2}\sum^n_{i=1}\EE_{\zeta^{''}}\left[\frac{\x_{\H_i}\check\Q_{-sk}\check\z_k\check\z_k^\top\check\Q_{-sk}\x_{\H_i} r^2_i}{m\pi_i}\right]\nonumber\\
 &+\frac{16384\theta^2_{\max}(\X)p^2}{m^2}\sum^n_{i=1}\EE_{\zeta^{''}}\left[\frac{\check a^2_k\check\z_k^\top\check\Q_{-sk}\x_{\H_i}\x_{\H_i}^\top\check\Q_{-sk}\check\z_k}{m}\right]\nonumber\\
 &= O\left(\frac{\theta^3_{\max}(\X)p^3}{m^3}\right)\cdot\|\rr\|^2.\label{eq:w_82}
\end{align}
We now proceed to bound the term $W_9$. Similar to \eqref{eq:exp_zq2zaa} and \eqref{eq:var_zq2zaa}, 
we  have 
\begin{align*}
   & \EE_{\zeta^{''}}\left[\left|\frac{\check\z_k^\top\check\Q^2_{-sk}\check\z_k\check\z_k^\top\check\Q_{-sk}\check\z_sF_{i_si_s}F^2_{i_ki_k}\check a_s\check a_k}{m\check\gamma_{i_s}\check\gamma^2_{i_{s,k}}}\right|\right]\leq  \frac{512\theta_{\max}(\X)p}{m} \EE_{\zeta^{''}}\left[|\check\z_k^\top\check\Q_{-sk}\check\z_s\check a_s\check a_k|\right]\\
   %%%
   &\leq\frac{512\theta_{\max}(\X)p}{m} \sqrt{\EE_{\zeta^{''}}\left[\check\z_k^\top\check\Q_{-sk}\check\z_s\check\z_s\check\Q_{-sk}\check\z_k\right]\EE_{\zeta^{''}}\left[\check  a^2_s\check a^2_k\right]} \\ 
   &=\frac{512\theta_{\max}(\X)p\|\rr\|^2}{m} \sqrt{\EE_{\zeta^{''}}\left[\check\z_k^\top\check\Q_{-sk}\X_{\H}^\top\X_{\H}\check\Q_{-sk}\check\z_k\right]}\\
   %%%
   &\leq \frac{512\theta_{\max}(\X)p\|\rr\|^2}{m} \sqrt{\sum^n_{j=1}\EE_{\zeta^{''}}\left[\x_{\H_j}^\top\check\Q_{-sk}\X_{\H}^\top\X_{\H}\check\Q_{-sk}\x_{\H_j}\right]}\leq \frac{4096\theta_{\max}(\X)p^{3/2}\|\rr\|^2}{m},
\end{align*}
and 
\begin{align*}
   & \Var_{\zeta^{''}}\left[\left(\frac{\check\z_k^\top\check\Q^2_{-sk}\check\z_k\check\z_k^\top\check\Q_{-sk}\check\z_sF_{i_si_s}F^2_{i_ki_k}\check a_s\check a_k}{m\check\gamma_{i_s}\check\gamma^2_{i_{s,k}}}\right)^2\right]\\ 
   %%%
   &\leq \sum^n_{j=1}\sum^n_{i=1}\EE_{\zeta^{''}}\left[\frac{(\x_{\H_j}^\top\check\Q^2_{-sk}\x_{\H_j}\x_{\H_j}^\top\Q_{-sk}\x_{\H_i} )^2F^2_{ii}F^4_{jj} r^2_ir^2_j}{m^2\pi^3_j\pi_i\bar \gamma^2_{i}\bar\gamma_{i_{s,k,j}}^4}\right]\\
   %%%
    &\leq 64\sum^n_{j=1}\sum^n_{i=1}\frac{1}{m^2\pi^3_j\pi_i}\EE_{\zeta^{''}}\left[(\x_{\H_j}^\top\check\Q^2_{-sk}\x_{\H_j}\x_{\H_j}^\top\check \Q_{-sk}\x_{\H_i} )^2 r^2_i r^2_j\right]\leq \frac{8^8\theta_{\max}^4(\X)p^4}{m^2} \|\rr\|^4,
\end{align*}
Similarly to the bound of $W_7$,  we obtain
\begin{align*}
    W_9= O\left(\frac{1}{m}+\frac{\theta_{\max}^4(\X)p^4}{m^4}\right) \cdot\|\rr\|^2.
\end{align*}
In the following, we turn to explore the  term $W_{10}$. We first rewrite
\begin{align*}
    W_{10}&=\sum^m_{s\neq k}\frac{1}{1-\delta_{\zeta^{'''}}}\EE_{\zeta^{''}}\left[\frac{\check\z_k^\top\check\Q_{-sk}\check\z_s\check\z_s^\top\check\Q_{-s}^2\check\z_sF^2_{i_si_s}F_{i_ki_k}\check a_s\check a_k}{m^3\check\gamma_{i_s}^2\check\gamma_{i_{s,k}}}\right]\\
    &=\frac{m(m-1)}{m^2(1-\delta_{\zeta^{'''}})}\EE_{\zeta^{''}}\left[\frac{\check\z_k^\top\check\Q_{-sk}\check\z_s\check\z_s^\top\check\Q_{-s}^2\check\z_sF^2_{i_si_s}F_{i_ki_k}\check a_s\check a_k}{m\check\gamma_{i_s}^2\check\gamma_{i_{s,k}}}\right]\\
    &=\underbrace{\frac{m(m-1)}{m^2(1-\delta_{\zeta^{'''}})}\EE_{\zeta^{''}}\left[\frac{\check\z_k^\top\check\Q_{-sk}\check\z_s\check\z_s^\top\check\Q_{-s}^2\check\z_sF^2_{i_si_s}F_{i_ki_k}\check a_s\check a_k}{m\check\gamma^2_{i_{k,s}}\check\gamma_{i_{s,k}}}\right]}_{W_{101}}\\
    &+\underbrace{\frac{m(m-1)}{m^2(1-\delta_{\zeta^{'''}})}\EE_{\zeta^{''}}\left[\frac{\check\z_k^\top\check\Q_{-sk}\check\z_s\check\z_s^\top\check\Q_{-s}^2\check\z_sF^2_{i_si_s}F_{i_ki_k}\check a_s\check a_k}{m\check\gamma_{i_{s,k}}}\left(\frac{1}{\check\gamma^2_{i_{s}}}-\frac{1}{\check\gamma^2_{i_{k,s}}}\right)\right]}_{W_{102}}.
\end{align*}
Then, we write the term $W_{101}$:
\begin{align*}
 W_{101}&=\underbrace{\frac{m(m-1)}{m^2(1-\delta_{\zeta^{'''}})}\EE_{\zeta^{''}}\left[\frac{\check\z_k^\top\check\Q_{-sk}\check\z_s\check\z_s^\top\check\Q_{-s}^2\check\z_sF_{i_ki_k}\check a_s\check a_k}{m\check\gamma_{i_{s,k}}}\left(\frac{F^2_{i_si_s}}{\check\gamma^2_{i_{k,s}}}-1\right)\right]}_{W_{1011}}\\
 %%%
 &+\underbrace{\frac{m(m-1)}{m^2(1-\delta_{\zeta^{'''}})}\EE_{\zeta^{''}}\left[\frac{\check\z_k^\top\check\Q_{-sk}\check\z_s\check\z_s^\top\check\Q_{-s}^2\check\z_s\check a_s\check a_k}{m}\left(\frac{F_{i_ki_k}}{\check\gamma_{i_{s,k}}}-1\right)\right]}_{W_{1012}} \\
 %%%
 &+\underbrace{\frac{m(m-1)}{m^2(1-\delta_{\zeta^{'''}})}\EE_{\zeta^{''}}\left[\frac{\check\z_k^\top\check\Q_{-sk}\check\z_s\check\z_s^\top\check\Q_{-s}^2\check\z_s\check a_s\check a_k}{m}\right]}_{W_{1013}}.
\end{align*}

 For the term $W_{1013}$, we have 
\begin{align}
W_{1013}&=\frac{m(m-1)}{m^2(1-\delta_{\zeta^{'''}})}\EE_{\zeta^{''}}\left[\frac{\check\z_k^\top\check\Q_{-sk}\check\z_s\check\z_s^\top\check\Q_{-sk}\check\Q_{-s}\check\z_s\check a_s\check a_k}{m}\right]\nonumber\\
&-\frac{m(m-1)}{m^2(1-\delta_{\zeta^{'''}})}\EE_{\zeta^{''}}\left[\frac{\check\z_k^\top\check\Q_{-sk}\check\z_s\check\z_s^\top\check\Q_{-sk}\check\z_k\check\z_k^\top\check\Q_{-sk}\check\Q_{-s}\check\z_sF_{i_ki_k}\check a_s\check a_k}{m^2\check\gamma_{i_{s,k}}}\right]\nonumber\\
&=\underbrace{-\frac{m(m-1)}{m^2(1-\delta_{\zeta^{'''}})}\EE_{\zeta^{''}}\left[\frac{\check\z_k^\top\check\Q_{-sk}\check\z_s\check\z_s^\top\check\Q^2_{-sk}\check\z_k\check\z_k^\top\check\Q_{-sk}\check\z_sF_{i_ki_k}\check a_s\check a_k}{m^2\check\gamma_{i_{s,k}}}\right]}_{W_{10131}}\nonumber\\
&+\underbrace{\frac{m(m-1)}{m^2(1-\delta_{\zeta^{'''}})}\EE_{\zeta^{''}}\left[\frac{\check\z_k^\top\check\Q_{-sk}\check\z_s\check\z_s^\top\check\Q^2_{-sk}\check\z_s\check a_s\check a_k}{m}\right]}_{W_{10132}}\nonumber\\
&+\underbrace{\frac{m(m-1)}{m^2(1-\delta_{\zeta^{'''}})}\EE_{\zeta^{''}}\left[\frac{\check\z_k^\top\check\Q_{-sk}\check\z_s\check\z_s^\top\check\Q_{-sk}\check\z_k\check\z_k^\top\check\Q^2_{-sk}\check\z_k\check\z_k^\top\check\Q_{-sk}\check\z_sF_{i_ki_k}^2\check a_s\check a_k}{m^3\check\gamma_{i_{s,k}}^2}\right]}_{W_{10133}}\nonumber\\
&\underbrace{-\frac{m(m-1)}{m^2(1-\delta_{\zeta^{'''}})}\EE_{\zeta^{''}}\left[\frac{\check\z_k^\top\check\Q_{-sk}\check\z_s\check\z_s^\top\check\Q_{-sk}\check\z_k\check\z_k^\top\check\Q^2_{-sk}\check\z_sF_{i_ki_k}\check a_s\check a_k}{m^2\check\gamma_{i_{s,k}}}\right]}_{W_{10134}}. \label{eq:w_1013}
\end{align}
Applying  \eqref{eq_W_62} yields 
\begin{align*}
    W_{10131}&= W_{10134}=-\frac{m(m-1)}{m^2(1-\delta_{\zeta^{'''}})}\EE_{\zeta^{''}}\left[\frac{\check\z_k^\top\check\Q^2_{-sk}\check\z_s\check\z_s^\top\check\Q_{-sk}\check\z_k\check\z_k^\top\check\Q_{-sk}\check\z_sF_{i_ki_k}\check a_s\check a_k}{m^2\check\gamma_{i_{s,k}}}\right]\\
    & = O\left(\frac{\theta^2_{\max}(\X)p^2}{m^2}\right)\cdot\|\rr\|^2.
\end{align*}
Taking the fact $\X_{\H}^\top\rr=0$, we further obtain
\begin{align*}
W_{10132}= \frac{m(m-1)}{m^2(1-\delta_{\zeta^{'''}})}\EE_{\zeta^{''}}\left[\rr^\top \X_{\H}\check\Q_{-sk}\X_{\H}^\top\diag\left\{\frac{\x_{\H_i}^\top\check\Q^2_{-sk}\x_{\H_i}}{m\pi_i}\right\}^n_{i=1}\rr\right]=0.
\end{align*}
Recalling  \eqref{eq:w_82}, we 
 bound the term $W_{10133}$:
\begin{align*}
    W_{10133} &=\frac{m(m-1)}{m^2(1-\delta_{\zeta^{'''}})}\EE_{\zeta^{''}}\left[\frac{\check\z_k^\top\check\Q^2_{-sk}\check\z_k\check\z_k^\top\check\Q_{-sk}\check\z_s\check\z_s^\top\check\Q_{-sk}\check\z_k\check\z_k^\top\check\Q_{-sk}\check\z_sF^2_{i_ki_k}\check a_s\check a_k}{m^3\check\gamma_{i_{s,k}}^2}\right] \\ 
    &= O\left(\frac{\theta^3_{\max}(\X)p^3}{m^3}\right)\cdot\|\rr\|^2.
\end{align*}
We thus conclude 
\begin{align*}
     W_{1013} = O\left(\frac{\theta^2_{\max}(\X)p^2}{m^2}\right)\cdot\|\rr\|^2.
\end{align*}
Subsequently, we continue to bound the term  $ W_{1011}$.
We rewrite 
\begin{align*}
    W_{1011}&=\underbrace{-\frac{m(m-1)}{m^2(1-\delta_{\zeta^{'''}})}\EE_{\zeta^{''}}\left[\frac{\check\z_k^\top\check\Q_{-sk}\check\z_s\check\z_s^\top\check\Q^2_{-sk}\check\z_k\check\z_k^\top\check\Q_{-sk}\check\z_sF^2_{i_ki_k}\check a_s\check a_k}{m^2\check\gamma^2_{i_{s,k}}}\left(\frac{F_{i_si_s}^2}{\check\gamma^2_{i_{k,s}}}-1\right)\right]}_{W_{10111}}\\
&+\underbrace{\frac{m(m-1)}{m^2(1-\delta_{\zeta^{'''}})}\EE_{\zeta^{''}}\left[\frac{\check\z_k^\top\check\Q_{-sk}\check\z_s\check\z_s^\top\check\Q^2_{-sk}\check\z_sF_{i_ki_k}\check a_s\check a_k}{m\check\gamma_{i_{s,k}}}\left(\frac{F_{i_si_s}^2}{\check\gamma^2_{i_{k,s}}}-1\right)\right]}_{W_{10112}}\\
&+\underbrace{\frac{m(m-1)}{m^2(1-\delta_{\zeta^{'''}})}\EE_{\zeta^{''}}\left[\frac{\check\z_k^\top\check\Q_{-sk}\check\z_s\check\z_s^\top\check\Q_{-sk}\check\z_k\check\z_k^\top\check\Q^2_{-sk}\check\z_k\check\z_k^\top\check\Q_{-sk}\check\z_sF^3_{i_ki_k}\check a_s\check a_k}{m^3\check\gamma_{i_{s,k}}^3}\left(\frac{F_{i_si_s}^2}{\check\gamma^2_{i_{k,s}}}-1\right)\right]}_{W_{10113}}\\
&\underbrace{-\frac{m(m-1)}{m^2(1-\delta_{\zeta^{'''}})}\EE_{\zeta^{''}}\left[\frac{\check\z_k^\top\check\Q_{-sk}\check\z_s\check\z_s^\top\check\Q_{-sk}\check\z_k\check\z_k^\top\check\Q^2_{-sk}\check\z_sF^2_{i_ki_k}\check a_s\check a_k}{m^2\check\gamma^2_{i_{s,k}}}\left(\frac{F_{i_si_s}^2}{\check\gamma^2_{i_{k,s}}}-1\right)\right]}_{W_{10114}}.
\end{align*}
Following  \eqref{eq:check_Q_quadratic_two} in \Cref{lem:niu_auxiliary} and the bound on \eqref{eq:w_82}, we explore  the terms $W_{10111}$ and $W_{10114}$:
\begin{align}\label{eq:w_10111_10114}
 &W_{10111}=    W_{10114}=-\frac{m(m-1)}{m^2(1-\delta_{\zeta^{'''}})}\EE_{\zeta^{''}}\left[\frac{\check\z_k^\top\check\Q^2_{-sk}\check\z_s\check\z_s^\top\check\Q_{-sk}\check\z_k\check\z_k^\top\check\Q_{-sk}\check\z_sF^2_{i_ki_k}\check a_s\check a_k}{m^2\check\gamma^2_{i_{s,k}}}\left(\frac{F_{i_si_s}^2}{\check\gamma^2_{i_{k,s}}}-1\right)\right]\nonumber\\
 &=- \frac{m(m-1)}{m^2(1-\delta_{\zeta^{'''}})}\sum^n_{i=1}\EE_{\zeta^{''}}\left[\frac{\check\z_k^\top\check\Q_{-sk}^2\x_{\H_i}\x_{\H_i}^\top\check\Q_{-sk}\check\z_k\check\z_k^\top\check\Q_{-sk}\x_{\H_i}F^2_{i_ki_k}r_i\check a_k}{m^2\check\gamma^2_{i_{s,k}}\pi_i}\left(\frac{F_{ii}^2}{\bar\gamma_{i_{k,s,i}}^2}-1\right)\right]
 \nonumber\\
 &\leq \sum^n_{i=1}\EE_{\zeta^{''}}\left[\frac{\x_{\H_i}\check\Q_{-sk}^2\check\z_k\check\z_k^\top\check\Q_{-sk}^2\x_{\H_i}(\x_{\H_i}^\top\check\Q_{-sk}\check\z_k\check\z_k^\top\check\Q_{-sk}\x_{\H_i})^{1/2} F^2_{i_ki_k} r^2_i}{m^2\check\gamma^2_{i_{s,k}}\pi^{3/2}_i}\left|\frac{F_{ii}^2}{\bar\gamma_{i_{k,s,i}}^2}-1\right|\right]\nonumber\\
 &+\sum^n_{i=1}\EE_{\zeta^{''}}\left[\frac{F^2_{i_ki_k}\check a^2_k(\x_{\H_i}^\top\check\Q_{-sk}\check\z_k\check\z_k^\top\check\Q_{-sk}\x_{\H_i})^{3/2}}{m^2\check\gamma^2_{i_{s,k}}\pi^{1/2}_i}\left|\frac{ F_{ii}^2}{\bar\gamma_{i_{k,s,i}}^2}-1\right|\right]\nonumber\\
 &\leq \sum^n_{i=1}\EE_{\zeta^{''}}\left[\frac{\x_{\H_i}\check\Q_{-sk}^2\check\z_k\check\z_k^\top\check\Q_{-sk}^2\x_{\H_i}(\x_{\H_i}^\top\check\Q_{-sk}\check\z_k\check\z_k^\top\check\Q_{-sk}\x_{\H_i})^{1/2}F^2_{i_ki_k} r^2_i}{m^2\check\gamma^2_{i_{s,k}}\pi^{3/2}_i}\left(\frac{F_{ii}}{\bar\gamma_{i_{k,s,i}}}\left|\frac{ F_{ii}}{\bar\gamma_{i_{k,s,i}}}-1\right|+\left|\frac{ F_{ii}}{\bar\gamma_{i_{k,s,i}}}-1\right|\right)\right]\nonumber\\
 &+\sum^n_{i=1}\EE_{\zeta^{''}}\left[\frac{F^2_{i_ki_k}\check a^2_k(\x_{\H_i}^\top\check\Q_{-sk}\check\z_k\check\z_k^\top\check\Q_{-sk}\x_{\H_i})^{3/2}}{m^2\check\gamma^2_{i_{s,k}}\pi^{1/2}_i}\left(\frac{F_{ii}}{\bar\gamma_{i_{k,s,i}}}\left|\frac{ F_{ii}}{\bar\gamma_{i_{k,s,i}}}-1\right|+\left|\frac{ F_{ii}}{\bar\gamma_{i_{k,s,i}}}-1\right|\right)\right]\nonumber\\
 &\leq 3\sum^n_{i=1}\EE_{\zeta^{''}}\left[\frac{\x_{\H_i}\check\Q_{-sk}^2\check\z_k\check\z_k^\top\check\Q_{-sk}^2\x_{\H_i}(\x_{\H_i}^\top\check\Q_{-sk}\check\z_k\check\z_k^\top\check\Q_{-sk}\x_{\H_i})^{1/2}F^2_{i_ki_k} r^2_i}{m^2\check\gamma^2_{i_{s,k}}\pi^{3/2}_i}\left|\frac{ F_{ii}}{\bar\gamma_{i_{k,s,i}}}-1\right|\right]\nonumber\\
 &+3\sum^n_{i=1}\EE_{\zeta^{''}}\left[\frac{F^2_{i_ki_k}\check a^2_k(\x_{\H_i}^\top\check\Q_{-sk}\check\z_k\check\z_k^\top\check\Q_{-sk}\x_{\H_i})^{3/2}}{m^2\check\gamma^2_{i_{s,k}}\pi^{1/2}_i}\left|\frac{F_{ii}}{\bar\gamma_{i_{k,s,i}}}-1\right|\right]\nonumber\\
 &\leq \frac{24\theta_{\max}(\X)p}{m}\sum^n_{i=1}\EE_{\zeta^{''}}\left[\frac{\x_{\H_i}\check\Q_{-sk}^2\check\z_k\check\z_k^\top\check\Q_{-sk}^2\x_{\H_i}F^2_{i_ki_k}r^2_i}{m\pi_i}\left|\frac{F_{ii}}{\bar\gamma_{i_{k,s,i}}}-1\right|\right]\nonumber\\
 &+\frac{24\theta_{\max}(\X)p}{m}\sum^n_{i=1}\EE_{\zeta^{''}}\left[\frac{F^2_{i_ki_k}\check a^2_k\check\z_k^\top\check\Q_{-sk}\x_{\H_i}\x_{\H_i}^\top\check\Q_{-sk}\check\z_k}{m}\left|\frac{F_{ii}}{\bar\gamma_{i_{k,s,i}}}-1\right|\right]\nonumber\\
 &\leq \frac{96\theta_{\max}(\X)p}{m}\sum^n_{i=1}\sum^n_{j=1}\EE_{\zeta^{''}}\left[\frac{\x_{\H_i}\check\Q_{-sk}^2\x_{\H_j}\x_{\H_j}^\top\check\Q_{-sk}^2\x_{\H_i}r^2_i}{m\pi_i}\left|\frac{ F_{ii}}{\bar\gamma_{i_{k,s,i}}}-1\right|\right]\nonumber\\
 &+\frac{24\theta_{\max}(\X)p}{m}\EE_{\zeta^{''}}\left[\frac{F^2_{i_ki_k}\check a^2_k\check\z_k^\top\check\Q_{-sk}\X_{\H}^\top\X_{\H}\check\Q_{-sk}\check\z_k}{m}\max\limits_{1\leq i\leq n}\left|\frac{ F_{ii}}{\bar\gamma_{i_{k,s,i}}}-1\right|\right]\nonumber\\
 &\leq \frac{96\theta_{\max}(\X)p}{m}\sum^n_{i=1}\EE_{\zeta^{''}}\left[\frac{\x_{\H_i}\check\Q_{-sk}^2\X_{\H}^\top\X_{\H}\check\Q_{-sk}^2\x_{\H_i}r^2_i}{m\pi_i}\left|\frac{ F_{ii}}{\bar\gamma_{i_{k,s,i}}}-1\right|\right]\nonumber\\
 &+\frac{1536\theta^2_{\max}(\X)p^2}{m^2}\EE_{\zeta^{''}}\left[F^2_{i_ki_k}\check a^2_k\max\limits_{1\leq i\leq n}\left|\frac{ F_{ii}}{\bar\gamma_{i_{k,s,i}}}-1\right|\right]\nonumber\\
 &\leq \frac{12 \times 8^5\theta^2_{\max}(\X)p^2\|\rr\|^2}{m^2}\EE_{\zeta^{''}}\left[\max\limits_{1\leq i\leq n}\left|\frac{ F_{ii}}{\bar\gamma_{i_{k,s,i}}}-1\right|\right]+\frac{6144\theta^2_{\max}(\X)p^2}{m^2}\EE_{\zeta^{''}}\left[\check a^2_k\max\limits_{1\leq i\leq n}\left|\frac{F_{ii}}{\bar\gamma_{i_{k,s,i}}}-1\right|\right]\nonumber\\
 &\leq \frac{12 \times 8^5\theta^2_{\max}(\X)p^2\|\rr\|^2}{m^2}\EE_{\zeta^{''}}\left[\max\limits_{1\leq i\leq n}\left|\frac{ F_{ii}}{\bar\gamma_{i_{k,s,i}}}-1\right|\right]+\frac{6144\theta^2_{\max}(\X)p^2\|\rr\|^2}{m^2}\EE_{\zeta^{''}}\left[\max\limits_{1\leq i\leq n}\left|\frac{F_{ii}}{\bar\gamma_{i_{k,s,i}}}-1\right|\right]\nonumber\\
 &= O\left(\frac{\log n}{\loglog n}\sqrt{\frac{\theta^7_{\max}(\X)p^7}{m^7}}\right)\cdot\|\rr\|^2.
\end{align}
Similarly, we get
\begin{align*}
    W_{10113}&=\frac{m(m-1)}{m^2(1-\delta_{\zeta^{'''}})}\EE_{\zeta^{''}}\left[\frac{\check\z_k^\top\check\Q^2_{-sk}\check\z_k\check\z_k^\top\check\Q_{-sk}\check\z_s\check\z_s^\top\check\Q_{-sk}\check\z_k\check\z_k^\top\check\Q_{-sk}\check\z_sF_{i_ki_k}^3\check a_s\check a_k}{m^3\check\gamma_{i_{s,k}}^3}\left(\frac{F_{i_si_s}^2}{\check\gamma^2_{i_{k,s}}}-1\right)\right]\\
    &= O\left(\frac{\log n}{\loglog n}\sqrt{\frac{\theta^9_{\max}(\X)p^9}{m^9}}\right)\cdot\|\rr\|^2.
\end{align*}
We then derive  
\begin{align}
&    W_{10112}=\frac{m(m-1)}{m^2(1-\delta_{\zeta^{'''}})}\EE_{\zeta^{''}}\left[\rr^\top \diag\left\{\frac{F_{jj}}{\bar\gamma_{i_{s,k,j}}}\right\}^n_{j=1}\X_{\H}\check\Q_{-sk}\X_{\H}^\top\diag\left\{\frac{\x_{\H_i}^\top\check\Q_{-sk}\x_{\H_i}}{m\pi_i}\left(\frac{F_{ii}^2}{\bar\gamma_{i_{k,s,i}}^2}-1\right)\right\}^n_{i=1}\rr\right]\nonumber\\
   &\leq \frac{256\theta_{\max}(\X)p\|\rr\|^2}{m}\EE_{\zeta^{''}}\left[\max\limits_{1\leq i\leq n}\left|\frac{F_{ii}^2}{\bar\gamma^2_{i_{k,s,i}}}-1\right|\right] \nonumber\\
   &\leq \frac{256\theta_{\max}(\X)p\|\rr\|^2}{m}\left(\EE_{\zeta^{''}}\left[\max\limits_{1\leq i\leq n}\frac{F_{ii}}{\bar\gamma_{i_{k,s,i}}}\left|\frac{F_{ii}}{\bar\gamma_{i_{k,s,i}}}-1\right|\right]+\EE_{\zeta^{''}}\left[\max\limits_{1\leq i\leq n}\left|\frac{F_{ii}}{\bar\gamma_{i_{k,s,i}}}-1\right|\right]\right)\nonumber\\
   & \leq \frac{768\theta_{\max}(\X)p\|\rr\|^2}{m}\EE_{\zeta^{''}}\left[\max\limits_{1\leq i\leq n}\left|\frac{F_{ii}}{\bar\gamma_{i_{k,s,i}}}-1\right|\right] = O\left(\frac{\log n}{\loglog n}\sqrt{\frac{\theta^5_{\max}(\X)p^5}{m^5}}\right)\cdot\|\rr\|^2.\label{eq:w_10112}
\end{align}
Thus, we get
\begin{align*}
     W_{1011}= O\left(\frac{\log n}{\loglog n}\sqrt{\frac{\theta^5_{\max}(\X)p^5}{m^5}}\right)\cdot\|\rr\|^2.
\end{align*}
Now, we turn to bound the term $W_{1012}$:
\begin{align*}
    W_{1012}&=-\underbrace{\frac{m(m-1)}{m^2(1-\delta_{\zeta^{'''}})}\EE_{\zeta^{''}}\left[\frac{\check\z_k^\top\check\Q_{-sk}\check\z_s\check\z_s^\top\check\Q^2_{-sk}\check\z_k\check\z_k^\top\check\Q_{-sk}\check\z_sF_{i_ki_k}\check a_s\check a_k}{m^2\check\gamma_{i_{s,k}} }\left(\frac{F_{i_ki_k}}{\check\gamma_{i_{s,k}}}-1\right)\right]}_{W_{10121}}\\
&+\underbrace{\frac{m(m-1)}{m^2(1-\delta_{\zeta^{'''}})}\EE_{\zeta^{''}}\left[\frac{\check\z_k^\top\check\Q_{-sk}\check\z_s\check\z_s^\top\check\Q^2_{-sk}\check\z_s\check a_s\check a_k}{m}\left(\frac{F_{i_ki_k}}{\check\gamma_{i_{s,k}}}-1\right)\right]}_{W_{10122}}\\
&+\underbrace{\frac{m(m-1)}{m^2(1-\delta_{\zeta^{'''}})}\EE_{\zeta^{''}}\left[\frac{\check\z_k^\top\check\Q_{-sk}\check\z_s\check\z_s^\top\check\Q_{-sk}\check\z_k\check\z_k^\top\check\Q^2_{-sk}\check\z_k\check\z_k^\top\check\Q_{-sk}\check\z_sF_{i_ki_k}^2\check a_s\check a_k}{m^3\check\gamma_{i_{s,k}}^2}\left(\frac{F_{i_ki_k}}{\check\gamma_{i_{s,k}}}-1\right)\right]}_{W_{10123}}\\
&-\underbrace{\frac{m(m-1)}{m^2(1-\delta_{\zeta^{'''}})}\EE_{\zeta^{''}}\left[\frac{\check\z_k^\top\check\Q_{-sk}\check\z_s\check\z_s^\top\check\Q_{-sk}\check\z_k\check\z_k^\top\check\Q^2_{-sk}\check\z_sF_{i_ki_k}\check a_s\check a_k}{m^2\check\gamma_{i_{s,k}} }\left(\frac{F_{i_ki_k}}{\check\gamma_{i_{s,k}}}-1\right)\right]}_{W_{10124}}.
\end{align*}
We begin by bounding the  term $W_{10121}$.
Using \eqref{eq:check_Q_quadratic_two} in \Cref{lem:niu_auxiliary} and \eqref{eq:w_10111_10114}, we get
\begin{align*}
 W_{10121}&=    W_{10124}=\frac{m(m-1)}{m^2(1-\delta_{\zeta^{'''}})}\EE_{\zeta^{''}}\left[\frac{\check\z_k^\top\check\Q^2_{-sk}\check\z_s\check\z_s^\top\check\Q_{-sk}\check\z_k\check\z_k^\top\check\Q_{-sk}\check\z_sF_{i_ki_k}\check a_s\check a_k}{m^2\check\gamma_{i_{s,k}}}\left(\frac{F_{i_ki_k}}{\check\gamma_{i_{s,k}}}-1\right)\right]\\
 &= \frac{m(m-1)}{m^2(1-\delta_{\zeta^{'''}})}\sum^n_{i=1}\EE_{\zeta^{''}}\left[\frac{\check\z_k^\top\check\Q_{-sk}^2\x_{\H_i}\x_{\H_i}^\top\check\Q_{-sk}\check\z_k\check\z_k^\top\check\Q_{-sk}\x_{\H_i}F_{i_ki_k}r_i\check a_k}{m^2\check\gamma_{i_{s,k}}\pi_i}\left(\frac{F_{i_ki_k}}{\check\gamma_{i_{s,k}}}-1\right)\right]
 \nonumber\\
 &\leq \sum^n_{i=1}\EE_{\zeta^{''}}\left[\frac{\x_{\H_i}\check\Q_{-sk}^2\check\z_k\check\z_k^\top\check\Q_{-sk}^2\x_{\H_i}(\x_{\H_i}^\top\check\Q_{-sk}\check\z_k\check\z_k^\top\check\Q_{-sk}\x_{\H_i})^{1/2}F_{i_ki_k}r^2_i}{m^2\check\gamma_{i_{s,k}}\pi^{3/2}_i}\left|\frac{F_{i_ki_k}}{\check\gamma_{i_{s,k}}}-1\right|\right]\nonumber\\
 &+\sum^n_{i=1}\EE_{\zeta^{''}}\left[\frac{F_{i_ki_k}\check a^2_k(\x_{\H_i}^\top\check\Q_{-sk}\check\z_k\check\z_k^\top\check\Q_{-sk}\x_{\H_i})^{3/2}}{m^2\check\gamma_{i_{s,k}}\pi^{1/2}_i}\left|\frac{F_{i_ki_k}}{\check\gamma_{i_{s,k}}}-1\right|\right]\nonumber\\
 &\leq \frac{8\theta_{\max}(\X)p}{m}\sum^n_{i=1}\EE_{\zeta^{''}}\left[\frac{\x_{\H_i}\check\Q_{-sk}^2\check\z_k\check\z_k^\top\check\Q_{-sk}^2\x_{\H_i}F_{i_ki_k}r^2_i}{m\pi_i}\left|\frac{F_{i_ki_k}}{\check\gamma_{i_{s,k}}}-1\right|\right]\\
 &+\frac{8\theta_{\max}(\X)p}{m}\sum^n_{i=1}\EE_{\zeta^{''}}\left[\frac{F_{i_ki_k}\check a^2_k\check\z_k^\top\check\Q_{-sk}\x_{\H_i}\x_{\H_i}^\top\check\Q_{-sk}\check\z_k}{m}\left|\frac{F_{i_ki_k}}{\check\gamma_{i_{s,k}}}-1\right|\right]\nonumber\\
 &\leq \frac{8\theta_{\max}(\X)p}{m}\sum^n_{i=1}\sum^n_{j=1}\EE_{\zeta^{''}}\left[\frac{\x_{\H_i}\check\Q_{-sk}^2\x_{\H_j}\x_{\H_j}^\top\check\Q_{-sk}^2\x_{\H_i}F_{jj}r^2_i}{m\pi_i}\left|\frac{F_{jj}}{\bar\gamma_{i_{s,k,j}}}-1\right|\right]\\
 &+\frac{8\theta_{\max}(\X)p}{m}\EE_{\zeta^{''}}\left[\frac{F_{i_ki_k}\check a^2_k\check\z_k^\top\check\Q_{-sk}\X_{\H}^\top\X_{\H}\check\Q_{-sk}\check\z_k}{m}\left|\frac{F_{i_ki_k}}{\check\gamma_{i_{s,k}}}-1\right|\right]\nonumber\\
 &\leq \frac{16\theta_{\max}(\X)p}{m}\sum^n_{i=1}\EE_{\zeta^{''}}\left[\frac{\x_{\H_i}\check\Q_{-sk}^2\X_{\H}^\top\X_{\H}\check\Q_{-sk}^2\x_{\H_i}r^2_i}{m\pi_i}\max\limits_{1\leq j\leq n}\left|\frac{F_{jj}}{\bar\gamma_{i_{s,k,j}}}-1\right|\right]\nonumber\\
 &+\frac{512\theta^2_{\max}(\X)p^2}{m^2}\EE_{\zeta^{''}}\left[F_{i_ki_k}\check a^2_k\left|\frac{F_{i_ki_k}}{\check\gamma_{i_{s,k}}}-1\right|\right]\nonumber\\
 &\leq \frac{2\times 8^5\theta^2_{\max}(\X)p^2\|\rr\|^2}{m^2}\EE_{\zeta^{''}}\left[\max\limits_{1\leq j\leq n}\left|\frac{F_{jj}}{\bar\gamma_{i_{s,k,j}}}-1\right|\right]+\frac{512\theta^2_{\max}(\X)p^2}{m^2}\EE_{\zeta^{''}}\left[F_{i_ki_k}\check a^2_k\left|\frac{F_{i_ki_k}}{\check\gamma_{i_{s,k}}}-1\right|\right]\nonumber\\
 &= O\left(\frac{\log n}{\loglog n}\sqrt{\frac{\theta^7_{\max}(\X)p^7}{m^7}}\right)\|\rr\|^2.
\end{align*}
Analogously,  we derive
\begin{align*}
    W_{10123}&=\frac{m(m-1)}{m^2(1-\delta_{\zeta^{'''}})}\EE_{\zeta^{''}}\left[\frac{\check\z_k^\top\check\Q^2_{-sk}\check\z_k\check\z_k^\top\check\Q_{-sk}\check\z_s\check\z_s^\top\check\Q_{-sk}\check\z_k\check\z_k^\top\check\Q_{-sk}\check\z_sF_{i_ki_k}^2\check a_s\check a_k}{m^3\check\gamma_{i_{s,k}}^2}\left(\frac{F_{i_ki_k}}{\check\gamma_{i_{s,k}}}-1\right)\right]\\
    &= O\left(\frac{\log n}{\loglog n}\sqrt{\frac{\theta^9_{\max}(\X)p^9}{m^9}}\right)\cdot\|\rr\|^2.
\end{align*}
Recalling  $\X_{\H}^\top\rr= \mathbf{0}$, together with the bound on $W_{10112}$ in \eqref{eq:w_10112}, we attain
\begin{align*}
   W_{10122}& =\frac{m(m-1)}{m^2(1-\delta_{\zeta^{'''}})}\EE_{\zeta^{''}}\left[\rr^\top \diag\left\{\frac{F_{jj}}{\bar\gamma_{i_{s,k,j}}}-1\right\}^n_{j=1}\X_{\H}\check\Q_{-sk}\X_{\H}^\top\diag\left\{\frac{\x_{\H_i}^\top\check\Q_{-sk}\x_{\H_i}}{m\pi_i}\right\}^n_{i=1}\rr\right]\\
   &\leq \frac{128\theta_{\max}(\X)p\|\rr\|^2}{m}\EE_{\zeta^{''}}\left[\max\limits_{1\leq j\leq n}\left|\frac{F_{jj}}{\bar\gamma_{i_{s,k,j}}}-1\right|\right] = O\left(\frac{\log n}{\loglog n}\sqrt{\frac{\theta^5_{\max}(\X)p^5}{m^5}}\right)\cdot\|\rr\|^2.
\end{align*}
Then, we get 
\begin{align*}
      W_{1012}= O\left(\frac{\log n}{\loglog n}\sqrt{\frac{\theta^5_{\max}(\X)p^5}{m^5}}\right)\cdot\|\rr\|^2.
\end{align*}
Consequently, it follows that 
\begin{align*}
    W_{101}= O\left(\frac{\log n}{\loglog n} \cdot \frac{\theta^2_{\max}(\X)p^2}{m^2}\right)\cdot\|\rr\|^2.
\end{align*}
Now, we proceed to  explore the term $W_{102}$:
\begin{align*}
  W_{102}&=  \frac{m(m-1)}{m^2(1-\delta_{\zeta^{'''}})}\EE_{\zeta^{''}}\left[\frac{\check\z_k^\top\check\Q_{-sk}\check\z_s\check\z_s^\top\check\Q_{-s}^2\check\z_sF^2_{i_si_s}F_{i_ki_k}\check a_s\check a_k}{m\check\gamma_{i_{s,k}}}\left(\frac{1}{\check\gamma^2_{i_{s}}}-\frac{1}{\check\gamma^2_{i_{k,s}}}\right)\right]\\
  &=\underbrace{\frac{m(m-1)}{m^2(1-\delta_{\zeta^{'''}})}\EE_{\zeta^{''}}\left[\frac{2\check\z_k^\top\check\Q_{-sk}\check\z_s\check\z_s^\top\check\Q_{-s}^2\check\z_sF^3_{i_si_s}F^2_{i_ki_k}\check a_s\check a_k\check\z_s^\top\check\Q_{-sk}\check\z_k\check\z_k^\top\check\Q_{-sk}\check\z_s}{m^3\check\gamma^2_{i_{s,k}}\check\gamma^2_{i_{s}}\check\gamma_{i_{k,s}}}\right]}_{W_{1021}}\\
  &-\underbrace{\frac{m(m-1)}{m^2(1-\delta_{\zeta^{'''}})}\EE_{\zeta^{''}}\left[\frac{\check\z_k^\top\check\Q_{-sk}\check\z_s\check\z_s^\top\check\Q_{-s}^2\check\z_sF^4_{i_si_s}F^3_{i_ki_k}\check a_s\check a_k(\check\z_s^\top\check\Q_{-sk}\check\z_k\check\z_k^\top\check\Q_{-sk}\check\z_s)^2}{m^5\check\gamma^3_{i_{s,k}}\check\gamma^2_{i_{s}}\check\gamma^2_{i_{k,s}}}\right]}_{W_{1022}}.
\end{align*}
Considering the bound on  $W_{82}$ in \eqref{eq:w_82}, we have 

\begin{align*}
  W_{1021}
 &= \frac{2m(m-1)}{m^2(1-\delta_{\zeta^{'''}})}\sum^n_{i=1}\EE_{\zeta^{''}}\left[\frac{\check\z_k^\top\check\Q_{-sk}\x_{\H_i}\x_{\H_i}^\top\check\Q^2_{-s}\x_{\H_i}\x_{\H_i}^\top\check\Q_{-sk}\check\z_k\check\z_k^\top\check\Q_{-sk}\x_{\H_i}F^3_{ii}F^2_{i_ki_k}r_i\check a_k}{m^3\bar\gamma_{i_{k,s,i}}\check\gamma^2_{i_{s,k}}\bar\gamma^2_{i}\pi_i^2}\right]
 \\
 &\leq 2\sum^n_{i=1}\EE_{\zeta^{''}}\left[\frac{\x_{\H_i}\check\Q_{-sk}\check\z_k\check\z_k^\top\check\Q_{-sk}\x_{\H_i}\x_{\H_i}^\top\check\Q^2_{-s}\x_{\H_i}(\x_{\H_i}^\top\check\Q_{-sk}\check\z_k\check\z_k^\top\check\Q_{-sk}\x_{\H_i})^{1/2}F^3_{ii}F^2_{i_ki_k}r^2_i}{m^3\bar\gamma_{i_{k,s,i}}\check\gamma^2_{i_{s,k}}\bar\gamma^2_{i}\pi^{5/2}_i}\right]\\
 &+2\sum^n_{i=1}\EE_{\zeta^{''}}\left[\frac{F^3_{ii}F^2_{i_ki_k}\check a^2_k\x_{\H_i}^\top\check\Q^2_{-s}\x_{\H_i}(\x_{\H_i}^\top\check\Q_{-sk}\check\z_k\check\z_k^\top\check\Q_{-sk}\x_{\H_i})^{3/2}}{m^3\bar\gamma_{i_{k,s,i}}\check\gamma^2_{i_{s,k}}\bar\gamma^2_{i}\pi^{3/2}_i}\right]\\
 &\leq \frac{8^5\theta^2_{\max}(\X)p^2}{m^2}\sum^n_{i=1}\EE_{\zeta^{''}}\left[\frac{\x_{\H_i}^\top\check\Q_{-sk}\check\z_k\check\z_k^\top\check\Q_{-sk}\x_{\H_i}r^2_i}{m\pi_i}\right]\nonumber\\
 &+\frac{8^5\theta^2_{\max}(\X)p^2}{m^2}\sum^n_{i=1}\EE_{\zeta^{''}}\left[\frac{\check a^2_k\check\z_k^\top\check\Q_{-sk}\x_{\H_i}\x_{\H_i}^\top\check\Q_{-sk}\check\z_k}{m}\right]= O\left(\frac{\theta^3_{\max}(\X)p^3}{m^3}\right)\cdot\|\rr\|^2.%\label{eq:w_82}
\end{align*}
Similarly, we bound the term $W_{1022}$:
\begin{align*}
    W_{1022}
 &= \frac{m(m-1)}{m^2(1-\delta_{\zeta^{'''}})}\sum^n_{i=1}\EE_{\zeta^{''}}\left[\frac{\check\z_k^\top\check\Q_{-sk}\x_{\H_i}\x_{\H_i}^\top\check\Q^2_{-s}\x_{\H_i}(\x_{\H_i}^\top\check\Q_{-sk}\check\z_k\check\z_k^\top\check\Q_{-sk}\x_{\H_i})^2F^4_{ii}F^3_{i_ki_k}r_i\check a_k}{m^5\bar\gamma_{i_{k,s,i}}^2\check\gamma^3_{i_{s,k}}\bar\gamma^2_{i}\pi_i^3}\right]
 \\
 &\leq \sum^n_{i=1}\EE_{\zeta^{''}}\left[\frac{\x^\top_{\H_i}\check\Q_{-sk}\check\z_k\check\z_k^\top\check\Q_{-sk}\x_{\H_i}\x_{\H_i}^\top\check\Q^2_{-s}\x_{\H_i}(\x_{\H_i}^\top\check\Q_{-sk}\check\z_k\check\z_k^\top\check\Q_{-sk}\x_{\H_i})^{3/2}F^4_{ii}F^3_{i_ki_k}r^2_i}{m^5\bar\gamma_{i_{k,s,i}}^2\check\gamma^3_{i_{s,k}}\bar\gamma^2_{i}\pi^{7/2}_i}\right]\\
 &+\sum^n_{i=1}\EE_{\zeta^{''}}\left[\frac{F^4_{ii}F^3_{i_ki_k}\check a^2_k\x_{\H_i}^\top\check\Q^2_{-s}\x_{\H_i}(\x_{\H_i}^\top\check\Q_{-sk}\check\z_k\check\z_k^\top\check\Q_{-sk}\x_{\H_i})^{5/2}}{m^5\bar\gamma_{i_{k,s,i}}^2\check\gamma^3_{i_{s,k}}\bar\gamma^2_{i}\pi^{5/2}_i}\right]\\
 &\leq \frac{2\times 8^7\theta^4_{\max}(\X)p^4}{m^4}\sum^n_{i=1}\EE_{\zeta^{''}}\left[\frac{\x_{\H_i}^\top\check\Q_{-sk}\check\z_k\check\z_k^\top\check\Q_{-sk}\x_{\H_i}r^2_i}{m\pi_i}\right]\\
 &+\frac{2\times 8^7\theta^4_{\max}(\X)p^4}{m^4}\sum^n_{i=1}\EE_{\zeta^{''}}\left[\frac{\check a^2_k\check\z_k^\top\check\Q_{-sk}\x_{\H_i}\x_{\H_i}^\top\check\Q_{-sk}\check\z_k}{m}\right]\nonumber\\
 &= O\left(\frac{\theta^5_{\max}(\X)p^5}{m^5}\right)\cdot\|\rr\|^2. 
\end{align*}
This yields 
\begin{align*}
    W_{102}= O\left(\frac{\theta^3_{\max}(\X)p^3}{m^3}\right)\cdot\|\rr\|^2. 
\end{align*}
We further conclude 
\begin{align*}
    W_{10}= O\left(\frac{\log n}{\loglog n} \cdot \frac{\theta^2_{\max}(\X)p^2}{m^2} \right) \cdot\|\rr\|^2. 
\end{align*}
Now, it remains to derive the bound of the term   $W_{11}$.  
Recalling \eqref{eq:exp_zq2zaa} and \eqref{eq:var_zq2zaa} again, 
we get 
\begin{align*}
&\EE_{\zeta^{''}}\left[\left|\frac{\check\z_k^\top\check\Q_{-sk}\check\z_s\check\z_s^\top\check\Q_{-s}^2\check\z_sF^2_{i_si_s}F_{i_ki_k}\check a_s\check a_k}{m\check\gamma_{i_s}^2\check\gamma_{i_{s,k}}}\right|\right]\leq   \frac{8^3\theta_{\max}(\X)p}{m}\EE_{\zeta^{''}}\left[|\check\z_k^\top\check\Q_{-sk}\check\z_s\check a_s\check a_k|\right] \\
%%%
&\leq\frac{8^3\theta_{\max}(\X)p}{m} \sqrt{\EE_{\zeta^{''}}\left[\check\z_k^\top\check\Q_{-sk}\check\z_s\check\z_s\check\Q_{-sk}\check\z_k\right]\EE_{\zeta^{''}}\left[\check  a^2_s\check a^2_k\right]}\\ 
&=\frac{8^3\theta_{\max}(\X)p\|\rr\|^2}{m} \sqrt{\EE_{\zeta^{''}}\left[\check\z_k^\top\check\Q_{-sk}\X_{\H}^\top\X_{\H}\check\Q_{-sk}\check\z_k\right]}\\
&= \frac{8^3\theta_{\max}(\X)p\|\rr\|^2}{m} \sqrt{\sum^n_{j=1}\EE_{\zeta^{''}}\left[\x_{\H_j}^\top\check\Q_{-sk}\X_{\H}^\top\X_{\H}\check\Q_{-sk}\x_{\H_j}\right]}\leq \frac{8^4\theta_{\max}(\X)p^{3/2}\|\rr\|^2}{m},
\end{align*}
and 
\begin{align*}
  &  \Var_{\zeta^{''}}\left[\frac{\check\z_k^\top\check\Q_{-sk}\check\z_s\check\z_s^\top\check\Q_{-s}^2\check\z_sF^2_{i_si_s}F_{i_ki_k}\check a_s\check a_k}{m\check\gamma_{i_s}^2\check\gamma_{i_{s,k}}}\right]\leq  \EE_{\zeta^{''}}\left[\frac{(\check\z_k^\top\check\Q_{-sk}\check\z_s\check\z_s^\top\check\Q_{-s}^2\check\z_sF^2_{i_si_s}F_{i_ki_k}\check a_s\check a_k)^2}{m^2\check\gamma_{i_s}^4\check\gamma_{i_{s,k}}^2}\right]\\
  &\leq \frac{8^6\theta^2_{\max}(\X)p^2}{m^2}\EE_{\zeta^{''}}\left[\check\z_k^\top\check\Q_{-sk}\check\z_s\check\z_s^\top\check\Q_{-s}\check\z_k\check  a^2_s\check a^2_k\right]\leq \frac{8^8\theta^4_{\max}(\X)p^4}{m^2}\|\rr\|^4.
\end{align*}
Following  the bound on the term  $W_{7}$,
we similarly obtain, for $\delta_{\zeta^{'''}}\leq m^{-3}$,
\begin{align*}
 W_{11}=O\left(\frac{1}{m}+\frac{\theta^4_{\max}(\X)p^4}{m^4}\right)\cdot \|\rr\|^2.
\end{align*}
We thus complete the proof.

\subsection{Proof  of Corollary~\ref{coro:sub_ols_lev}}\label{subsec:auxillary_resul_theodesamols_lev}

% \begin{proof}
    The proof of \Cref{coro:sub_ols_lev} is based on the same proof strategy developed for \Cref{theo:de_sam_ols}. For completeness, we briefly outline the necessary modifications.
The only substantive departure from the proof of \Cref{theo:de_sam_ols} lies in the auxiliary estimates used to control several key terms. In the leverage-score sampling setting,  we rely on the refined bounds  $\left\| \EE_\zeta[\check\Q^2]-\I_p-\X_{\H}^\top\bar \F\X_{\H}\right\| $,  $\| \EE_\zeta[\check\Q^2]- \EE_\zeta[\check\Q^2_{-s}]\|$, $\EE_{\zeta^{'}}\left[\max\limits_{1\leq i\leq n}\left|\frac{F_{ii}}{\bar \gamma_{i}}-1\right|\right]$  and $\EE_{\zeta^{''}}\left[\max\limits_{1\leq i\leq n}\left|\frac{F_{ii}}{\bar \gamma_{i_{k,s,i}}}-1\right|\right]$, which are established in \Cref{lem:niu_auxiliary_lev}.
Once these bounds are invoked, the remainder of the argument proceeds identically to that of \Cref{theo:de_sam_ols}.
% \end{proof}

\subsection{Proof  of Corollary~\ref{coro:debiasing_srht}}\label{subsec:auxillary_resul_theodesamols_srht}

  The proof of \Cref{coro:debiasing_srht} follows the same overall strategy as that of \Cref{theo:de_sam_ols}. While the general structure of the argument remains unchanged, the technical distinctions arise from the auxiliary bounds employed. In the SRHT setting, the required estimates for $\left\| \EE_\zeta[\check\Q^2]-\I_p-\X_{\H}^\top\bar \F\X_{\H}\right\| $, $\| \EE_\zeta[\check\Q^2]- \EE_\zeta[\check\Q^2_{-s}]\|$, $\EE_{\zeta^{'}}\left[\max\limits_{1\leq i\leq n}\left|\frac{F_{ii}}{\bar \gamma_{i}}-1\right|\right]$ and $\EE_{\zeta^{''}}\left[\max\limits_{1\leq i\leq n}\left|\frac{F_{ii}}{\bar \gamma_{i_{k,s,i}}}-1\right|\right]$, are taken from \Cref{lem:niu_auxiliary_srht}.
Aside from substituting these bounds, the remainder of the proof follows the same sequence of arguments as in \Cref{theo:de_sam_ols}.

\subsection{Auxiliary results from \cite{niu2025fundamental}}\label{subsec:auxil_results_proof_de_sam_pls}

Here, we recall several auxiliary lemmas from \cite{niu2025fundamental}, which serve as key technical ingredients in the proof of \Cref{theo:de_sam_ols}.

We first present the following two lemmas, which are derived from Proposition~3.2 and Proposition~F.1, along with their corresponding proofs given in Sections~E.2 and~F.1 of \cite{niu2025fundamental}.

\begin{lemma}[{\cite{niu2025fundamental}}]\label{lem:niu_auxiliary}
Under the settings and notations of \Cref{theo:de_sam_ols}, let $\X_{\H}=\X\H^{-1/2}$. Then,   
there exists $C > 0$ independent of $n, p$ so that for $m \geq C\theta_{\max}(\X) p\log (p/\delta)$,  $ \delta\leq m^{-3}$, $\theta_{\max}(\X)$ in \Cref{def:approx_factor}, when  conditioned on an event $\zeta$ that holds with probability at least $1-\delta$, we have
\begin{align}
   & \left\| \EE_\zeta[\check\Q^2]-\I_p-\X_{\H}^\top\bar \F\X_{\H}\right\| =O\left(\frac{\log n}{\loglog n} \sqrt{\frac{\theta^3_{\max}(\X)p^3}{m^3}}\right),\label{eq:check_Q_second_moment}\\
   &\| \EE_\zeta[\check\Q^2]- \EE_\zeta[\check\Q^2_{-s}]\|=O\left(\frac{1}{m}\right),\label{eq:check_Q_check_Q_s}\\
   &\EE_{\zeta^{'}}\left[\max\limits_{1\leq i\leq n}\left|\frac{F_{ii}}{\bar \gamma_{i}}-1\right|\right]=O\left(\frac{\log n}{\loglog n} \sqrt{\frac{\theta^3_{\max}(\X)p^3}{m^3}}\right), \label{eq:check_Q_quadratic}\\
    &\EE_{\zeta^{''}}\left[\max\limits_{1\leq i\leq n}\left|\frac{F_{ii}}{\bar\gamma_{i_{k,s,i}}}-1\right|\right]=O\left(\frac{\log n}{\loglog n} \sqrt{\frac{\theta^3_{\max}(\X)p^3}{m^3}}\right). \label{eq:check_Q_quadratic_two}
\end{align}    
\end{lemma}

\begin{lemma}[{\cite{niu2025fundamental}}]\label{lem:niu_proba_space}
Under the settings and notations of \Cref{theo:de_sam_ols}, suppose that the events $\zeta_1$, $\zeta_2$, $\zeta_3$ are independent of $\z_{s}$. Define $\zeta^{'} =\bigcap^3_{j=1} \zeta_j$ and  $\delta_4=\Pr(\neg\zeta_4)$. 
For a p.s.d.\@ random matrix $\K$ (or any non-negative random variable) depending only on  random sampling $\S$,  we obtain, 
\begin{align}
    \EE_{\zeta}[\K]=\frac{\EE[(\prod^4_{j=1}\mathbf{1}_{\zeta_j})\K]}{\Pr(\zeta)} \preceq \frac{1}{1-\delta}  \EE[\mathbf{1}_{\zeta^{'}}\K] \preceq 2 \EE_{\zeta^{'}}[\K], \label{eq:condition_expec_three}
\end{align}
where $\mathbf{1}_{\zeta_j}$ is the indicator of the event $\zeta_j$.
\end{lemma}

We next state two lemmas specializing \Cref{lem:niu_auxiliary} to the settings of leverage-score sampling and the subsampled randomized Hadamard transform (SRHT). These results follow from Corollary~3.4, Corollary~3.7, Corollary~F.2, and   Corollary~F.3, and  their proofs appear in Section~E.3.2 and Section~E.3.4 of \cite{niu2025fundamental}.

\begin{lemma}[{\cite{niu2025fundamental}}] \label{lem:niu_auxiliary_lev}
 Under the settings and notations of \Cref{coro:sub_ols_lev},  then, 
there exists $C > 0$ independent of $n, p$ so that for $m \geq C\theta_{\max}(\X) p\log (p/\delta)$,  $ \delta\leq m^{-3}$, $\theta_{\max}(\X)$ in \Cref{def:approx_factor}, when  conditioned on an event $\zeta$ that holds with probability at least $1-\delta$, we have
\begin{align*}
   % &\| \EE_\zeta[\check\Q]-\I_p\|=O\left(\sqrt{\frac{\theta^3_{\max}(\X)p^3}{m^3}}+\frac{\epsilon_{\theta}\theta_{\max}(\X) p}{m} \right)\\
   & \left\| \EE_\zeta[\check\Q^2]-\I_p-\X_{\H}^\top\bar \F\X_{\H}\right\| =O\left(\frac{\log n}{\loglog n}\sqrt{\frac{\theta^3_{\max}(\X)p^3}{m^3}}+\frac{\epsilon_{\theta}\theta_{\max}(\X) p}{m} \right),\\
   &\| \EE_\zeta[\check\Q^2]- \EE_\zeta[\check\Q^2_{-s}]\|=O\left(\frac{1}{m}\right),\\
   &\EE_{\zeta^{'}}\left[\max\limits_{1\leq i\leq n}\left|\frac{F_{ii}}{\bar \gamma_{i}}-1\right|\right]=O\left(\frac{\log n}{\loglog n}\sqrt{\frac{\theta^3_{\max}(\X)p^3}{m^3}}+\frac{\epsilon_{\theta}\theta_{\max}(\X) p}{m}  \right), \\
    &\EE_{\zeta^{''}}\left[\max\limits_{1\leq i\leq n}\left|\frac{F_{ii}}{\bar\gamma_{i_{k,s,i}}}-1\right|\right]=O\left(\frac{\log n}{\loglog n}\sqrt{\frac{\theta^3_{\max}(\X)p^3}{m^3}}+\frac{\epsilon_{\theta}\theta_{\max}(\X) p}{m}  \right), 
\end{align*} 
where $F_{ii}=m/(m-p)$.
\end{lemma}

\begin{lemma}[{\cite{niu2025fundamental}}]\label{lem:niu_auxiliary_srht}
 Under the settings and notations of \Cref{coro:debiasing_srht},  then, 
 then there exists $C> 0$,  $n \exp(-p) < \delta \leq  m^{-3}$ such that
for   $m\geq C p \log (p/\delta)$, when   conditioned on the event $\zeta$ that holds with the probability at least $1-\delta$,  we have
\begin{align*}
   % &\| \EE_\zeta[\check\Q]-\I_p\|=O\left(\sqrt{\frac{p^3}{m^3}}+\frac{\sqrt{p\log(n/\delta)}}{m} \right)\\
   & \left\| \EE_\zeta[\check\Q^2]-\I_p-\X_{\H}^\top\bar \F\X_{\H}\right\| =O\left(\frac{\log n}{\loglog n}\sqrt{\frac{p^3}{m^3}}+\frac{\sqrt{p\log(n/\delta)}}{m} \right),\\
   &\| \EE_\zeta[\check\Q^2]- \EE_\zeta[\check\Q^2_{-s}]\|=O\left(\frac{1}{m}\right),\\
   &\EE_{\zeta^{'}}\left[\max\limits_{1\leq i\leq n}\left|\frac{F_{ii}}{\bar \gamma_{i}}-1\right|\right]=O\left(\frac{\log n}{\loglog n}\sqrt{\frac{p^3}{m^3}}+\frac{\sqrt{p\log(n/\delta)}}{m}\right), \\
    &\EE_{\zeta^{''}}\left[\max\limits_{1\leq i\leq n}\left|\frac{F_{ii}}{\bar\gamma_{i_{k,s,i}}}-1\right|\right]=O\left(\frac{\log n}{\loglog n}\sqrt{\frac{p^3}{m^3}}+\frac{\sqrt{p\log(n/\delta)}}{m}\right), 
\end{align*} 
where $F_{ii}=m/(m-p)$.
\end{lemma}

\section{Proof of Theorem~\ref{theo:de_sam_obliqueprojec}}\label{sec:proof_theo_de_sam_obliqueprojec}

We begin by recalling some notations used in  \Cref{sec:proof_theo_sam_ols}.
Let  $\check\z^\top_s=\ee^\top_{i_s}/\sqrt{\pi_{i_s}}\X_{\H} \in \RR^p$ with $\X_{\H}=\X\H^{-1/2}$ and $\H=\X^\top\X$, and let 
\begin{equation*}
\check\Q=(\X_{\H}^\top\check \S^\top\check\S \X_{\H} )^{-1}= \left(\sum^{m}_{s=1}\frac{1}{m}F_{i_s,i_s}\check\z_s\check{\z}^{\top}_s \right)^{-1},~~  \check\Q_{-s}=\left(\sum_{k\neq s}\frac{1}{m}F_{i_k,i_k}\check\z_k\check{\z}^{\top}_k  \right)^{-1}.
\end{equation*}
Note that for $\mathrm{rank}(\X)=p$,
\begin{align*}
\X( \check\S \X)^{\dagger} \check\S =\X(\X^\top\check\S^\top \check\S \X)^{-1}\X^\top\check\S^\top  \check\S,~~\X\X^{\dagger}=\X(\X^\top \X)^{-1}\X^\top,
\end{align*}
which implies 
\begin{align*}
&\|\E[\X( \check\S \X)^{\dagger} \check\S ] -\X\X^{\dagger}\|^2_F=\|\H^{\frac{1}{2}}\EE_{\zeta}[(\X^\top\check\S^\top\check\S\X)^{-1}\X^\top\check\S^\top\check\S\P_{\perp}]\|^2_F\\
&=\|\EE_{\zeta}[\check\Q\X_\H^\top\check\S^\top\check\S\P_{\perp}]\|^2_F=\sum^n_{l=1}\|\EE_{\zeta}[\check\Q\X_\H^\top\check\S^\top\check\S\rr_l]\|^2
\end{align*}
and 
\begin{align*}
   &  \E_{\zeta}[\|\X( \check\S \X)^{\dagger} \check\S  -\X\X^{\dagger}\|^{2}_{F}]= \E_{\zeta}[\|\H^{\frac{1}{2}}(\X^\top\check\S^\top\check\S\X)^{-1}\X^\top\check\S^\top\check\S\P_{\perp}\|^{2}_{F}]\\
     &=\E_{\zeta}[\|\check\Q\X_\H^\top\check\S^\top\check\S\P_{\perp}\|^{2}_{F}]=\sum^n_{l=1}\EE_{\zeta}[\|\check\Q\X_\H^\top\check\S^\top\check\S\rr_l\|^2]
\end{align*}
with $\rr_l$ the $l^{\text{th}}$ column of $\P_{\perp}$.
Following the same line of arguments as in the proof of \Cref{theo:de_sam_ols}, we obtain that, for each column $\rr_l$, 
 \begin{align*}
  \| \EE_{\zeta}[\check\Q\X_\H^\top\check\S^\top\check\S\rr_l]\|^2=\epsilon^2\|\rr_l\|^2,
 \end{align*}
 and 
 \begin{align*}
     \EE_{\zeta}[\|\check\Q\X_\H^\top\check\S^\top\check\S\rr_l\|^2]=\rr_l^\top\diag\Big\{\frac{\ell_i(\X)}{m\pi_{i}}\Big\}^{n}_{i=1}\rr_l+\epsilon\|\rr_l\|^2_2.
 \end{align*}
 This implies 
 \begin{align*}
     &\|\E[\X( \check\S \X)^{\dagger} \check\S ] -\X\X^{\dagger}\|^2_F=\epsilon^2\sum^n_{l=1}\|\rr_l\|^2=\epsilon^2\|\P_{\perp}\|^2_F,
 \end{align*}
 and 
 \begin{align*}
      \E_{\zeta}[\|\X( \check\S \X)^{\dagger} \check\S -\X\X^{\dagger}\|^{2}_{F}]&=\sum^n_{l=1}\rr_l^\top\diag\Big\{\frac{\ell_i(\X)}{m\pi_{i}}\Big\}^{n}_{i=1}\rr_l + \epsilon \sum^n_{l=1} \|\rr_l\|^2_2\\
     &=\tr\left(\P_{\perp}\diag\Big\{\frac{\ell_i(\X)}{m\pi_{i}}\Big\}^{n}_{i=1}\right)+\epsilon \|\P_{\perp}\|^2_F.
 \end{align*}
 This completes the proof.

\section{The full bias-variance characterizations  on fast CUR} \label{sec:proof_theo_cur}

In this section, we provide the full bias--variance characterizations in \Cref{theo:cur_bias_variance}, supplementing \Cref{theo:cur} with the corresponding variance result. We then prove \Cref{theo:cur_bias_variance} in \Cref{subsec:detail_proof_theo:cur}, which also includes the proof of \Cref{theo:cur}. Additional results related to \Cref{theo:cur_bias_variance} are provided  in \Cref{subsec:cur_lev_srht}.

Before showing  \Cref{theo:cur_bias_variance}, we first present \Cref{alg:fast_cur}, which summarizes the fast CUR decomposition procedure and the computation of the standard sampled CUR estimator $\tilde{\U}$.

\begin{algorithm}[htb]
  \caption{Fast CUR Decomposition~\cite{wang2016towards}}
  \label{alg:fast_cur}
  \begin{algorithmic}
    \STATE {\bfseries Input:} $\X\in\RR^{n\times p}$,    $\C \in \RR^{n \times c}$ consisting of $c$ columns of $\X$,   $\R \in \RR^{r \times p}$ consisting of $r$ rows of $\X$,   sample sizes $m_c$ and $m_r$, sampling probabilities  $\{\pi_{\C,i}\}^{n}_{i=1}$ and $\{\pi_{\R,j}\}^{p}_{j=1}$ with   $\pi_{\C,i} \geq 0$, $\sum^n_{i=1}\pi_{\C,i}=1$ and $\pi_{\R,j} \geq 0$, $\sum^p_{j=1}\pi_{\R,j}=1$.
    \STATE{\bfseries Output:} $\tilde{\U}\in\RR^{c\times r}$.
     % \vspace{-0.7em}
     \begin{enumerate}
     \setlength{\itemsep}{2pt}
  \setlength{\parskip}{0pt}
  \setlength{\topsep}{2pt}
  \setlength{\partopsep}{0pt}
         \item Draw sampling matrices $\S_{\C}\in \RR^{m_c\times n}$ with probabilities $\{\pi_{\C,i}\}^{n}_{i=1}$  and $\S_{\R} \in\RR^{m_r\times p}$ with probabilities $\{\pi_{\R,j}\}^{n}_{j=1}$ as in \Cref{def:RS};
         \item  Compute the sampled matrices  $\S_{\C}\C $, $\S_{\R}\R^\top$,  $ \S_{\C} \X \S^\top_{\R}$; 
         \item  Compute the solution of the sketched problem $\min\limits_{\U}  \| \S_{\C} ( \C \U \R- \X) \S^\top_{\R} \|_F^2 $ as 
         % \vspace{-1.3em}
     \begin{equation}
    \tilde{\U}=(\S_{\C}\C)^\dagger 
   \S_{\C} \X \S^\top_{\R} 
   (\R\S^\top_{\R})^\dagger.\label{eq:faster_U}  
\end{equation}
     \end{enumerate}
 % \vspace{-1.1em}
     % \STATE{\bfseries Output:} $\tilde{\U}=$.
  \end{algorithmic}
\end{algorithm}

In the following, we   provide precise statistical characterizations of the proposed debiased solution  $\check{\U}$.

\begin{theorem}[\textbf{Precise characterizations for debiased fast CUR decomposition}]\label{theo:cur_bias_variance}
For a data matrix $\X \in \RR^{n \times p}$, let $\C\in  \RR^{n \times c}$ and $\R\in \RR^{r \times p}$ be the pre-selected columns and rows of $\X$, and let $\U_{\CUR}$ be defined as in \eqref{eq:def_exact_U_CUR}.  Assume that $ \min_{i,\|\cc_i\|>0}\|\cc_i\|^2 /\|\C\|_F^2 \geq n^{-{\alpha_c}}$  and $ \min_{j,\|\rr_j\|>0}\|\rr_j\|^2 /\|\R\|_F^2 \geq p^{-{\alpha_r}}$ for some constants $\alpha_c>0$ and $\alpha_r>0$.
For standard random sampling matrices $\S_{\C}$ and $\S_{\R} $ used in \Cref{alg:fast_cur}, define their corresponding \emph{debiased} sampling matrices $ \check\S_{\C} \in \RR^{ m_{c}  \times n}$ and $ \check\S_{\R} \in \RR^{ m_{r}\times p}$ as in \eqref{eq:debias_check_S}. 
Then, there exists $C > 0$ independent of $n, p, c, r$ so that for $m_{c} \geq C\theta_{\max}(\C)c\log (c/\delta)$ and $m_{r} \geq C \theta_{\max}(\R^\top)r\log (r/\delta)$, $\delta\leq \min\{n^{-(3+\alpha_c)}, p^{-(3+\alpha_r)}\} $, when conditioned on an event $\zeta$ that holds with probability at least $1-\delta$, the debiased matrix $\check{\U}$ defined in \eqref{eq:def_check_U_fast} satisfies
\begin{equation*}
  \Bias_{\zeta}(\check{\U}) = \epsilon^2 \cdot L(\U_{\CUR}), \quad \Var_{\zeta}(\check{\U}) = \Delta_{\CUR} + \epsilon \cdot L(\U_{\CUR}),
\end{equation*}
with $\epsilon= O(\sqrt{ (\log n/\loglog n)^2\cdot\theta_{\max}^3(\C)c^3/ m_c^3}+\sqrt{ (\log p/\loglog p)^2\cdot \theta^3_{\max}(\R^\top)r^3/m_r^3})$,  for $L(\U) \equiv \|\C \U \R-\X\|^2_F$, the associated bias $\Bias_{\zeta}(\check{\U}) \equiv L(\EE_{\zeta}[\check{\U}]) - L(\U_{\CUR})$ and variance $\Var_{\zeta}(\check{\U}) \equiv \EE_{\zeta}[ L(\check{\U}) ] - L(\U_{\CUR})$, in a spirit similar to those of subsampled OLS in \Cref{def:bias_and_variance},  where we define $\Delta_{\CUR} = (\Delta_{\CUR, 1} + \Delta_{\CUR, 2})^2 + \sqrt{\epsilon \cdot L(\U_{\CUR})} (\Delta_{\CUR, 1} + \Delta_{\CUR, 2})$, with $\Delta_{\CUR, 1}= \|\diag\{\sqrt{ \ell_{i}(\C)/m_c\pi_{\C, i}} \}^{n}_{i=1}\B_1 \|_F$ and $\Delta_{\CUR, 2}= \|\diag \{\sqrt{\ell_{j}(\R^\top)/m_r\pi_{\R, j}} \}^p_{j=1}\B_2\|_F$ for   $\B_1=(\I_n-\W_{\C}\W_{\C}^\top)\X \in \RR^{n \times p}$ and $\B_2=(\I_p-\V_{\R}\V_{\R}^\top)\X^\top\in \RR^{p \times n}$, $\W_{\C} \in \RR^{n\times k_c}$ left singular vector of $\C$, and $\V_{\R}  \in \RR^{p\times k_r}$ right singular vector of $\R$.
Here, $\cc_i^\top$ is the $i^{\text{th}}$ row of $\C$, $\rr_j^\top$ is the  $j^{\text{th}}$  row of $\R^\top$,
$\theta_{\max}(\C)$ and $\theta_{\max}(\R^\top)$ as in \Cref{def:approx_factor} denote the maximum importance sampling approximation
factors for $\C$ and $\R^\top$, respectively.
\end{theorem}

\subsection{Detailed proof of Theorem~\ref{theo:cur_bias_variance}}\label{subsec:detail_proof_theo:cur}
We first recall some notations.
Let $k_c= \rm{rank}(\C)\leq c$ and $k_r= \rm{rank}(\R)\leq r$. 

For $\C$, consider its  singular value decomposition (SVD) 
$\C=\W_{\C}\Sigma_{\C}\V_{\C}^\top$,  where $\W_{\C} \in \RR^{n\times k_c}$ and  $\V_{\C} \in \RR^{c\times k_c}$  contain the left and right singular vectors of 
$\C$, respectively, and $\Sigma_{\C}=\diag(\sigma_{\C,1},\ldots,\sigma_{\C,k_c} )$, where $\sigma_{\C,1}\geq \sigma_{\C,2}\geq \ldots \geq \sigma_{\C,k_c}$ are its  nonzero singular values. Similarly, the SVD of $\R$ is $\R=\W_{\R}\Sigma_{\R}\V_{\R}^\top$, where $\W_{\R} \in \RR^{r\times k_r}$ and  $\V_{\R} \in \RR^{p\times k_r}$  contain the left and right singular vectors of  $\R$, and $\Sigma_{\R}=\diag(\sigma_{\R,1},\ldots,\sigma_{\R,k_r} )$, where $\sigma_{\R,1}\geq \sigma_{\R,2}\geq \ldots \geq \sigma_{\R,k_r}$ are its  nonzero singular values.
We  rewrite 
\begin{align*}
    \C\U_{\CUR}  \R
= \C \C^\dagger \X \R^\dagger\R=\W_{\C}\W_{\C}^\top\X\V_{\R}\V_{\R}^\top%\equiv \W_{\C}\M^\star\V_{\R}^\top.
\end{align*} 
and 
\begin{align*}
   \C\check{\U}\R=\W_{\C}\Sigma_{\C}\V_{\C}^\top(\check\S_{\C} \W_{\C} \Sigma_{\C}\V_{\C}^\top)^{\dagger}\check\S_{\C}\X\check\S_{\R}^\top  (\W_{\R}\Sigma_{\R}\V_{\R}^\top \check\S_{\R}^\top  )^{\dagger} \W_{\R}\Sigma_{\R}\V_{\R}^\top.
\end{align*}

Following the proof strategies developed in \Cref{theo:de_sam_ols} and \Cref{theo:de_sam_obliqueprojec}, the proof of \Cref{theo:cur_bias_variance} proceeds in two main steps:
\begin{enumerate}
  \item 
  construct a high probability event $\zeta$, based on subspace-embedding-type results from \Cref{lem:sub_embed}, under which the matrices $(\W_{\C} ^\top \check\S_{\C}^\top \check\S_{\C} \W_{\C} )^{-1}$ and $(\V_{\R}^\top \check\S_{\R}^\top \check\S_{\R} \V_{\R} )^{-1}$  are well conditioned; and
  \item conditioning  on that event $\zeta$, bound the quantities  $  \Bias_{\zeta}(\check{\U})$ and   $\Var_{\zeta}(\check{\U})-\Delta_{\mathrm{CUR}}$ using a ``leave-one-out'' type analysis.
\end{enumerate}

For random sampling matrix $\S_{\C}$ in \Cref{alg:fast_cur}, denote $\z^\top_{\C,s_c}=\ee^\top_{i_{s_c}}/\sqrt{\pi_{\C,i_{s_c}}}\W_{\C} \in \RR^{k_c}$ the $i_{s_c}^{th}$ row of the sketch $\S_{\C}\W_{\C} $, where $\{\pi_{\C,i}\}^n_{i=1}$ are the sampling probabilities of $\S_\C$, so that $\EE[\z_{\C,s_c}{\z}^{\top}_{\C,s_c}]=\W^\top_\C\W_{\C}$. Similarly, for random sampling matrix $\S_{\R}$ in \Cref{alg:fast_cur}, denote $\z^\top_{\R,s_r}=\ee^\top_{i_{s_r}}/\sqrt{\pi_{\R,i_{s_r}}}\V_{\R} \in \RR^{k_r}$ the $i_{s_r}^{th}$ row of the sketch $\S_{\R}\V_{\R}$, where $\{\pi_{\R,i}\}^p_{i=1}$ are the sampling probabilities of $\S_\R$, so that $\EE[\z_{\R,s_r}{\z}^{\top}_{\R,s_r}]=\V^\top_\R\V_{\R}$.
Without loss of generality, assume that $t_c=m_c/4$ and $t_r=m_r/4$ are integers, and define the event $\zeta=\zeta_{\C}\bigcap\zeta_{\R}$, where $\zeta_{\C}$ depends only on $\S_{\C}$ and $\zeta_{\R}$ depends only on $\S_{\R}$, constructed as in \eqref{eq:events}.  Specifically, for $\zeta_{\C}$, we set
\begin{align}\label{eq:events_c}
    \zeta_{\C,j}:\sum^{t_cj}_{s_c=t_c(j-1)+1}\frac{1}{t_c}\z_{\C,s_c}{\z}^{\top}_{\C, s_c}\succeq \frac{1}{2}   \W_{\C}^\top\W_{\C},~~~j=1,2,3,4,~~~ \zeta_{\C}=\bigcap^4_{j=1}\zeta_{\C,j},
\end{align}
and analogously, for $\zeta_{\R}$,
\begin{align}\label{eq:events_r}
    \zeta_{\R,j}:\sum^{t_rj}_{s_r=t_r(j-1)+1}\frac{1}{t_r}\z_{\R,s_r}{\z}^{\top}_{\R, s_r}\succeq \frac{1}{2}   \V_{\R}^\top\V_{\R},~~~j=1,2,3,4,~~~ \zeta_{\R}=\bigcap^4_{j=1}\zeta_{\R,j}.
\end{align}
Thus, we have, on  $\zeta_{\C,j}$ 
\begin{align}
 \sum^{t_cj}_{s_c=t_c(j-1)+1}\frac{F_{i_{ s_c}i_{ s_c}}}{t_c}\z_{\C, s_c}{\z}^{\top}_{\C, s_c} \succeq  \sum^{t_cj}_{s_c=t_c(j-1)+1}\frac{1}{t_c}\z_{\C,s_c}{\z}^{\top}_{\C, s_c} \succeq \frac{1}{2}   \W_{\C}^\top\W_{\C}, ~~~j=1,2,3,4,\nonumber
\end{align}
and 
  similarly, on $\zeta_{\R,j}$,
\begin{align}
 \sum^{t_rj}_{s_r=t_r(j-1)+1}\frac{F_{i_{ s_r}i_{ s_r}}}{t_r}\z_{\R, s_r}{\z}^{\top}_{\R, s_r}  \succeq \sum^{t_rj}_{s_r=t_r(j-1)+1}\frac{1}{t_r}\z_{\R,s_r}{\z}^{\top}_{\R, s_r} \succeq \frac{1}{2}   \V_{\R}^\top\V_{\R}, ~~~j=1,2,3,4.\nonumber
\end{align}
Since random sampling matrices $\S_{\C}$ and $\S_{\R}$ are independent, the events $\zeta_{\C,1}$, $\zeta_{\C,2}$, $\zeta_{\C,3}$, $\zeta_{\C,4}$, $\zeta_{\R,1}$, $\zeta_{\R,2}$, $\zeta_{\R,3}$, and $\zeta_{\R,4}$  are independent.   
% Consequently, for any index $s_c\in \{1,\ldots, m_c\}$,    there exists a block index $j=j(s)\in\{1,2,3,4\}$ such that
% \begin{enumerate}
%     \item the event $  \zeta_j$ is independent of $\z_{s}$; and 
%     \item conditioning  on $  \zeta_j$, it follows that  $ \Q\preceq \Q_{-s}\preceq 8\H^{-1}$.
% \end{enumerate} 
% and for  for any index  $s_r\in \{1,\ldots, m_r\}$    there exists a block index $j=j(s)\in\{1,2,3,4\}$ such that
% \begin{enumerate}
%     \item the event $  \zeta_{C,j}$ is independent of $\z_{s_c}$; and 
%     \item conditioning  on $  \zeta_{C,j}$, it follows that  $ \preceq \Q_{-s}\preceq 8\H^{-1}$.
% \end{enumerate} 

In this following, we begin by bounding $\|\EE_{\zeta}[\C \check{\U}\R]-\X\|^2_F-\|\C \U_{\CUR}\R-\X\|^2_F$.
Note that conditioned on the event $\zeta$, $\check\S_{\C}\W_\C$ and $\W_{\R}\Sigma_{\R}$  have full column rank, and $\Sigma_{\C}\V_{\C}^\top$  and $\V_{\R}^\top \check\S_{\R}^\top $ have full row rank.  Together with  \Cref{lemm:matrix_produc_pesedo}, we further rewrite
\begin{align*}
  \C\check{\U}\R& = \W_{\C}\Sigma_{\C}\V_{\C}^\top(\Sigma_{\C}\V_{\C}^\top)^\dagger(\check\S_{\C} \W_{\C} )^{\dagger}\check\S_{\C}\X\check\S_{\R}^\top  (\V_{\R}^\top \check\S_{\R}^\top  )^{\dagger} (\W_{\R}\Sigma_{\R})^\dagger \W_{\R}\Sigma_{\R}\V_{\R}^\top\\
 &= \W_{\C}(\check\S_{\C} \W_{\C} )^{\dagger}\check \S_{\C} \X\check\S_{\R}^\top(\V_{\R}^\top\check \S_{\R}^\top  )^{\dagger}   \V_{\R}^\top.%\equiv   \W_{\C}\check \M   \V_{\R}^\top.
\end{align*}

Considering 
\begin{align*}
  & \| \C\U_{\CUR}  \R-\X\|^2_F=\tr(\V_{\R}\V_{\R}^\top\X^\top\W_{\C}\W_{\C}^\top\W_{\C}\W_{\C}^\top\X\V_{\R}\V_{\R}^\top)-\tr(\V_{\R}\V_{\R}^\top\X^\top\W_{\C}\W_{\C}^\top\X)\\
     &-\tr(\X^\top\W_{\C}\W_{\C}^\top\X\V_{\R}\V_{\R}^\top)+\tr(\X^\top\X)=\tr(\V_{\R}\V_{\R}^\top\X^\top\W_{\C}\W_{\C}^\top\X\V_{\R}\V_{\R}^\top)\\
     &-\tr(\V_{\R}\V_{\R}^\top\V_{\R}\V_{\R}^\top\X^\top\W_{\C}\W_{\C}^\top\X)-\tr(\X^\top\W_{\C}\W_{\C}^\top\X\V_{\R}\V_{\R}^\top\V_{\R}\V_{\R}^\top)\\
     &+\tr(\X^\top\X)=\tr(\V_{\R}\V_{\R}^\top\X^\top\W_{\C}\W_{\C}^\top\X\V_{\R}\V_{\R}^\top)-\tr(\V_{\R}\V_{\R}^\top\X^\top\W_{\C}\W_{\C}^\top\X\V_{\R}\V_{\R}^\top)\\
     &-\tr(\V_{\R}\V_{\R}^\top\X^\top\W_{\C}\W_{\C}^\top\X\V_{\R}\V_{\R}^\top)+\tr(\X^\top\X)=\tr(\X^\top\X)-\tr(\V_{\R}\V_{\R}^\top\X^\top\W_{\C}\W_{\C}^\top\X\V_{\R}\V_{\R}^\top)\\
     &=\tr(\X^\top\X)-\tr(\V_{\R}\V_{\R}^\top\X^\top\W_{\C}\W_{\C}^\top\W_{\C}\W_{\C}^\top\X\V_{\R}\V_{\R}^\top)
\end{align*}
and 
\begin{align*}
     & \|\EE_{\zeta}[\C  \check{\U}\R]-\X\|^2_F=\tr(\EE_{\zeta}[\V_{\R} (\check \S_{\R}\V_{\R}  )^{\dagger}\check\S_{\R}\X^\top\check \S_{\C}^\top  ( \W_{\C}^\top\check\S_{\C}^\top )^{\dagger}\W_{\C}^\top ] \EE_{\zeta}[\W_{\C}(\check\S_{\C} \W_{\C} )^{\dagger}\check \S_{\C} \X\check\S_{\R}^\top(\V_{\R}^\top\check \S_{\R}^\top  )^{\dagger}   \V_{\R}^\top])\\
     &-\tr((\EE_{\zeta}[\V_{\R} (\check \S_{\R}\V_{\R}  )^{\dagger}\check\S_{\R}\X^\top\check \S_{\C}^\top  ( \W_{\C}^\top\check\S_{\C}^\top )^{\dagger}\W_{\C}^\top]\X)-\tr(\X^\top\EE_{\zeta}[\W_{\C}(\check\S_{\C} \W_{\C} )^{\dagger}\check \S_{\C} \X\check\S_{\R}^\top(\V_{\R}^\top\check \S_{\R}^\top  )^{\dagger}   \V_{\R}^\top])\\
     &+\tr(\X^\top\X)=\tr(\V_{\R}\EE_{\zeta}[ (\check \S_{\R}\V_{\R}  )^{\dagger}\check\S_{\R}\X^\top\check \S_{\C}^\top  ( \W_{\C}^\top\check\S_{\C}^\top )^{\dagger} ] \EE_{\zeta}[(\check\S_{\C} \W_{\C} )^{\dagger}\check \S_{\C} \X\check\S_{\R}^\top(\V_{\R}^\top\check \S_{\R}^\top  )^{\dagger}  ] \V_{\R}^\top)\\
     &-\tr(\V_{\R}\V_{\R}^\top\V_{\R}(\EE_{\zeta}[ (\check \S_{\R}\V_{\R}  )^{\dagger}\check\S_{\R}\X^\top\check \S_{\C}^\top  ( \W_{\C}^\top\check\S_{\C}^\top )^{\dagger}\W_{\C}^\top]\X)\\
     &-\tr(\X^\top\EE_{\zeta}[\W_{\C}(\check\S_{\C} \W_{\C} )^{\dagger}\check \S_{\C} \X\check\S_{\R}^\top(\V_{\R}^\top\check \S_{\R}^\top  )^{\dagger}   ]\V_{\R}^\top\V_{\R}\V_{\R}^\top)\\
     &+\tr(\X^\top\X)=\tr(\V_{\R}\EE_{\zeta}[ (\check \S_{\R}\V_{\R}  )^{\dagger}\check\S_{\R}\X^\top\check \S_{\C}^\top  ( \W_{\C}^\top\check\S_{\C}^\top )^{\dagger} \W_{\C}^\top] \EE_{\zeta}[\W_{\C}(\check\S_{\C} \W_{\C} )^{\dagger}\check \S_{\C} \X\check\S_{\R}^\top(\V_{\R}^\top\check \S_{\R}^\top  )^{\dagger}  ] \V_{\R}^\top)\\
     &-\tr(\V_{\R}\EE_{\zeta}[ (\check \S_{\R}\V_{\R}  )^{\dagger}\check\S_{\R}\X^\top\check \S_{\C}^\top  ( \W_{\C}^\top\check\S_{\C}^\top )^{\dagger}]\W_{\C}^\top\W_{\C}\W_{\C}^\top\X\V_{\R}\V_{\R}^\top)\\
     &-\tr(\V_{\R}\V_{\R}^\top\X^\top\W_{\C} \W_{\C}^\top\W_{\C}\EE_{\zeta}[(\check\S_{\C} \W_{\C} )^{\dagger}\check \S_{\C} \X\check\S_{\R}^\top(\V_{\R}^\top\check \S_{\R}^\top  )^{\dagger}   ]\V_{\R}^\top)+\tr(\X^\top\X),
\end{align*}
  we have
\begin{align*}
   & \Bias_{\zeta}(\check{\U})=   \|\EE_{\zeta}[\W_{\C}(\check\S_{\C} \W_{\C} )^{\dagger}\check \S_{\C} \X\check\S_{\R}^\top(\V_{\R}^\top\check \S_{\R}^\top  )^{\dagger}   \V_{\R}^\top]-\W_{\C}\W_{\C}^\top\X\V_{\R}\V_{\R}^\top\|^2_F
   \\
   &=\|\EE_{\zeta}[\W_{\C}(\check\S_{\C} \W_{\C} )^{\dagger}\check \S_{\C} (\I_n-\W_{\C}\W_{\C}^\top)\X(\I_p-\V_{\R}\V_{\R}^\top)\check\S_{\R}^\top(\V_{\R}^\top\check \S_{\R}^\top  )^{\dagger}   \V_{\R}^\top]\\
   &+\EE_{\zeta}[\W_{\C}\W_{\C}^\top\X(\I_p-\V_{\R}\V_{\R}^\top)\check\S_{\R}^\top(\V_{\R}^\top\check \S_{\R}^\top  )^{\dagger}   \V_{\R}^\top]+\EE_{\zeta}[\W_{\C}(\check\S_{\C} \W_{\C} )^{\dagger}\check \S_{\C} (\I_n-\W_{\C}\W_{\C}^\top)\X\V_{\R} \V_{\R}^\top]\|^2_F\\
   &\leq \|\EE_{\zeta}[\W_{\C}(\check\S_{\C} \W_{\C} )^{\dagger}\check \S_{\C} (\I_n-\W_{\C}\W_{\C}^\top)\X(\I_p-\V_{\R}\V_{\R}^\top)\check\S_{\R}^\top(\V_{\R}^\top\check \S_{\R}^\top  )^{\dagger}   \V_{\R}^\top]\|^2_F\\
   &+\|\EE_{\zeta}[\W_{\C}\W_{\C}^\top\X(\I_p-\V_{\R}\V_{\R}^\top)\check\S_{\R}^\top(\V_{\R}^\top\check \S_{\R}^\top  )^{\dagger}   \V_{\R}^\top]\|^2_F+\|\EE_{\zeta}[\W_{\C}(\check\S_{\C} \W_{\C} )^{\dagger}\check \S_{\C} (\I_n-\W_{\C}\W_{\C}^\top)\X\V_{\R} \V_{\R}^\top]\|^2_F\\
   &+2\|\EE_{\zeta}[\W_{\C}(\check\S_{\C} \W_{\C} )^{\dagger}\check \S_{\C} (\I_n-\W_{\C}\W_{\C}^\top)\X(\I_p-\V_{\R}\V_{\R}^\top)\check\S_{\R}^\top(\V_{\R}^\top\check \S_{\R}^\top  )^{\dagger}   \V_{\R}^\top]\|_F\\
   &\cdot\|\EE_{\zeta}[\W_{\C}\W_{\C}^\top\X(\I_p-\V_{\R}\V_{\R}^\top)\check\S_{\R}^\top(\V_{\R}^\top\check \S_{\R}^\top  )^{\dagger}   \V_{\R}^\top]\|_F\\
   &+2\|\EE_{\zeta}[\W_{\C}(\check\S_{\C} \W_{\C} )^{\dagger}\check \S_{\C} (\I_n-\W_{\C}\W_{\C}^\top)\X(\I_p-\V_{\R}\V_{\R}^\top)\check\S_{\R}^\top(\V_{\R}^\top\check \S_{\R}^\top  )^{\dagger}   \V_{\R}^\top]\|_F\\
   &\cdot\|\EE_{\zeta}[\W_{\C}(\check\S_{\C} \W_{\C} )^{\dagger}\check \S_{\C} (\I_n-\W_{\C}\W_{\C}^\top)\X\V_{\R} \V_{\R}^\top]\|_F\\
   &+2\|\EE_{\zeta}[\W_{\C}\W_{\C}^\top\X(\I_p-\V_{\R}\V_{\R}^\top)\check\S_{\R}^\top(\V_{\R}^\top\check \S_{\R}^\top  )^{\dagger}   \V_{\R}^\top]\|_F\|\EE_{\zeta}[\W_{\C}(\check\S_{\C} \W_{\C} )^{\dagger}\check \S_{\C} (\I_n-\W_{\C}\W_{\C}^\top)\X\V_{\R} \V_{\R}^\top]\|_F.
\end{align*}
Recalling the result in  \eqref{eq:expec_obliquepro_bias} in  \Cref{theo:de_sam_obliqueprojec},  we get
\begin{align*}  
&\|\EE_{\zeta}[\W_{\C}(\check\S_{\C} \W_{\C} )^{\dagger}\check \S_{\C} (\I_n-\W_{\C}\W_{\C}^\top)\X\V_{\R} \V_{\R}^\top]\|^2_F= O\left(\frac{(\log n)^2\theta^3_{\max}(\C)c^3}{(\loglog n)^2m_c^3}\right)\|(\I_n-\W_{\C}\W_{\C}^\top)\X\V_{\R} \V_{\R}^\top\|^2_F\\
&= O\left(\frac{(\log n)^2\theta^3_{\max}(\C)c^3}{(\loglog n)^2m_c^3}\right)\|(\I_n-\W_{\C}\W_{\C}^\top)\X\|^2_F \|\V_{\R} \V_{\R}^\top\|^2= O\left(\frac{(\log n)^2\theta^3_{\max}(\C)c^3}{(\loglog n)^2m_c^3}\right)\| \C\U_{\CUR}  \R-\X\|^2_F
\end{align*}
and 
\begin{align*}
  & \|\EE_{\zeta}[\W_{\C}\W_{\C}^\top\X(\I_p-\V_{\R}\V_{\R}^\top)\check\S_{\R}^\top(\V_{\R}^\top\check \S_{\R}^\top  )^{\dagger}   \V_{\R}^\top]\|^2_F =  \|\EE_{\zeta}[\V_{\R} (\check \S_{\R}\V_{\R}  )^{\dagger}\check\S_{\R}(\I_p-\V_{\R}\V_{\R}^\top)\X^\top\W_{\C}\W_{\C}^\top]\|^2_F\\
  &= O\left(\frac{(\log p)^2\theta^3_{\max}(\R^\top)r^3}{(\loglog p)^2m_r^3}\right)\|(\I_p-\V_{\R}\V_{\R}^\top)\X^\top\W_{\C}\W_{\C}^\top\|^2_F \\
  &=O\left(\frac{(\log p)^2\theta^3_{\max}(\R^\top)r^3}{(\loglog p)^2m_r^3}\right)\|\W_{\C}\W_{\C}^\top\X(\I_p-\V_{\R}\V_{\R}^\top)\|^2_F \\
  &= O\left(\frac{(\log p)^2\theta^3_{\max}(\R^\top)r^3}{(\loglog p)^2m_r^3}\right)\|\X(\I_p-\V_{\R}\V_{\R}^\top)\|^2_F  = O\left(\frac{(\log p)^2\theta^3_{\max}(\R^\top)r^3}{(\loglog p)^2m_r^3}\right)\| \C\U_{\CUR}  \R-\X\|^2_F,
\end{align*}
where we use the inequalities
\begin{align*}
   \|(\W_{\C}\W_{\C}^\top-\I_n)\X\|^2_F \leq  \| \C\U_{\CUR}  \R-\X\|^2_F,~~\|\X(\V_{\R}\V_{\R}^\top-\I_p)\|^2_F\leq  \| \C\U_{\CUR}  \R-\X\|^2_F.
\end{align*}
We also similarly obtain
\begin{align*}
 &\|\EE_{\zeta}[\W_{\C}(\check\S_{\C} \W_{\C} )^{\dagger}\check \S_{\C} (\I_n-\W_{\C}\W_{\C}^\top)\X(\I_p-\V_{\R}\V_{\R}^\top)\check\S_{\R}^\top(\V_{\R}^\top\check \S_{\R}^\top  )^{\dagger}   \V_{\R}^\top]\|^2_F  \\
 &=\|\EE_{\zeta_{\C}}[\W_{\C}(\check\S_{\C} \W_{\C} )^{\dagger}\check \S_{\C}] \EE_{\zeta_{\R}}[(\I_n-\W_{\C}\W_{\C}^\top)\X(\I_p-\V_{\R}\V_{\R}^\top)\check\S_{\R}^\top(\V_{\R}^\top\check \S_{\R}^\top  )^{\dagger}   \V_{\R}^\top]\|^2_F  \\
 &= O\left(\frac{(\log n)^2\theta^3_{\max}(\C)c^3}{(\loglog n)^2m_c^3}\right)\|\EE_{\zeta_{\R}}[ (\I_n-\W_{\C}\W_{\C}^\top)\X(\I_p-\V_{\R}\V_{\R}^\top)\check\S_{\R}^\top(\V_{\R}^\top\check \S_{\R}^\top  )^{\dagger}   \V_{\R}^\top]\|^2_F \\
 &= O\left(\frac{(\log n\log p)^2\theta^3_{\max}(\C)\theta^3_{\max}(\R^\top)c^3r^3}{(\loglog n\loglog p)^2m_c^3m_r^3}\right)\|(\I_n-\W_{\C}\W_{\C}^\top)\X(\I_p-\V_{\R}\V_{\R}^\top)\|^2_F\\
 &= O\left(\frac{(\log n\log p)^2\theta^3_{\max}(\C)\theta^3_{\max}(\R^\top)c^3r^3}{(\loglog n\loglog p)^2m_c^3m_r^3}\right)(\|\X(\I_p-\V_{\R}\V_{\R}^\top)\|^2_F+\|\W_{\C}\W_{\C}^\top\X(\I_p-\V_{\R}\V_{\R}^\top)\|^2_F\\
 &+2\|\W_{\C}\W_{\C}^\top\X(\I_p-\V_{\R}\V_{\R}^\top)\|_F\|\X(\I_p-\V_{\R}\V_{\R}^\top)\|_F)\\
 &= O\left(\frac{(\log n\log p)^2\theta^3_{\max}(\C)\theta^3_{\max}(\R^\top)c^3r^3}{(\loglog n\loglog p)^2m_c^3m_r^3}\right)\| \C\U_{\CUR}  \R-\X\|^2_F,
\end{align*}
where in the first equality we take $\zeta=\zeta_{\C}\cap\zeta_{\R}$, with $\zeta_{\C}$ and $\zeta_{\R}$ independent.
Thus, we conclude 
\begin{align*}
 &    \Bias_{\zeta}(\check{\U})= O\left(\frac{(\log n\log p)^2\theta^3_{\max}(\C)\theta^3_{\max}(\R^\top)c^3r^3}{(\loglog n\loglog p)^2m_c^3m_r^3}+\frac{(\log n)^2\theta^3_{\max}(\C)c^3}{(\loglog n)^2m_c^3}\sqrt{\frac{(\log p)^2\theta^3_{\max}(\R^\top)r^3}{(\loglog p)^2m_r^3}}\right.\\
    &\left.+\frac{(\log p)^2\theta^3_{\max}(\R^\top)r^3}{(\loglog p)^2m_r^3}\sqrt{\frac{(\log n)^2\theta^3_{\max}(\C)c^3}{(\loglog n)^2m_c^3}}+\frac{\log n\log p}{\loglog n\loglog p}\sqrt{\frac{\theta^3_{\max}(\C)\theta^3_{\max}(\R^\top)c^3r^3}{m_c^3m_r^3}}\right.\\
    &\left.+\frac{(\log n)^2\theta^3_{\max}(\C)c^3}{(\loglog n)^2m_c^3}+\frac{(\log p)^2\theta^3_{\max}(\R^\top)r^3}{(\loglog p)^2m_r^3}\right)\| \C\U_{\CUR}  \R-\X\|^2_F. 
\end{align*}

Subsequently, we proceed to bound the another term $\Var_{\zeta}(\check{\U})-\Delta_{\mathrm{CUR}}$. We first rewrite
\begin{align*}
 &\Var_{\zeta}(\check{\U}) =\EE[\|\W_{\C}(\check\S_{\C} \W_{\C} )^{\dagger}\check \S_{\C} \X\check\S_{\R}^\top(\V_{\R}^\top\check \S_{\R}^\top  )^{\dagger}   \V_{\R}^\top-\W_{\C}\W_{\C}^\top\X\V_{\R}\V_{\R}^\top\|^2_F]  \\
  &=\EE_{\zeta}[ \|\W_{\C}(\check\S_{\C} \W_{\C} )^{\dagger}\check \S_{\C} (\I_n-\W_{\C}\W_{\C}^\top)\X(\I_p-\V_{\R}\V_{\R}^\top)\check\S_{\R}^\top(\V_{\R}^\top\check \S_{\R}^\top  )^{\dagger}   \V_{\R}^\top\\
   &+\W_{\C}\W_{\C}^\top\X(\I_p-\V_{\R}\V_{\R}^\top)\check\S_{\R}^\top(\V_{\R}^\top\check \S_{\R}^\top  )^{\dagger}   \V_{\R}^\top+\W_{\C}(\check\S_{\C} \W_{\C} )^{\dagger}\check \S_{\C} (\I_n-\W_{\C}\W_{\C}^\top)\X\V_{\R} \V_{\R}^\top\|^2_F]\\
   &\leq \EE_{\zeta}[\|\W_{\C}(\check\S_{\C} \W_{\C} )^{\dagger}\check \S_{\C} (\I_n-\W_{\C}\W_{\C}^\top)\X(\I_p-\V_{\R}\V_{\R}^\top)\check\S_{\R}^\top(\V_{\R}^\top\check \S_{\R}^\top  )^{\dagger}   \V_{\R}^\top\|^2_F]\\
   &+\EE_{\zeta}[\|\W_{\C}\W_{\C}^\top\X(\I_p-\V_{\R}\V_{\R}^\top)\check\S_{\R}^\top(\V_{\R}^\top\check \S_{\R}^\top  )^{\dagger}   \V_{\R}^\top\|^2_F]+\EE_{\zeta}[\|\W_{\C}(\check\S_{\C} \W_{\C} )^{\dagger}\check \S_{\C} (\I_n-\W_{\C}\W_{\C}^\top)\X\V_{\R} \V_{\R}^\top\|^2_F]\\
   &+2\EE_{\zeta}[\|\W_{\C}(\check\S_{\C} \W_{\C} )^{\dagger}\check \S_{\C} (\I_n-\W_{\C}\W_{\C}^\top)\X(\I_p-\V_{\R}\V_{\R}^\top)\check\S_{\R}^\top(\V_{\R}^\top\check \S_{\R}^\top  )^{\dagger}   \V_{\R}^\top\|_F\\
   &\cdot\|\W_{\C}\W_{\C}^\top\X(\I_p-\V_{\R}\V_{\R}^\top)\check\S_{\R}^\top(\V_{\R}^\top\check \S_{\R}^\top  )^{\dagger}   \V_{\R}^\top\|_F]\\
   &+2\EE_{\zeta}[\|\W_{\C}(\check\S_{\C} \W_{\C} )^{\dagger}\check \S_{\C} (\I_n-\W_{\C}\W_{\C}^\top)\X(\I_p-\V_{\R}\V_{\R}^\top)\check\S_{\R}^\top(\V_{\R}^\top\check \S_{\R}^\top  )^{\dagger}   \V_{\R}^\top\|_F\\
   &\cdot\|\W_{\C}(\check\S_{\C} \W_{\C} )^{\dagger}\check \S_{\C} (\I_n-\W_{\C}\W_{\C}^\top)\X\V_{\R} \V_{\R}^\top\|_F]\\
   &+2\EE_{\zeta}[\|\W_{\C}\W_{\C}^\top\X(\I_p-\V_{\R}\V_{\R}^\top)\check\S_{\R}^\top(\V_{\R}^\top\check \S_{\R}^\top  )^{\dagger}   \V_{\R}^\top\|_F\|\W_{\C}(\check\S_{\C} \W_{\C} )^{\dagger}\check \S_{\C} (\I_n-\W_{\C}\W_{\C}^\top)\X\V_{\R} \V_{\R}^\top\|_F]\\
    &\leq \EE_{\zeta}[\|\W_{\C}(\check\S_{\C} \W_{\C} )^{\dagger}\check \S_{\C} (\I_n-\W_{\C}\W_{\C}^\top)\X(\I_p-\V_{\R}\V_{\R}^\top)\check\S_{\R}^\top(\V_{\R}^\top\check \S_{\R}^\top  )^{\dagger}   \V_{\R}^\top\|^2_F]\\
   &+\EE_{\zeta}[\|\W_{\C}\W_{\C}^\top\X(\I_p-\V_{\R}\V_{\R}^\top)\check\S_{\R}^\top(\V_{\R}^\top\check \S_{\R}^\top  )^{\dagger}   \V_{\R}^\top\|^2_F]+\EE_{\zeta}[\|\W_{\C}(\check\S_{\C} \W_{\C} )^{\dagger}\check \S_{\C} (\I_n-\W_{\C}\W_{\C}^\top)\X\V_{\R} \V_{\R}^\top\|^2_F]\\
   &+2\sqrt{\EE_{\zeta}[\|\W_{\C}(\check\S_{\C} \W_{\C} )^{\dagger}\check \S_{\C} (\I_n-\W_{\C}\W_{\C}^\top)\X(\I_p-\V_{\R}\V_{\R}^\top)\check\S_{\R}^\top(\V_{\R}^\top\check \S_{\R}^\top  )^{\dagger}   \V_{\R}^\top\|_F^2]}\\
   &\cdot\sqrt{\EE_{\zeta}[\|\W_{\C}\W_{\C}^\top\X(\I_p-\V_{\R}\V_{\R}^\top)\check\S_{\R}^\top(\V_{\R}^\top\check \S_{\R}^\top  )^{\dagger}   \V_{\R}^\top\|^2_F]}\\
   &+2\sqrt{\EE_{\zeta}[\|\W_{\C}(\check\S_{\C} \W_{\C} )^{\dagger}\check \S_{\C} (\I_n-\W_{\C}\W_{\C}^\top)\X(\I_p-\V_{\R}\V_{\R}^\top)\check\S_{\R}^\top(\V_{\R}^\top\check \S_{\R}^\top  )^{\dagger}   \V_{\R}^\top\|_F^2]}\\
   &\cdot \sqrt{\EE_{\zeta}[\|\W_{\C}(\check\S_{\C} \W_{\C} )^{\dagger}\check \S_{\C} (\I_n-\W_{\C}\W_{\C}^\top)\X\V_{\R} \V_{\R}^\top\|^2_F]}\\
   &+2\sqrt{\EE_{\zeta}[\|\W_{\C}\W_{\C}^\top\X(\I_p-\V_{\R}\V_{\R}^\top)\check\S_{\R}^\top(\V_{\R}^\top\check \S_{\R}^\top  )^{\dagger}   \V_{\R}^\top\|_F^2]}\\
   &\cdot \sqrt{\EE_{\zeta}[\|\W_{\C}(\check\S_{\C} \W_{\C} )^{\dagger}\check \S_{\C} (\I_n-\W_{\C}\W_{\C}^\top)\X\V_{\R} \V_{\R}^\top\|^2_F]}.
\end{align*}
Using the result \eqref{eq:expec_obliqueprovariance} in \Cref{theo:de_sam_obliqueprojec}, it follows that
\begin{align*}  
&\EE_{\zeta}[\|\W_{\C}(\check\S_{\C} \W_{\C} )^{\dagger}\check \S_{\C} (\I_n-\W_{\C}\W_{\C}^\top)\X\V_{\R} \V_{\R}^\top\|^2_F]\leq \EE_{\zeta}[\|\W_{\C}(\check\S_{\C} \W_{\C} )^{\dagger}\check \S_{\C} (\I_n-\W_{\C}\W_{\C}^\top)\X\|^2_F]\\
&= O\left(\frac{\log n}{\loglog n}\sqrt{\frac{\theta^3_{\max}(\C)c^3}{m_c^3}}\right)\|(\I_n-\W_{\C}\W_{\C}^\top)\X\|^2_F+\left\|\diag\Big\{\frac{(\ell_{i}(\C))^{1/2}}{\sqrt{m_c\pi_{\C, i}}}\Big\}^{n}_{i=1}(\I_n-\W_{\C}\W_{\C}^\top)\X\right\|^2_F\\
&= O\left(\frac{\log n}{\loglog n}\sqrt{\frac{\theta^3_{\max}(\C)c^3}{m_c^3}}\right)\| \C\U_{\CUR}  \R-\X\|^2_F+\left\|\diag\Big\{\frac{(\ell_{i}(\C))^{1/2}}{\sqrt{m_c\pi_{\C, i}}}\Big\}^{n}_{i=1}(\I_n-\W_{\C}\W_{\C}^\top)\X\right\|^2_F\\
  &= O\left(\frac{\log n}{\loglog n}\sqrt{\frac{\theta^3_{\max}(\C)c^3}{m_c^3}}+\frac{\theta_{\max}(\C)c}{m_c}\right)\| \C\U_{\CUR}  \R-\X\|^2_F
\end{align*}
and 
\begin{align*}
  & \EE_{\zeta}[\|\W_{\C}\W_{\C}^\top\X(\I_p-\V_{\R}\V_{\R}^\top)\check\S_{\R}^\top(\V_{\R}^\top\check \S_{\R}^\top  )^{\dagger}   \V_{\R}^\top\|^2_F ]= \EE_{\zeta}[ \|\V_{\R} (\check \S_{\R}\V_{\R}  )^{\dagger}\check\S_{\R}(\I_p-\V_{\R}\V_{\R}^\top)\X^\top\W_{\C}\W_{\C}^\top\|^2_F]\\
  &\leq \EE_{\zeta}[ \|\V_{\R} (\check \S_{\R}\V_{\R}  )^{\dagger}\check\S_{\R}(\I_p-\V_{\R}\V_{\R}^\top)\X^\top\|^2_F]\\
  &= O\left(\frac{\log p}{\loglog p}\sqrt{\frac{\theta^3_{\max}(\R^\top)r^3}{m_r^3}}\right)\|(\I_p-\V_{\R}\V_{\R}^\top)\X^\top\|^2_F +\left\|\diag\Big\{\frac{(\ell_{j}(\R^\top))^{1/2}}{\sqrt{m_r\pi_{\R, j}}}\Big\}^{p}_{j=1}(\I_p-\V_{\R}\V_{\R}^\top)\X^\top\right\|^2_F\\
  &= O\left(\frac{\log p}{\loglog p}\sqrt{\frac{\theta^3_{\max}(\R^\top)r^3}{m_r^3}}\right)\| \C\U_{\CUR}  \R-\X\|^2_F +\left\|\diag\Big\{\frac{(\ell_{j}(\R^\top))^{1/2}}{\sqrt{m_r\pi_{\R, j}}}\Big\}^{p}_{j=1}(\I_p-\V_{\R}\V_{\R}^\top)\X^\top\right\|^2_F\\
 &=   O\left(\frac{\log p}{\loglog p}\sqrt{\frac{\theta^3_{\max}(\R^\top)r^3}{m_r^3}}+\frac{\theta_{\max}(\R^\top)r}{m_r}\right)\| \C\U_{\CUR}  \R-\X\|^2_F,
\end{align*}
where we apply the inequalities
\begin{align*}
   \|(\W_{\C}\W_{\C}^\top-\I_n)\X\|^2_F \leq  \| \C\U_{\CUR}  \R-\X\|^2_F,~~\|\X(\V_{\R}\V_{\R}^\top-\I_p)\|^2_F\leq  \| \C\U_{\CUR}  \R-\X\|^2_F.
\end{align*}
Recalling  \eqref{eq:expec_obliqueprovariance} in \Cref{theo:de_sam_obliqueprojec} again, we similarly obtain, 
\begin{align*}
 &\EE_{\zeta}[\|\W_{\C}(\check\S_{\C} \W_{\C} )^{\dagger}\check \S_{\C} (\I_n-\W_{\C}\W_{\C}^\top)\X(\I_p-\V_{\R}\V_{\R}^\top)\check\S_{\R}^\top(\V_{\R}^\top\check \S_{\R}^\top  )^{\dagger}   \V_{\R}^\top\|^2_F ] \\
 &= O\left(\frac{\log n}{\loglog n}\sqrt{\frac{\theta^3_{\max}(\C)c^3}{m_c^3}}\right)\EE_{\zeta}[\|\X(\I_p-\V_{\R}\V_{\R}^\top)\check\S_{\R}^\top(\V_{\R}^\top\check \S_{\R}^\top  )^{\dagger}   \V_{\R}^\top\|^2_F ]\\
 &+\EE_{\zeta}\left[\left\|\diag\Big\{\frac{(\ell_{i}(\C))^{1/2}}{\sqrt{m_c\pi_{\C, i}}}\Big\}^{n}_{i=1} (\I_n-\W_{\C}\W_{\C}^\top)\X(\I_p-\V_{\R}\V_{\R}^\top)\check\S_{\R}^\top(\V_{\R}^\top\check \S_{\R}^\top  )^{\dagger}   \V_{\R}^\top\right\|^2_F\right ]\\
 &= O\left(\frac{\log n}{\loglog n}\sqrt{\frac{\theta^3_{\max}(\C)c^3}{m_c^3}}\right)\EE_{\zeta}[\|\X(\I_p-\V_{\R}\V_{\R}^\top)\check\S_{\R}^\top(\V_{\R}^\top\check \S_{\R}^\top  )^{\dagger}   \V_{\R}^\top\|^2_F ]\\
 &+\frac{4\theta_{\max}(\C)c}{m_c} \EE_{\zeta}[\|\X(\I_p-\V_{\R}\V_{\R}^\top)\check\S_{\R}^\top(\V_{\R}^\top\check \S_{\R}^\top  )^{\dagger}   \V_{\R}^\top\|^2_F ].
 % \leq 2 \EE_{\zeta}[\|\X(\I_p-\V_{\R}\V_{\R}^\top)\check\S_{\R}^\top(\V_{\R}^\top\check \S_{\R}^\top  )^{\dagger}   \V_{\R}^\top\|^2_F ].
\end{align*}
We further similarly  gain, 
\begin{align*}
    &\EE_{\zeta}[\|\W_{\C}(\check\S_{\C} \W_{\C} )^{\dagger}\check \S_{\C} (\I_n-\W_{\C}\W_{\C}^\top)\X(\I_p-\V_{\R}\V_{\R}^\top)\check\S_{\R}^\top(\V_{\R}^\top\check \S_{\R}^\top  )^{\dagger}   \V_{\R}^\top\|^2_F ] \\
    &= O\left(\frac{\log p}{\loglog p}\sqrt{\frac{\theta^3_{\max}(\R^\top)r^3}{m_r^3}}\right) \EE_{\zeta}[\|\W_{\C}(\check\S_{\C} \W_{\C} )^{\dagger}\check \S_{\C} (\I_n-\W_{\C}\W_{\C}^\top)\X\|^2_F ]\\
    &+ \EE_{\zeta}\left[\left\|\W_{\C}(\check\S_{\C} \W_{\C} )^{\dagger}\check \S_{\C} (\I_n-\W_{\C}\W_{\C}^\top)\X(\I_p-\V_{\R}\V_{\R}^\top)\diag\Big\{\frac{(\ell_{j}(\R^\top))^{1/2}}{\sqrt{m_r\pi_{\R, j}}}\Big\}^{p}_{j=1}\right\|^2_F\right]\\
    &= O\left(\frac{\log p}{\loglog p}\sqrt{\frac{\theta^3_{\max}(\R^\top)r^3}{m_r^3}}\right) \EE_{\zeta}[\|\W_{\C}(\check\S_{\C} \W_{\C} )^{\dagger}\check \S_{\C} (\I_n-\W_{\C}\W_{\C}^\top)\X\|^2_F ]\\
    &+ \frac{4\theta_{\max}(\R^\top)r}{m_r} \EE_{\zeta}[\|\W_{\C}(\check\S_{\C} \W_{\C} )^{\dagger}\check \S_{\C} (\I_n-\W_{\C}\W_{\C}^\top)\X\|^2_F].
    % &\leq 2 \EE_{\zeta}[\left\|\W_{\C}(\check\S_{\C} \W_{\C} )^{\dagger}\check \S_{\C} (\I_n-\W_{\C}\W_{\C}^\top)\X\right\|^2_F]
\end{align*}
Putting the above together, we conclude
\begin{align*}
  &\Var_{\zeta}(\check{\U})=O\left(\frac{\log n}{\loglog n}\sqrt{\frac{\theta^3_{\max}(\C)c^3}{m_c^3}} +\frac{\log p}{\loglog p}\sqrt{\frac{\theta^3_{\max}(\R^\top)r^3}{m_r^3}}\right)\| \C\U_{\CUR}  \R-\X\|^2_F\\
  &+\left\|\diag\Big\{\frac{(\ell_{i}(\C))^{1/2}}{\sqrt{m_c\pi_{\C, i}}}\Big\}^{n}_{i=1}(\I_n-\W_{\C}\W_{\C}^\top)\X\right\|^2_F+\left\|\diag\Big\{\frac{(\ell_{j}(\R^\top))^{1/2}}{\sqrt{m_r\pi_{\R, j}}}\Big\}^{p}_{j=1}(\I_p-\V_{\R}\V_{\R}^\top)\X^\top\right\|^2_F\\
  &+ O\left(\sqrt[4]{\frac{(\log n)^2\theta^3_{\max}(\C)c^3}{(\loglog n)^2m_c^3}} \right)\| \C\U_{\CUR}  \R-\X\|_F \left\|\diag\Big\{\frac{(\ell_{j}(\R^\top))^{1/2}}{\sqrt{m_r\pi_{\R, j}}}\Big\}^{p}_{j=1}(\I_p-\V_{\R}\V_{\R}^\top)\X^\top\right\|_F\\
  &+  O\left(\sqrt[4]{\frac{(\log p)^2\theta^3_{\max}(\R^\top)r^3}{(\loglog p)^2m_r^3}}\right)\| \C\U_{\CUR}  \R-\X\|_F\left\|\diag\Big\{\frac{(\ell_{i}(\C))^{1/2}}{\sqrt{m_c\pi_{\C, i}}}\Big\}^{n}_{i=1}(\I_n-\W_{\C}\W_{\C}^\top)\X\right\|_F\\
  &+2\left\|\diag\Big\{\frac{(\ell_{j}(\R^\top))^{1/2}}{\sqrt{m_r\pi_{\R, j}}}\Big\}^{p}_{j=1}(\I_p-\V_{\R}\V_{\R}^\top)\X^\top\right\|_F\left\|\diag\Big\{\frac{(\ell_{i}(\C))^{1/2}}{\sqrt{m_c\pi_{\C, i}}}\Big\}^{n}_{i=1}(\I_n-\W_{\C}\W_{\C}^\top)\X\right\|_F.
\end{align*}
This completes the proof.

\subsection{Special cases of Theorem~\ref{theo:cur_bias_variance} under leverage score sampling and SRHT}\label{subsec:cur_lev_srht}

In this section, we establish fine-grained accuracy characterizations for fast CUR  under leverage-score sampling and the SRHT.
The proofs follow directly by combining the arguments from \Cref{theo:de_sam_ols}, \Cref{coro:sub_ols_lev}, and \Cref{coro:debiasing_srht}, and are therefore omitted.

\begin{corollary}[Bias--variance characterizations for fsat  CUR  under  leverage score sampling]\label{coro:sub_cur_lev}
      Under the settings and notations of \Cref{theo:cur_bias_variance},  standard random sampling matrices $\S_{\C}\in \RR^{m_{c} \times n}$  with sampling probabilities
$\pi_{\C,i}\in [\ell_{i}(\C)/(c \theta_{\max}(\C)), \ell_{i}(\C)/(c \theta_{\min}(\C))]$ with $\theta_{\min}(\C) \in [1/2,1]$    and $\S_{\R}\in \RR^{m_{r} \times p} $   with  sampling probabilities
$\pi_{\R,i}\in [\ell_{i}(\R^\top)/(r \theta_{\max}(\R^\top)), \ell_{i}(\R^\top)/(r \theta_{\min}(\R^\top))]$ with $\theta_{\min}(\R^\top)\in [1/2,1]$, there exists $C > 0$, $\alpha_c>0$, $\alpha_r>0$, independent of $n, p, c, r$ so that for $m_{c} \geq C\theta_{\max}(\C)c\log (c/\delta)$, $m_{r} \geq C \theta_{\max}(\R^\top)r\log (r/\delta)$,  $\delta\leq \min\{n^{-(3+\alpha_c)}, p^{-(3+\alpha_r)}\} $, when conditioned on an event $\zeta$ that holds with the probability at least $1-\delta$, the standard sampled estimator $\tilde{\U}$ in \eqref{eq:faster_U} satisfies 
\begin{align*}
 \Bias_{\zeta}(\tilde{\U}) = \epsilon^2 \cdot L(\U_{\CUR}),~\Var_{\zeta}(\tilde{\U}) = \Delta_{\CUR} + \epsilon \cdot L(\U_{\CUR}),
\end{align*} 
with  $\epsilon=O(\sqrt{(\log n/\loglog n)^2\cdot \theta^3_{\max}(\C) c^3/m_c^3}+\epsilon_{\theta,c}\cdot\theta_{\max}(\C)c/m_c+\sqrt{(\log p/\loglog p)^2\cdot \theta^3_{\max}(\R^\top) r^3/m_r^3}+\epsilon_{\theta,r}\cdot\theta_{\max}(\R^\top)r/m_r)$, $\epsilon_{\theta,c}=\max\{\theta^{-1}_{\min}(\C)-1,1 - \theta^{-1}_{\max}(\C)\} $,  and  $\epsilon_{\theta,r}=\max\{\theta^{-1}_{\min}(\R^\top)-1,1 - \theta^{-1}_{\max}(\R^\top)\} $.
Here, $(\theta_{\max}(\C),\theta_{\min}(\C))$ and $(\theta_{\max}(\R^\top),\theta_{\min}(\R^\top))$ in \Cref{def:approx_factor} denote the maximum and minimum importance sampling approximation
factors for $\C$ and $\R^\top$, respectively. 
\end{corollary}

\begin{corollary}[Bias--variance characterizations  for fast CUR decomposition under SRHT]\label{coro:debiasing_cur_srht} 
      Under the setting and notations of \Cref{theo:cur_bias_variance}, 
      let $\tilde{\C}_{\SRHT}= \S_{\C}\C_{\HD} \in \RR^{m_c \times c}$ the SRHT of $\C$ with $\C_{\HD}=\H_n \D_n\C/\sqrt n$,  $\tilde{\R}^\top_{\SRHT}= \S_{\R}\R^\top_{\HD}\in \RR^{m_r \times r}$ the SRHT of $\R^\top$ with $ \R^\top_{\HD}=\H_p \D_p \R^\top/\sqrt p $, $\tilde \X_{\SRHT_{\C,\R} }=\S_{\C} \H_n \D_n \X\D_p\H_p  \S^\top_{\R}/\sqrt{np}$.   Assume that  $ \min_{i,\|\cc_{\HD,i}\|>0}\|\cc_{\HD,i}\|^2 /\|\C_{\HD}\|_F^2 \geq n^{-{\alpha_c}}$  and $ \min_{j,\|\rr_{\HD,j}\|>0}\|\rr_{\HD,j}\|^2 /\|\R_{\HD}\|_F^2 \geq p^{-{\alpha_r}}$ for some constants $\alpha_c>0$ and $\alpha_r>0$.
       Then, there exists $C> 0$,     $\max\{n \exp(-c), p \exp(-r)\} < \delta < \min\{n^{-(3+\alpha_c)}, p^{-(3+\alpha_r)}\} $ such that
for   $m_{c} \geq Cc\log (c/\delta)$, $m_{r} \geq C  r\log (r/\delta)$, when  conditioned on an event $\zeta$ that holds with the probability at least $1-\delta$, the  sampled estimator  $\tilde{\U}=\tilde{\C}_{\SRHT}^\dagger 
   \tilde \X_{\SRHT_{\C,\R} }
   \tilde{\R}_{\SRHT}^\dagger$ satisfies 
\begin{align*}
\Bias_{\zeta}(\tilde{\U}) = \epsilon^2 \cdot L(\U_{\CUR}),~\Var_{\zeta}(\tilde{\U}) = \Delta_{\CUR} + \epsilon \cdot L(\U_{\CUR}),
\end{align*}  
with $ \epsilon=  O(\sqrt{(\log n/\loglog n)^2\cdot c^3/m_c^3}+ \sqrt{c\log(n/\delta)}/m_c+\sqrt{(\log p/\loglog p)^2\cdot r^3/m_r^3}+\sqrt{r\log(p/\delta)}/m_r)$, 
 where we define
 $\Delta_{\CUR} = (\Delta_{\CUR, 1} + \Delta_{\CUR, 2})^2+\sqrt{\epsilon \cdot L(\U_{\CUR})} (\Delta_{\CUR, 1} + \Delta_{\CUR, 2}) $, with $\Delta_{\CUR, 1}= \|\diag\{\sqrt{ \ell_{i}(\C_{\HD}  )/m_c\pi_{\C, i}} \}^{n}_{i=1}\B_1 \|_F$ and $\Delta_{\CUR, 2}= \|\diag \{\sqrt{\ell_{j}( \R_{\HD}^\top )/m_r\pi_{\R, j}} \}^p_{j=1}\B_2\|_F$ for $\B_1=\H_n \D_n(\I_n-\W_{\C}\W_{\C}^\top)\X /\sqrt n \in  \RR^{n \times p}$ and   $\B_2= \H_p \D_p (\I_p-\V_{\R}\V_{\R}^\top)\X^\top/\sqrt p  \in  \RR^{p \times n}$.
Here, 
$\S_{\C}\in \RR^{m_c \times n}$ and $\S_{\R}\in \RR^{m_r \times p}$ are the uniform random sampling matrices applied to $\C_{\HD}  $ and $\R_{\HD}^\top  $, respectively.
\end{corollary}

\Cref{table:compare_cur} summarizes the existing upper bounds on the bias of  fast CUR decomposition established in prior work under uniform sampling, leverage score sampling, and  SRHT~\cite{wang2016towards,ye2019fast}, and contrasts them with our results on debiased fast CUR decomposition.
We see that the proposed debiased approach consistently achieves a substantially smaller bias across all three sampling schemes.

\begin{table*}[htb!] 
  \caption{ Bias characterizations for classical fast CUR  established in previous efforts versus that for the proposed debiased fast CUR in \Cref{theo:cur_bias_variance}, 
 under uniform sampling (\textbf{UNI}), (exact/ approximated) leverage-score sampling (\textbf{Lev}) with   $\theta_{\min}(\C)/(2\theta_{\min}(\C)-1)\leq \theta_{\max}(\C) \leq  1/(1 - \sqrt{\theta_{\max}(\C)c/m_c})$ and  $\theta_{\min}(\R^\top)/(2\theta_{\min}(\R^\top)-1)\leq \theta_{\max}(\R^\top) \leq  1/(1 - \sqrt{\theta_{\max}(\R^\top)r/m_r})$, and \textbf{SRHT} with $c>\log(n/\delta)$ and $r> \log(p/\delta)$.   We report the regime $m_c/m_r \leq \min\{\nu_1, \nu_2,\nu_3, \nu_4,\nu_5, \nu_6\}$, where  $\nu_1=c/r$, $\nu_2=c\theta_{\max}(\C) / (r\theta_{\max}(\R^\top))$, $\nu_3=c\theta_{\max}(\C)(\log n\loglog p )^{2/3}/(r\theta_{\max}(\R^\top)(\log p \loglog n )^{2/3})$,   $\nu_4=\sqrt{cn(\log n)^2/( rp(\log p)^2)}$, $\nu_5=c(\log n\loglog p )^{2/3}/(r(\log p \loglog n )^{2/3})$, $\nu_6=\sqrt{c\log(n/\delta)/( r\log(p/\delta))}$;
  the complementary regime, where the ratio $m_c/m_r $ exceeds this minimum, yields analogous conclusions and is omitted for brevity.
} %with mean accuracy and standard deviation
  \label{table:compare_cur}
  \centering
  \begin{threeparttable}
  \begin{tabular}{lcccccccc}
    \toprule
     & \textbf{UNI} &   \textbf{Lev} & \textbf{SRHT} \\
    \midrule    
\cite{wang2016towards}   & $O(\theta_{\max}(\C)^2c^2q/m_c^2)$  &  $O( c^2q/m_c^2)$   &   $O(cn(\log n)^2/m_c^2)$ \\
\cite{ye2019fast}  & ---  &  $O( \rho c/m_c)$   &   $O( \rho c/m_c)$    \\
\underline{\textbf{This work}}  &  $O\Big(\frac{(\log n)^2\theta^3_{\max}(\C)c^3}{(\loglog n)^2m^3_c}\Big)$   &  $O\Big(\frac{(\log n)^2c^3}{(\loglog n)^2m^3_c}\Big)$   &   $ O\Big(\frac{(\log n)^2c^3}{\loglog n)^2m^3_c}+\frac{c\log(n/\delta)}{m_c^2}\Big)$ \\
    \bottomrule
  \end{tabular}
\begin{tablenotes}[flushleft]
\footnotesize
\item[a]
Here, $q$ in the error bounds of \cite{wang2016towards} is defined as $q=\min\{n,p\}$. Moreover,  the bias error bounds in \cite{ye2019fast} are established under the additional condition $\rho=\bigl(\|(\I_n-\C \C^{\dagger}) \X \R^{\dagger}\R\|_F
      + \|\C \C^{\dagger}\X(\I_p-\R^{\dagger}\R)\|_F\bigr)^2/\|\C \C^{\dagger}\X \R^{\dagger}\R - \X\|_F^2
\geq\epsilon.
$
\end{tablenotes}
\end{threeparttable}
\end{table*}

\section{Additional numerical experiments and implementation details}
\label{sec:imple_detail_nmerical_exper}

In this section, we first provide implementation details for the numerical experiments on subsampled OLS and fast CUR decomposition in \Cref{sec:sketch_matrix} and \Cref{sec:Datasets}. We then present additional numerical results for subsampled OLS in \Cref{subsec:addition_num_ols}, followed by numerical results for fast CUR decomposition in \Cref{subsec:numerical_cur}.

\subsection{Sketching matrices}
\label{sec:sketch_matrix}

Given a data matrix $\X\in\RR^{n\times p}$, leverage scores are estimated following the procedure described in \cite{Drineas2012fast}. 
The LESS-uniform sketching matrix is generated according to the construction detailed in \cite[Section~E.1]{derezinski2021newtonless}. Implementation specifics, including those for the SRHT, follow~\cite{derezinski2021newtonless}.

Experiments were conducted on a server with an AMD EPYC 7452 32-Core Processor, 256GB RAM, and  NVIDIA GeForce RTX 3090 GPUs. Code is publicly available
at \url{https://github.com/chengmeiniu/debiased-oblique-projections.git}.

\subsection{Datasets}\label{sec:Datasets}

 The Flight Delay dataset used in \Cref{subsec:addition_num_ols} is obtained from the U.S. Department of Transportation. The dataset contains 539747 U.S. domestic weekday flights in January 2025, with five recorded variables per flight: arrival delay (difference in minutes between scheduled and actual arrival time, and early arrivals show negative numbers), arrival taxi in time (in minutes), departure taxi out time (in minutes), departure delays (difference in minutes between scheduled and actual departure time, and early departures show negative numbers), and computer reservation system based elapsed time of the flight (in minutes; a measure for the distance of the flight). Following \cite{ma2022asymptotic}, we take arrival delay as the response and the remaining four variables as predictors, and further augment the design with all quadratic terms and pairwise interactions, yielding 14 predictors in total.
All variables are standardized to zero mean and unit variance. In the experiments, we randomly select $n=2^{13}$ samples from the full dataset to construct the regression design matrix.

The MSD (Million Song Year Prediction) dataset used in \Cref{subsec:num_exper_ols} consists of 515,345 songs released between 1922 and 2011. Each song is represented by multiple audio segments, and each segment is described by 12 timbre features that capture perceptual  properties such as brightness and spectral flatness.
Our goal is to predict the release year using the complete set of timbre-based descriptors. Following~\cite{ma2022asymptotic}, we model the logarithm of the release year as the response variable and use all derived timbre features as predictors. All predictors are standardized to zero mean and unit variance. In the experiments, we randomly sample $n=2^{14}$ observations from the full dataset to form the regression design matrix. 

The CIFAR-10 dataset used in \Cref{subsec:numerical_cur} contains 60,000 images of size $32\times 32\times 3$.  Each image is flattened into a vector.
% , centered to have zero mean, and normalized so that its entries lie in  $[ 1, 1]$.
 We construct the data matrix using  $n = 2^{13}$ samples and $p =2^{11}$  features.

The ImageNet dataset used in \Cref{subsec:numerical_cur} is based on  ImageNet-64, a downsampled version of ImageNet. We use images from \texttt{Imagenet64\_train\_part1}, specifically \texttt{train\_data\_batch\_1}.  Each image has resolution $64\times 64\times 3$ and is flattened into a vector representation. 
% The resulting vectors are centered to have zero mean and normalized so that their entries lie in $[1, 1]$.
We construct the data matrix using $n = 2^{14}$ samples and $p =2^{11}$ features.

\subsection{Additional numerical experiments for subsampled OLS}\label{subsec:addition_num_ols}

In this section, we use the Flight Delay  (Airline) dataset, obtained from the website of the U.S. Department of Transportation\footnote{Airline On-Time Performance Data, \url{https://www.transtats.bts.gov/DL_SelectFields.asp?Table_ID=236}.}, to provide additional numerical evidence supporting our theoretical findings on classical and debiased subsampled OLS. The data matrix $\X \in \RR^{n \times p}$ is constructed from the Airline dataset, and $\y\in \RR^n$ denotes the corresponding response vector.

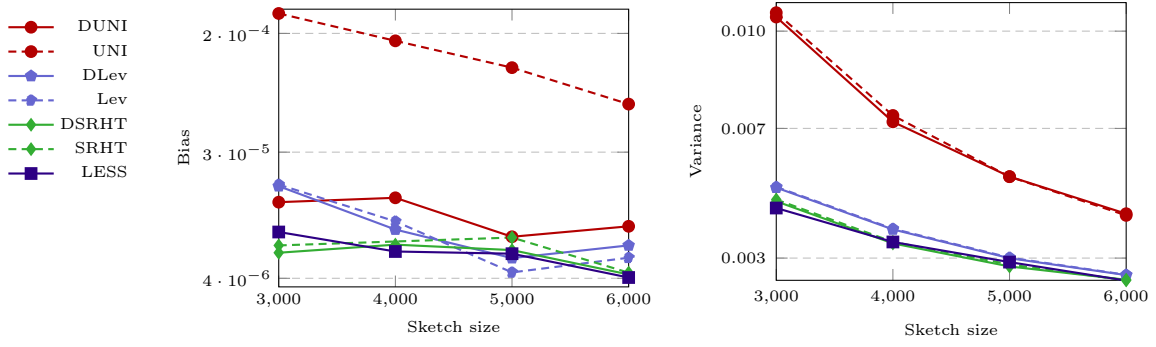
\begin{figure}[thb]
      \centering 
      \hspace*{-0.1\textwidth}
       \begin{subfigure}[c]{0.3\textwidth}
    \begin{tikzpicture}
    \renewcommand{\axisdefaulttryminticks}{5} 
    \pgfplotsset{every major grid/.style={densely dashed}}       
    \tikzstyle{every axis y label}+=[yshift=-10pt] 
    \tikzstyle{every axis x label}+=[yshift=5pt]
    \pgfplotsset{every axis legend/.append style={cells={anchor=east},fill=none,draw=none, at={(-0.8,1)}, anchor=north west, font=\tiny,legend columns=1,
        transpose legend }}
    \begin{axis}[
    width=1.3\columnwidth,
    height=1.1\columnwidth,
    xlabel style={font=\tiny},
        ylabel style={font=\tiny},
        tick label style={font=\tiny},
     xmin=3000,
      xmax=6000,
      ymin=3.487411e-06,
      ymax=2.944379e-04,
                    ymajorgrids=true,
                    scaled ticks=true,
                    xlabel = { Sketch size},
                    ylabel = { Bias},
                    ymode=log,
                    ytick={4e-6, 3e-5, 2e-4},
yticklabels={$4\cdot 10^{-6}$,$ 3\cdot 10^{-5}$,$2\cdot 10^{-4}$},
                    ]

     % revise 
     \addplot[mark=*,color=RED,line width=0.8pt,mark options={solid,fill=RED}] coordinates{
     (3000,1.349344e-05) (4000,1.446805e-05) (5000,7.764963e-06) (6000,9.175689e-06)
        };
        \addlegendentry{{DUNI}}[ font=\tiny];
    
              \addplot[densely dashed, mark=*,color=RED,line width=0.8pt,mark options={solid,fill=RED}] coordinates{
     (3000,2.750008e-04) (4000,1.775807e-04) (5000,1.157154e-04) (6000,6.464734e-05)
        };
        \addlegendentry{{UNI}}[ font=\tiny];

    % revise 
     \addplot[mark=pentagon*,color=BLUE!60!white,line width=0.8pt] coordinates{
      (3000,1.742441e-05) (4000,8.714470e-06) (5000,5.521574e-06) (6000,6.764934e-06)
          };
         \addlegendentry{{DLev}}[ font=\tiny];
                  
                  \addplot[densely dashed, mark=pentagon*,color=BLUE!60!white,line width=0.8pt] coordinates{
              (3000,1.787586e-05) (4000,9.924065e-06) (5000,4.393902e-06) (6000,5.556616e-06)
                    };
         \addlegendentry{{Lev}}[ font=\tiny];

    \addplot[mark=diamond*,color=GREEN!80!white,line width=0.8pt,mark options={solid,fill=GREEN!80!white}] coordinates{
  (3000,6.016039e-06) (4000,6.846337e-06) (5000,6.274772e-06) (6000,4.251345e-06)
    };
     \addlegendentry{{DSRHT}}[ font=\tiny];
     
    \addplot[densely dashed,mark=diamond*,color=GREEN!80!white,line width=0.8pt,mark options={solid,fill=GREEN!80!white}] coordinates{
     (3000,6.761861e-06) (4000) (6.566273e-06) (5000,7.670342e-06) (6000,4.375933e-06)
    };
     \addlegendentry{{SRHT}}[ font=\tiny];

      \addplot[mark=square*,color=RED!25!BLUE,line width=0.8pt] coordinates{
                (3000,8.378073e-06) (4000,6.136607e-06) (5000,5.923963e-06) (6000,4.049458e-06)
                    };
           \addlegendentry{{LESS}}[ font=\tiny];

    \end{axis}
    \end{tikzpicture}
    % \captionsetup{font=scriptsize}
    % \caption{MSD}
    % \label{subfig:SLev_and_SRHT}
    \end{subfigure}
    \hspace{0.25\textwidth}
      \begin{subfigure}[c]{0.3\textwidth}
      \begin{tikzpicture}
    \renewcommand{\axisdefaulttryminticks}{4} 
    \pgfplotsset{every major grid/.style={densely dashed}}       
    \tikzstyle{every axis y label}+=[yshift=-10pt] 
    \tikzstyle{every axis x label}+=[yshift=5pt]
    \pgfplotsset{every axis legend/.append style={cells={anchor=east},fill=none,draw=none, at={(1,0.95)}, anchor=north east, font=\tiny},legend columns=1,transpose legend}
    \begin{axis}[
    width=1.3\columnwidth,
    height=1.1\columnwidth,
    xlabel style={font=\tiny},
        ylabel style={font=\tiny},
        tick label style={font=\tiny},
     xmin=3000,
      xmax=6000,
      ymin=0.002314,
      ymax=0.010880,
      ymajorgrids=true,
      % scaled ticks=true,
       scaled y ticks=false,
                ytick={0.003,0.007,0.010},
yticklabels={$0.003$,$0.007$,$0.010$},
%          ytick={0.0840,0.0848,0.0856,0.0864},
% yticklabels={0.0840,0.0848,0.0856,0.0864},
% yticklabel style={
%     /pgf/number format/fixed,
%     /pgf/number format/precision=2
% },
      xlabel = { Sketch size},
      ylabel = { Variance},
      % ymode=log
      ]

              \addplot[mark=pentagon*,color=BLUE!60!white,line width=0.8pt] coordinates{
  (3000,0.005173) (4000,0.003884) (5000,0.002996) (6000,0.002475)
          };
         % \addlegendentry{{DLev}}[ font=\tiny];
         
                     \addplot[densely dashed, mark=pentagon*,color=BLUE!60!white,line width=0.8pt] coordinates{
              (3000,0.005195) (4000,0.003905) (5000,0.003015) (6000,0.002483)
                    };
         % \addlegendentry{{Lev}}[ font=\tiny];

         \addplot[mark=*,color=RED,line width=0.8pt,mark options={solid,fill=RED}] coordinates{
 (3000,0.010435) (4000,0.007203) (5000,0.005514) (6000,0.004368)
        };
        % \addlegendentry{{DUNI}}[ font=\tiny];
    
               \addplot[densely dashed, mark=*,color=RED,line width=0.8pt,mark options={solid,fill=RED}] coordinates{
     (3000,0.010569) (4000,0.007399) (5000,0.005510) (6000,0.004317)
        };
        % \addlegendentry{{UNI}}[ font=\tiny];

    \addplot[mark=diamond*,color=GREEN!80!white,line width=0.8pt,mark options={solid,fill=GREEN!80!white}] coordinates{
  (3000,0.004751) (4000,0.003453) (5000,0.002742) (6000,0.002331)
    };
     % \addlegendentry{{DSRHT}}[ font=\tiny];

    \addplot[densely dashed,mark=diamond*,color=GREEN!80!white,line width=0.8pt,mark options={solid,fill=GREEN!80!white}] coordinates{
 (3000,0.004806) (4000,0.003499) (5000,0.002803) (6000,0.002314)
    };
      % \addlegendentry{{SRHT}}[ font=\tiny];

                \addplot[mark=square*,color=RED!25!BLUE,line width=0.8pt] coordinates{
               (3000,0.004544) (4000,0.003497) (5000,0.002876) (6000,0.002305)
                    };
           % \addlegendentry{{LESS}}[ font=\tiny];

    \end{axis}
    \end{tikzpicture}
    % \captionsetup{font=scriptsize}
    % \caption{MSD}
    % \label{subfig:Uni}
    \end{subfigure}
    \captionsetup{skip=3pt}
\caption{{Bias and variance as functions of the sketch size $m$, comparing debiased sampling (in \textbf{solid} lines) and standard sampling (in \textbf{ dashed }lines) 
% under  uniform sampling (DUni/Uni), shrinkage leverage score (DSLev/SLev) sampling,  and SRHT (DSRHT/SRHT); and   LESS, 
on Airline dataset. \textbf{DUNI}, \textbf{DLev}, \textbf{DSRHT} are the corresponding debiased versions.  Expectation are estimated from $500$ independent runs. 
}}
% approximate leverage score (DLev/Lev),
\label{fig:error_ols_airline}
\end{figure}

Following \Cref{subsec:num_exper_ols}, \Cref{fig:error_ols_airline} examines how the sketch size $m$ affects the bias and variance of subsampled OLS on the Airline dataset. We report the bias by $(L(\hat\EE[\bbeta])-L(\bbeta_{\OLS}))/L(\bbeta_{\OLS})$  and the variance  by $(\hat\EE[L(\bbeta])]-L(\bbeta_{\OLS}))/L(\bbeta_{\OLS})$, where $L(\bbeta)=\|\y-\X\bbeta\|^2$ is  as in \Cref{eq:def_OLS_L}. We compare  uniform sampling (\textbf{UNI}), approximate leverage score sampling (\textbf{Lev}), \textbf{SRHT} (see \Cref{def:srht}), and LESS~\cite{derezinski2021newtonless,garg2024distributed}, together with their debiased variants \textbf{DUNI}, \textbf{DLev}, and \textbf{DSRHT}, using exact leverage scores.

The results on   the Airline dataset in \Cref{fig:error_ols_airline}  closely match those on the MSD dataset in \Cref{fig:error_ols}. UNI is the most biased,  while DUNI substantially reduces its  bias and,  approaches the comparable bias levels of Lev, DLev, SRHT, DSRHT, and LESS.
 Debiasing has  almost no visible effect on approximate leverage score sampling and SRHT, consistent with \Cref{coro:sub_ols_lev} and \Cref{coro:debiasing_srht}.  For variance, approximate leverage sampling, SRHT, their debiased variants, and LESS are comparable and uniformly outperform UNI and DUNI. Moreover, the debiased sampling schemes introduce essentially no variance increase relative to their corresponding standard counterparts.  These results further support that debiasing improves accuracy without sacrificing variance.  As expected, the standard non-debiased methods are also computationally cheaper; accordingly, \Cref{tab:ols_time_airline} reports timing results only for the standard schemes on the Airline and MSD datasets, where UNI, Lev, and SRHT are more efficient than LESS.

\begin{table}[htb]
  \caption{ Walk-clock time  versus sketch size $m$ for  various  methods on Airline and MSD datasets.  Results are obtained by averaging over $500$ independent runs. }
  \label{tab:ols_time_airline}
  \begin{center}
    \begin{small}
      %\begin{sc}
        \begin{tabular}{lccc}
          \toprule
      Dataset &     Method  &$m=3000$        & $m=5000$   \\
          \midrule
        \multirow{4}{*}{\textbf{Airline}}&   \textbf{UNI} & 0.0014& 0.0026  \\
        & \textbf{ Lev} & 0.0024& 0.0032  \\
        % Lev   & 0.0534 &0.1324   \\
        & \textbf{ SRHT}     & 0.3767&  0.4463\\
        &   \textbf{ LESS}   &0.5085&  0.5127\\
\midrule
 Dataset &     Method  &$m=5000$        & $m=7000$   \\
          \midrule
 \multirow{4}{*}{\textbf{MSD}} &   \textbf{UNI} & 0.0129 &  0.0181\\
          &  \textbf{ Lev} & 0.0240 &  0.0272  \\
        % Lev   & 0.0534 &0.1324   \\
      &   \textbf{ SRHT }    & 0.5864 &  0.6939\\
       &   \textbf{  LESS }  &  0.9814 &  1.0050 \\
          \bottomrule
        \end{tabular}
    %  \end{sc}
    \end{small}
  \end{center}
  \vskip -0.1in
\end{table}

\subsection{Numerical results for fast CUR decomposition} \label{subsec:numerical_cur}
% \subsection{CUR Decomposition} \label{subsec:num_exper_CUR}

In this section, we present numerical experiments for fast  CUR decomposition to support our theoretical findings on both classical and debiased CUR estimators. We evaluate CUR approximations on matrices $\X$ constructed from the CIFAR-10 dataset~\cite{krizhevsky2009Learning} and the ImageNet dataset~\cite{deng2009imagenet}.

In \Cref{fig:error_cur_cifar_10} and \cref{fig:error_cur_imagenet}, we study the effect of   sketch size $m$ on the bias and variance of fast CUR decomposition on CIFAR-10 and ImageNet datasets. We measure the bias by
$\|\hat\EE[\X-\C\U\R]\|_F^2/\|\X\|_F^2$ and the variance by $\hat\EE[\|\X-\C\U\R\|_F^2]/\|\X\|_F^2$.
We compare uniform sampling (\textbf{UNI}), approximate leverage score sampling (\textbf{Lev}),  and \textbf{SRHT} (see \Cref{def:srht}), together with their debiased counterparts \textbf{DUNI}, \textbf{DLev}, and \textbf{DSRHT}, using exact leverage scores. 
For each method, we first sample columns and rows of $\X$ to form $\C\in\RR^{n\times c}$ and $\R\in\RR^{r\times p}$, using column-norm and row-norm sampling probabilities, respectively. We fix $c=30$ and $r=60$ for CIFAR-10, and $c=40$ and $r=100$ for ImageNet. We then compute the core matrix $\U$ in \eqref{eq:def_U_fast} and \eqref{eq:def_check_U_fast}, with the standard sampling matrices
$\S_{\C}\in\RR^{m_c\times n}$ and $\S_{\R}\in\RR^{m_r\times p}$, and the debiased sampling matrices
$\check{\S}_{\C}\in\RR^{m_c\times n}$ and $\check{\S}_{\R}\in\RR^{m_r\times p}$.
\Cref{tab:cur_error_time} uses the same choices of  $c$  and  $r$  as in \Cref{fig:error_cur_cifar_10} and \Cref{fig:error_cur_imagenet}. 

As shown in \Cref{fig:error_cur_cifar_10} and \Cref{fig:error_cur_imagenet},  Lev and SRHT, together with their debiased variants, consistently outperform UNI and DUNI in the bias--variance performance, with SRHT achieving the best overall variance both with and without debiasing. Approximate leverage score sampling exhibits nearly identical bias and variance with and without debiasing; the same behavior is observed for SRHT. This is consistent with  \Cref{coro:sub_cur_lev} and \Cref{coro:debiasing_cur_srht}.
% Approximate leverage score sampling and SRHT exhibit nearly identical bias and variance with and without debiasing, consistent with \Cref{coro:sub_cur_lev} and \Cref{coro:debiasing_cur_srht}. 
In contrast, DUNI consistently reduces the bias of UNI, moving its performance closer to that of Lev, DLev, SRHT, and DSRHT, while maintaining nearly the same variance. Overall, debiasing improves accuracy without degrading variance.

In \Cref{tab:cur_error_time}, we compare the accuracy--efficiency trade-off of our method with the OSP methods of~\cite[Algorithm~2]{park2025accuracy} on CIFAR-10 and ImageNet datasets, where accuracy is measured by $\hat \EE[\|\X-\C\U\R\|_F^2/\|\X\|_F^2]$, denoted as variance in \Cref{fig:error_cur_cifar_10} and \Cref{fig:error_cur_imagenet}. The reported wall-clock time includes the computation of 
 $\C$, $\R$, $\tilde \U$ in \eqref{eq:faster_U} and $\C\tilde \U\R$ for Lev and SRHT, and the computation of  $\C$, $\R$, $\U^*=\X_{R,C}^{\dagger}$ and $\C\U^*\R$ for OSP-SRHT and OSP-SS. These OSP methods use randomized column-pivoted QR with oversampling to improve the stability and accuracy of the core matrix estimator $\U^\star=\X_{R,C}^{\dagger}$. Specifically, for OSP, we implement randomized column-pivoted QR using SRHT and sparse sign sketches, denoted by \textbf{OSP-SRHT} and \textbf{OSP-SS}, respectively. Since leverage-based methods and SRHT exhibit essentially identical performance with and without debiasing, while the standard versions are computationally cheaper, \Cref{tab:cur_error_time}  reports only the standard methods (Lev and SRHT) alongside OSP-SRHT and OSP-SS. 
 The results  show that Lev achieves a better accuracy--efficiency trade-off than both OSP baselines,  indicating that our CUR construction effectively balances computational cost and numerical stability. Although SRHT is slower than OSP-SS because SRHT-based sketching is more expensive than sparse sign sketching, it is still more accurate. Moreover, SRHT consistently outperforms OSP-SRHT in terms of the overall accuracy--efficiency trade-off, further demonstrating the advantage of our approach.

    {
    % \vspace{-0.7em}
    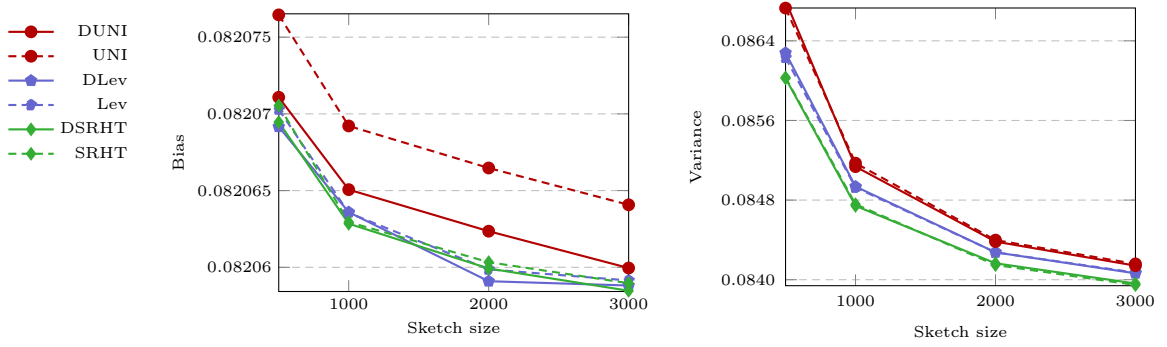
\begin{figure}[thb]
      \centering 
      \hspace*{-0.1\textwidth}
       \begin{subfigure}[c]{0.3\textwidth}
    \begin{tikzpicture}
    \renewcommand{\axisdefaulttryminticks}{5} 
    \pgfplotsset{every major grid/.style={densely dashed}}       
    \tikzstyle{every axis y label}+=[yshift=-10pt] 
    \tikzstyle{every axis x label}+=[yshift=5pt]
    \pgfplotsset{every axis legend/.append style={cells={anchor=east},fill=none,draw=none, at={(-0.8,1)}, anchor=north west, font=\tiny,legend columns=1,
        transpose legend }}
    \begin{axis}[
    width=1.3\columnwidth,
    height=1.1\columnwidth,
    xlabel style={font=\tiny},
        ylabel style={font=\tiny},
        tick label style={font=\tiny},
     xmin=500,
      xmax=3000,
      ymin=0.0820584,
      ymax=0.0820765,
                    ymajorgrids=true,
                    % scaled ticks=true,
                    scaled y ticks=false,
yticklabel style={
    /pgf/number format/fixed,
    /pgf/number format/precision=6
},
xtick={1000,2000,3000},
xticklabels={1000,2000,3000},
                    xlabel = { Sketch size},
                    ylabel = { Bias},
                    % ymode=log
                    ]

    \addplot[mark=*,color=RED,line width=0.8pt,mark options={solid,fill=RED}] coordinates{
         (500,0.08207106)  (1000,0.08206504)  (2000,0.08206233)  (3000,0.08205993
)
        };
        \addlegendentry{{DUNI}}[ font=\tiny];
    
              \addplot[densely dashed, mark=*,color=RED,line width=0.8pt,mark options={solid,fill=RED}] coordinates{
       (500,0.08207643)  (1000,0.08206919)  (2000,0.08206645)  (3000,0.08206405
)
        };
        \addlegendentry{{UNI}}[ font=\tiny];

      \addplot[mark=pentagon*,color=BLUE!60!white,line width=0.8pt] coordinates{
        (500,0.08206914)  (1000,0.08206357)  (2000,0.08205907)  (3000,0.08205879
)
          };
         \addlegendentry{{DLev}}[ font=\tiny];
                  
                  \addplot[densely dashed, mark=pentagon*,color=BLUE!60!white,line width=0.8pt] coordinates{
              (500,0.08207027)  (1000,0.08206353)  (2000,0.08205984)  (3000,0.08205912
)
                    };
         \addlegendentry{{Lev}}[ font=\tiny];

    \addplot[mark=diamond*,color=GREEN!80!white,line width=0.8pt,mark options={solid,fill=GREEN!80!white}] coordinates{
   (500,0.08206944)  (1000,0.08206282)  (2000,0.08205990)  (3000,0.08205845
)
    };
     \addlegendentry{{DSRHT}}[ font=\tiny];
     
    \addplot[densely dashed,mark=diamond*,color=GREEN!80!white,line width=0.8pt,mark options={solid,fill=GREEN!80!white}] coordinates{
     (500,0.08207052)  (1000,0.08206293)  (2000,0.08206031)  (3000,0.08205894
)
    };
     \addlegendentry{{SRHT}}[ font=\tiny];

    \end{axis}
    \end{tikzpicture}
    % \captionsetup{font=scriptsize}
    % \caption{MSD}
    % \label{subfig:SLev_and_SRHT}
    \end{subfigure}
    \hspace{0.25\textwidth}
      \begin{subfigure}[c]{0.3\textwidth}
      \begin{tikzpicture}
    \renewcommand{\axisdefaulttryminticks}{4} 
    \pgfplotsset{every major grid/.style={densely dashed}}       
    \tikzstyle{every axis y label}+=[yshift=-10pt] 
    \tikzstyle{every axis x label}+=[yshift=5pt]
    \pgfplotsset{every axis legend/.append style={cells={anchor=east},fill=none,draw=none, at={(1,0.95)}, anchor=north east, font=\tiny},legend columns=1,transpose legend}
    \begin{axis}[
    width=1.3\columnwidth,
    height=1.1\columnwidth,
    xlabel style={font=\tiny},
        ylabel style={font=\tiny},
        tick label style={font=\tiny},
     xmin=500,
      xmax=3000,
      ymin=0.08394,
      ymax=0.08673,
      ymajorgrids=true,
      % scaled ticks=true,
         scaled y ticks=false,
         ytick={0.0840,0.0848,0.0856,0.0864},
yticklabels={0.0840,0.0848,0.0856,0.0864},
% yticklabel style={
%     /pgf/number format/fixed,
%     /pgf/number format/precision=4
% },
xtick={1000,2000,3000},
xticklabels={1000,2000,3000},
      xlabel = { Sketch size},
      ylabel = { Variance},
      % ymode=log
      ]

     \addplot[mark=pentagon*,color=BLUE!60!white,line width=0.8pt] coordinates{
  (500,0.086279)  (1000,0.084936)  (2000,0.084275)  (3000,0.084064 )
          };
         % \addlegendentry{{DLev}}[ font=\tiny];
         
                     \addplot[densely dashed, mark=pentagon*,color=BLUE!60!white,line width=0.8pt] coordinates{
               (500,0.086231)  (1000,0.084926)  (2000,0.084275)  (3000,0.084071 )
                    };
         % \addlegendentry{{Lev}}[ font=\tiny];

         \addplot[mark=*,color=RED,line width=0.8pt,mark options={solid,fill=RED}] coordinates{
      (500,0.086802)  (1000,0.085136)  (2000,0.084383)  (3000,0.084141 )
        };
        % \addlegendentry{{DUNI}}[ font=\tiny];
    
               \addplot[densely dashed, mark=*,color=RED,line width=0.8pt,mark options={solid,fill=RED}] coordinates{
       (500,0.086729)  (1000,0.085171)  (2000,0.084398)  (3000,0.084159 )
        };
        % \addlegendentry{{UNI}}[ font=\tiny];

    \addplot[mark=diamond*,color=GREEN!80!white,line width=0.8pt,mark options={solid,fill=GREEN!80!white}] coordinates{
  (500,0.086024)  (1000,0.084742)  (2000,0.084165)  (3000,0.083962 )
    };
     % \addlegendentry{{DSRHT}}[ font=\tiny];

    \addplot[densely dashed,mark=diamond*,color=GREEN!80!white,line width=0.8pt,mark options={solid,fill=GREEN!80!white}] coordinates{
  (500,0.086032)  (1000,0.084758)  (2000,0.084148)  (3000,0.083947 )
    };
      % \addlegendentry{{SRHT}}[ font=\tiny];

    \end{axis}
    \end{tikzpicture}
    
    % \label{subfig:Uni}
    \end{subfigure}
    \captionsetup{skip=3pt}
   \caption{{Bias and variance as functions of the sketch size $m$, comparing debiased sampling (in \textbf{solid} lines) and standard sampling (in \textbf{dashed} lines) 
% under  uniform sampling (DUni/Uni), shrinkage leverage score (DSLev/SLev) sampling,  and SRHT (DSRHT/SRHT); and   LESS, 
on CIFAR-10  dataset. \textbf{DUNI}, \textbf{DLev}, \textbf{DSRHT} are the corresponding debiased versions.   Expectation are estimated from $200$ independent runs. 
% Results are obtained by averaging over $200$ independent runs. 
}}
% approximate leverage score (DLev/Lev),
\label{fig:error_cur_cifar_10}
    \end{figure}
    % \vspace{-1em}
    }

    {
    % \vspace{-0.7em}
    \begin{figure}[thb]
      \centering 
      \hspace*{-0.1\textwidth}
       \begin{subfigure}[c]{0.3\textwidth}
    \begin{tikzpicture}
    \renewcommand{\axisdefaulttryminticks}{5} 
    \pgfplotsset{every major grid/.style={densely dashed}}       
    \tikzstyle{every axis y label}+=[yshift=-10pt] 
    \tikzstyle{every axis x label}+=[yshift=5pt]
    \pgfplotsset{every axis legend/.append style={cells={anchor=east},fill=none,draw=none, at={(-0.8,1)}, anchor=north west, font=\tiny,legend columns=1,
        transpose legend }}
    \begin{axis}[
    width=1.3\columnwidth,
    height=1.1\columnwidth,
    xlabel style={font=\tiny},
        ylabel style={font=\tiny},
        tick label style={font=\tiny},
     xmin=1000,
      xmax=4000,
      ymin=0.11843809,
      ymax=0.11845761,
                    ymajorgrids=true,
                    % scaled ticks=true,
                    scaled y ticks=false,
yticklabel style={
    /pgf/number format/fixed,
    /pgf/number format/precision=6
},
                    xlabel = { Sketch size},
                    ylabel = { Bias},
                    % ymode=log
                    ]

     \addplot[mark=*,color=RED,line width=0.8pt,mark options={solid,fill=RED}] coordinates{
     (1000,0.11845242) (2000,0.11844463) (3000,0.11844302) (4000,0.11843919)
        };
        \addlegendentry{{DUNI}}[ font=\tiny];
    
              \addplot[densely dashed, mark=*,color=RED,line width=0.8pt,mark options={solid,fill=RED}] coordinates{
     (1000,0.11845761) (2000,0.11845178) (3000,0.11844821) (4000,0.11844564)
        };
        \addlegendentry{{UNI}}[ font=\tiny];

      \addplot[mark=pentagon*,color=BLUE!60!white,line width=0.8pt] coordinates{
       (1000,0.11844889)  (2000,0.11844255)  (3000,0.11844040)  (4000,0.11843830)
          };
         \addlegendentry{{DLev}}[ font=\tiny];
                  
                  \addplot[densely dashed, mark=pentagon*,color=BLUE!60!white,line width=0.8pt] coordinates{
              (1000,0.11845099)  (2000,0.11844253)  (3000,0.11844094)  (4000,0.11843835)
                    };
         \addlegendentry{{Lev}}[ font=\tiny];

    \addplot[mark=diamond*,color=GREEN!80!white,line width=0.8pt,mark options={solid,fill=GREEN!80!white}] coordinates{
 (1000,0.11844918)  (2000,0.11844181)  (3000,0.11843977)  (4000,0.11843809)
    };
     \addlegendentry{{DSRHT}}[ font=\tiny];
     
    \addplot[densely dashed,mark=diamond*,color=GREEN!80!white,line width=0.8pt,mark options={solid,fill=GREEN!80!white}] coordinates{
    (1000,0.11844977)  (2000,0.11844169)  (3000,0.11844005)  (4000,0.11843832)
    };
     \addlegendentry{{SRHT}}[ font=\tiny];

    \end{axis}
    \end{tikzpicture}
    % \captionsetup{font=scriptsize}
    % \caption{MSD}
    % \label{subfig:SLev_and_SRHT}
    \end{subfigure}
    \hspace{0.25\textwidth}
      \begin{subfigure}[c]{0.3\textwidth}
      \begin{tikzpicture}
    \renewcommand{\axisdefaulttryminticks}{4} 
    \pgfplotsset{every major grid/.style={densely dashed}}       
    \tikzstyle{every axis y label}+=[yshift=-10pt] 
    \tikzstyle{every axis x label}+=[yshift=5pt]
    \pgfplotsset{every axis legend/.append style={cells={anchor=east},fill=none,draw=none, at={(1,0.95)}, anchor=north east, font=\tiny},legend columns=1,transpose legend}
    \begin{axis}[
    width=1.3\columnwidth,
    height=1.1\columnwidth,
    xlabel style={font=\tiny},
        ylabel style={font=\tiny},
        tick label style={font=\tiny},
     xmin=1000,
      xmax=4000,
      ymin=0.122264,
      ymax=0.124952,
      ymajorgrids=true,
      % scaled ticks=true,
         scaled y ticks=false,
 ytick={0.1223,0.1230,0.1240,0.1249},
yticklabels={$0.1224$,$0.1230$,$0.1240$,$0.1247$},
% yticklabel style={
%     /pgf/number format/fixed,
%     /pgf/number format/precision=4
% },
% xtick={1000,2000,3000},
% xticklabels={1000,2000,3000},
      xlabel = { Sketch size},
      ylabel = { Variance},
      % ymode=log
      ]

  \addplot[mark=pentagon*,color=BLUE!60!white,line width=0.8pt] coordinates{
  (1000,0.124351) (2000,0.123024) (3000,0.122619) (4000,0.122401)
          };
         % \addlegendentry{{DLev}}[ font=\tiny];
         
                     \addplot[densely dashed, mark=pentagon*,color=BLUE!60!white,line width=0.8pt] coordinates{
              (1000,0.124378) (2000,0.123036) (3000,0.122599) (4000,0.122385)
                    };
         % \addlegendentry{{Lev}}[ font=\tiny];

         \addplot[mark=*,color=RED,line width=0.8pt,mark options={solid,fill=RED}] coordinates{
     (1000,0.125001) (2000,0.123358) (3000,0.122834) (4000,0.122571)
        };
        % \addlegendentry{{DUNI}}[ font=\tiny];
    
               \addplot[densely dashed, mark=*,color=RED,line width=0.8pt,mark options={solid,fill=RED}] coordinates{
     (1000,0.124952) (2000,0.123386) (3000,0.122820) (4000,0.122581)
        };
        % \addlegendentry{{UNI}}[ font=\tiny];

    \addplot[mark=diamond*,color=GREEN!80!white,line width=0.8pt,mark options={solid,fill=GREEN!80!white}] coordinates{
 (1000,0.124219) (2000,0.122916) (3000,0.122481) (4000,0.122264)
    };
     % \addlegendentry{{DSRHT}}[ font=\tiny];

    \addplot[densely dashed,mark=diamond*,color=GREEN!80!white,line width=0.8pt,mark options={solid,fill=GREEN!80!white}] coordinates{
 (1000,0.124214) (2000,0.122892) (3000,0.122483) (4000,0.122268)
    };
      % \addlegendentry{{SRHT}}[ font=\tiny];
    
    \end{axis}
    \end{tikzpicture}
    
    % \label{subfig:Uni}
    \end{subfigure}
\captionsetup{skip=3pt}
\caption{{Bias and variance as functions of the sketch size $m$, comparing debiased sampling (in \textbf{solid} lines) and standard sampling (in \textbf{dashed} lines) 
% under  uniform sampling (DUni/Uni), shrinkage leverage score (DSLev/SLev) sampling,  and SRHT (DSRHT/SRHT); and   LESS, 
on ImageNet dataset.  
\textbf{DUNI}, \textbf{DLev}, \textbf{DSRHT} are the corresponding debiased versions.   Expectation are estimated from $200$ independent runs. 
% Results are obtained by averaging over $200$ independent runs. 
}}
% approximate leverage score (DLev/Lev),
\label{fig:error_cur_imagenet}
\end{figure}
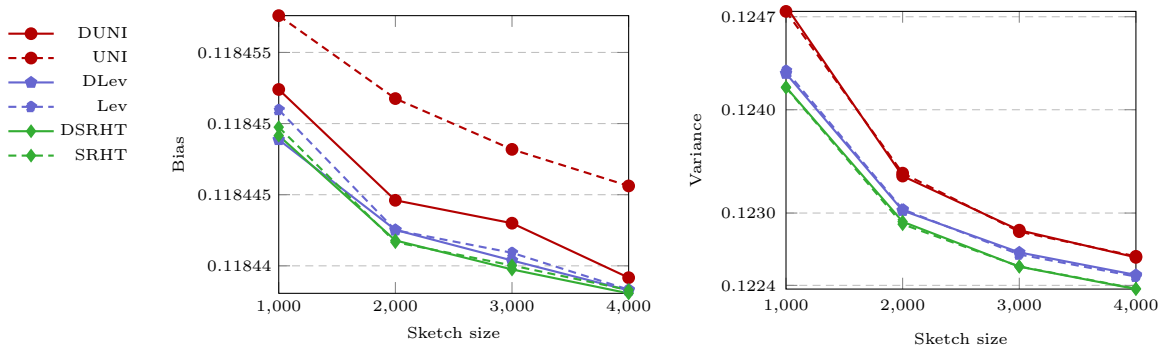

\begin{table}[htb]
  \caption{ Walk-clock time  versus relative error for Lev, SRHT, OSP-SRHT, and OSP-SS on CIFAR-10  and ImageNet datasets.   Expectation are estimated from $200$ independent runs. 
   % Results are averaged over 200 independent runs.  
   For \textbf{Lev} and \textbf{SRHT}, we vary  the sketch size $m_c$ of  $\S_c$ and fix the sketch size  $m_r=500$ of $\S_{\R}$, with $m_c\in\{2000,3000,4000\}$. %on CIFAR-10 and $m_c\in\{2000,3000,4000\}$ on ImageNet.
   For \textbf{OSP-SRHT} and \textbf{OSP-SS}, we vary the numbers of selected columns and rows of $\X$: on CIFAR-10, we use $c\in\{30,50,80\}$ and $r=c+r_0$ with row oversampling parameter $r_0=30$; on ImageNet, we use $c\in\{40,70,100\}$ and $r=c+r_0$ with $r_0=70$.
   } 
  \label{tab:cur_error_time}
  \begin{center}
    \begin{small}
      %\begin{sc}
        \begin{tabular}{lcccccccc}
          \toprule
Dataset & \multicolumn{2}{c}{\textbf{Lev}}  & \multicolumn{2}{c}{\textbf{SRHT} } & \multicolumn{2}{c}{\textbf{OSP-SRHT}} & \multicolumn{2}{c}{\textbf{OSP-SS}}\\
          \midrule
&Time(s)&  Error&  Time(s)&  Error &Time(s) &  Error  &Time(s) &  Error\\
 \midrule
 \multirowcell{4}{\textbf{CIFAR-10}}
 & 
0.0922&0.0843 &  0.6213&0.0842&0.7968&0.1042 & 0.1769&0.1044  \\
&0.1026&0.0841 & 0.7491&0.0840 & 0.9445&0.0946  & 0.1965& 0.0950\\
& 0.1120&0.0840 &  0.8088&0.0839 & 1.1439& 0.0885 &0.2390& 0.0888\\
\midrule

\multirowcell{4}{\textbf{ImageNet}}
& 
0.1869&0.1230 &     0.9936& 0.1229&  1.4659&0.1487&   0.2342&0.1491\\
& 0.1922& 0.1226 &   1.1565& 0.1225&   1.5233&0.1405 &  0.3434&0.1411\\
&0.2003&0.1224&   1.2704&0.1223&     1.5439&0.1398&  0.3845&0.1405 \\
       \bottomrule
        \end{tabular}
    %  \end{sc}
    \end{small}
  \end{center}
  \vskip -0.1in
\end{table}

\end{document}